\documentclass{stavanger}

\usepackage[utf8]{inputenc}
\usepackage{amssymb,amsmath,amsfonts,amsthm}
\usepackage{ dsfont }
\usepackage{hyperref}
\usepackage{natbib}
\usepackage{graphicx}
\usepackage{tikz}
\usepackage{cleveref}
\usepackage{mathtools}
\usepackage{placeins}
\usepackage{bm}
\usepackage{animate}
\usepackage[export]{adjustbox} 
\usepackage{mathbbol}

\startlocaldefs

\endlocaldefs

\begin{document}


\begin{frontmatter}

\begin{fmbox}


\dochead{Preprint}{FP}



\title{A symmetry and Noether charge preserving discretization of initial value problems}


\author[
   addressref={aff1},                   	  
   corref={aff1},                     		  
   email={alexander.rothkopf@uis.no}   		  
]{\inits{AR}\fnm{Alexander} \snm{Rothkopf}}
\author[
   addressref={aff2,aff3},                   	  
   email={jan.nordstrom@liu.se}   		  
]{\inits{JN}\fnm{Jan} \snm{Nordstr{\"o}m}}


\address[id=aff1]{
  \orgname{Faculty of Science and Technology}, 	 
  \street{University of Stavanger},                     		 
  \postcode{4021},                               			 
  \city{Stavanger},                              				 
  \cny{Norway}                                   				 
}

\address[id=aff2]{
  \orgname{Department of Mathematics}, 	 
  \street{Link{\"o}ping University},                     		 
  \postcode{SE-581 83},                               			 
  \city{Link{\"o}ping},                              				 
  \cny{Sweden}                                   				 
}

\address[id=aff3]{
  \orgname{Department of Mathematics and Applied Mathematics}, 	 
  \street{University of Johannesburg},                     		 
  \postcode{P.O. Box 524, Auckland Park 2006},                           
  \city{Johannesburg},                              				 
  \cny{South Africa}                                   				 
}



\end{fmbox}


\begin{abstractbox}
\begin{abstract} 
Taking insight from the theory of general relativity, where space and time are treated on the same footing, we develop a novel geometric variational discretization for second order initial value problems (IVPs). By discretizing the dynamics along a world-line parameter, instead of physical time directly, we retain manifest translation symmetry and conservation of the associated continuum Noether charge. A non-equidistant time discretization emerges dynamically, realizing a form of automatic adaptive mesh refinement (AMR), guided by the system symmetries. Using appropriately regularized summation by parts finite difference operators, the continuum Noether charge, defined via the Killing vector associated with translation symmetry, is shown to be exactly preserved in the interior of the simulated time interval. 
The convergence properties of the approach are demonstrated with two explicit examples.
\end{abstract}



\end{abstractbox}


\end{frontmatter}


\section{Introduction}
Symmetries play a central role in our understanding of dynamical processes in both classical \cite{goldstein1980classical,arnold2013mathematical} and quantum \cite{Coleman:1985rnk} physics. Emmy Noether achieved groundbreaking insight, when she proved that the presence of a global continuous symmetry in the action ${\cal S}$ of a system implies the existence of a conserved current, whenever the equations of motions are fulfilled \cite{noether1971invariant}. Via such a Noether current, one can define a quantity, which remains unchanged during the evolution of the system and which is referred to as Noether charge. Noether's theorem thus offers a fundamental understanding of central tenets of classical physics, such as energy and momentum conservation, which it relates to the invariance of physics under translations in time and space respectively.

In quantum theory, the presence of symmetries limits the type of quantum fluctuations which may occur \cite{Coleman:1985rnk}, with measurable consequences for the spectrum of elementary particles and their bound states. The four Noether currents associated with space and time translations are conventionally summarized in a quantity called the energy-momentum tensor $T^{\mu\nu}(x)$, where $\mu$ and $\nu$ refer to spatial and temporal components. It offers access to vital properties of a system, one pertinent example being the energy density profile \cite{landau2000classical} of a static charge distribution via the $\varepsilon(x)=T^{00}(x)$ component or the corresponding electric field-line configuration via the spatial components $T_{ij}(x)=F_{i\mu}F^{\mu}_j-\frac{1}{4}\delta_{ij}F_{\mu\nu}^2$ of the electromagnetic field $F_{\mu\nu}$, referred to as the Maxwell stress tensor (see e.g. \cite{Yanagihara:2018qqg}).

The simulation of dynamical phenomena in classical and quantum systems is often performed after discretizing space and time on a finite mesh (for a discussion of discretization in functional spaces see e.g \cite{cockburn2012discontinuous}). Finite difference schemes, formulated in their modern summation-by-parts (SBP) form (for reviews see e.g. \cite{svard2014review,fernandez2014review,lundquist2014sbp}) offer both conceptual and practical benefits. The SBP approach in both space and time \cite{lundquist2014sbp,nordstrom2013summation,nordstrom2016summation} offers proofs of stability based on the so-called energy method, which can be extended to high-order schemes in a straight forward fashion. Not only do SBP operators mimic integration by parts (IBP) exactly in the discretized setting, but in addition they constitute a cost effective approximation to differential operators on many mesh types.

The discretization of space and time in its conventional form, i.e. considering $x$ and $t$ as independent variables, necessarily affects the symmetry properties of the system at hand (see e.g. the discussion in \cite{JOHNSON1982147}). Where the continuum theory e.g. admits translations of any magnitude, i.e. in particular also infinitesimal ones, the discretized theory on a space-time mesh with grid spacing $\Delta_\mu$ only allows one to shift space and time by that finite amount. In general this entails that a central condition of Noether's theorem, the presence of a \textit{continuous} symmetry, does not hold and the corresponding continuum Noether charge fails to remain constant over time. This is particularly concerning with regards to time translation symmetry and energy conservation, which are closely related to the stability of the simulation. 

Artificial loss of energy is often considered benign, as it is simply a matter of loosing accuracy. An artificial increase of energy will, as energy is not bounded from above, eventually lead to a divergence of the simulated dynamics, characteristic of an unstable scheme. On the other hand, if energy is conserved, it puts stringent bounds on the growth of the solution. In the context of symplectic schemes, which conserve energy on average, one can relate energy conservation directly to the stability of the numerical scheme (see e.g. \cite{e6656c89-b0e4-3df4-8651-29dbb3d55273} and also \cite{nordstrom2023nonlinear}).

One strategy to retain energy conservation for systems with second order governing equations is to go over to a Hamiltonian approach, where only space is discretized, while time remains continuous. One converts the equation of motion of the Lagrange formalism, which is second order in the time derivative into a set of two equations of motion of first order, after replacing velocities with the so-called canonical momentum. After this step, a discrete phase-space volume preserving time stepping may be implemented (c.f. Verlet-Størmer \cite{verlet_computer_1967}). This approach crucially hinges on the availability of a Hamiltonian picture, i.e. whether the canonical momenta can be defined, which may face difficulties in systems with inherent constraints or requires the choice of a particular gauge, as in Maxwell's electrodynamics \cite{dirac_generalized_1950}.  Another strategy is to determine whether Noether's theorem may be salvaged in the presence of a finite grid spacing \cite{anerot_noethers-type_2020}. One may e.g. consider modifications to the continuum energy expression, which remain conserved, given a particular choice of difference approximation. However, as the necessary schemes are not of SBP type, they do not mimic other relevant properties of the continuum theory.  

In this study we develop a generic approach to discretize second order IVPs on the level of the system Lagrangian, while retaining the manifest translation invariance of the continuum theory. In order to do so we will take inspiration from the general theory of relativity (for a textbook see e.g. \cite{stephani2004relativity}), where \textit{space and time are treated on the same footing}. In this formalism the presence of translation symmetry is evident from the form of the Langrangian itself. We build upon our prior work on formulating IVPs directly via the action of the system, which allows us to avoid the need to derive their equation of motion. The action of the system is discretized using SBP finite difference operators with a physical null-space, developed in our previous paper \cite{Rothkopf:2022zfb}. These operators are crucial in mimicking the continuum derivation of Noether's theorem (and if one wishes to do so, the equations of motion).

The central outcome of this proof-of-principle study is a prescription of how to discretize second order IVPs directly on the level of the Lagrangian, while retaining the continuum time translation symmetry and thus exact conservation of the corresponding Noether charge. No reference to a Hamiltonian is required. We observe that a non-equidistant discretization emerges in the time coordinate, which represents a form of automatic adaptive mesh refinement (AMR) \cite{berger1984adaptive,lohner1987adaptive,berger1989local}, guided by the inherent symmetries of the system. Our results open up a novel route for obtaining optimal AMR procedures, where clustering and coarsening emerge as part of the solution process, thus avoiding the conventional use of sensors (see e.g. \cite{persson2006sub}), adjoint techniques (see e.g. \cite{nemec2008adjoint,offermans2023error}) or error estimates (see e.g. \cite{mavriplis1994adaptive,henderson1999adaptive,kompenhans2016adaptation}).

In \cref{sec:cont} we discuss the continuum formulation of our geometrized variational approach with time considered as dependent variable. In \cref{sec:disc} the discretized formalism is introduced and we present its efficacy in \cref{sec:num} using different example systems. We close with a summary and outlook in \cref{sec:conc}.

\section{Continuum formalism with manifest translation symmetry}
\label{sec:cont}

The common starting point for the formulation of the variational principle in classical point mechanics is to consider the dynamics of a system as boundary value problem (BVP). The system, which takes on position $x_i$ at $t_i$ evolves to position $x_f$ at $t_f$ and we wish to determine the trajectory it follows. Obviously this formulation is not causal, as we already need to know the end-point of the dynamics to determine the trajectory. As discussed in \cite{galley_classical_2013} and in our previous study \cite{Rothkopf:2022zfb} it is possible to formulate the variational problem as a genuine initial value problem through a doubling of the degrees of freedom of the system. 

In order to focus on the qualitatively novel ingredients of our variational approach, we first introduce it in the standard context of point mechanics as a BVP. The implementation for a genuine IVP is given in the subsequent subsection.

\subsection{Boundary value problem formulation}

Symmetry is a central mathematical pillar of the theory of relativity. In the special theory of relativity one formulates the laws of physics in a way that remains invariant under so-called Lorentz transformations of the coordinates, while in general relativity one constructs a description, which is invariant under an even larger class of transformations. Such a theory, invariant under arbitrary differentiable coordinate transformations, is called reparametrization invariant.

Reparameterization invariance is achieved by considering both space and time as dynamical degrees of freedom. In this study we are not interested in determining the dynamical evolution of space-time itself but will simply borrow this reparametrization invariant formalism of general relativity for our purposes of obtaining a symmetry preserving discretization. As our prime example, we set out to describe the dynamics of a point mass in the presence of a potential. The first step is to convert this physics question into a purely geometric problem.

In general relativity, the trajectory of a particle, traveling freely in (a not necessarily flat) space-time described by the metric tensor $g_{\mu\nu}$, is given by a path that generalizes the notion of the shortest path on the corresponding space-time manifold. This path is called a geodesic. While the particle may move in a $(1+d)$ dimensional space-time with $d$ space and one time direction, its path traces out a one-dimensional submanifold, which we can parameterize with a single, so called world-line parameter, denoted in the following by $\gamma$. We will restrict ourselves here to two dimensions, i.e. $d=1$, a system with one spatial and one temporal direction expressed in coordinates as ${\bf x}(\gamma)=(t(\gamma),x(\gamma))$.

A geodesic may be obtained from a variational principle \cite{jost1998calculus}, which asks for the critical point of the following action functional that measures the length of the path between two space-time points ${\bf x}(\gamma_i)$ and ${\bf x}(\gamma_f)$
\begin{align}
{\cal S}= \int_{\gamma_i}^{\gamma_f}\, d\gamma\, (-mc)\sqrt{ g_{\mu\nu} \frac{dx^{\mu}}{d\gamma} \frac{dx^{\nu}}{d\gamma}}, \qquad {\bf x}(\gamma_i)={\bf x}_i,\; {\bf x}(\gamma_f)={\bf x}_f.\label{eq:geodact}
\end{align}
Here Einstein's summation convention has been adopted and we have included the dimensionful prefactor $mc$, which, as we will show explicitly below, allows us to recover the usual action in the non-relativistic limit from \cref{eq:geodact}.

We refer to time $t(\gamma)$ as the zeroth component $x^0$ of the vector ${\bf x}$ and to the spatial coordinate $x(\gamma)$ as the first component $x^1$. Note that this functional is reparametrization invariant under any differentiable redefinition of the parameter $\gamma$. I.e. when converting from $\gamma\to\gamma^\prime$ the conversion of differentials under the square root produces terms $d\gamma^\prime/d\gamma$ that cancel with the conversion factor of the measure.

The geodesics of flat space-time, described by the diagonal metric tensor $g={\rm diag}[c^2,-1]$, which arise from the critical point of the action functional 
\begin{align}
    {\cal S}_{\rm flat}=\int_{\gamma_i}^{\gamma_f} d\gamma (-mc)\sqrt{c^2\Big(\frac{dt}{d\gamma}\Big)^2-\Big(\frac{dx}{d\gamma}\Big)^2},\label{eq:flataction}
\end{align}
are straight lines, which are traversed with constant speed (see chapter 3.4 of \cite{carroll2019spacetime}), in agreement with Newtonian mechanics. 

It is important to note that while our intuition of the concept of shortest path relies on geometries with positive definite metrics (Riemannian geometry), physical spacetime, as confirmed by experiment, has a metric with both positive and negative eigenvalues (pseudo-Riemannian geometry). In such a geometry the shortest path between two points can denote a saddle point of the action functional instead of a genuine minimum, as the temporal and spatial components enter relation (\ref{eq:geodact}) with opposite sign. 

To describe the presence of an external force acting on a point particle in flat spacetime, one conventionally amends the action ${\cal S}_{\rm flat}$ simply by adding the potential term $V(x)$ responsible for generating that force (see chapter 7.9 in \cite{goldstein1980classical}).

Let us now discuss how we can exploit the formalism of general relativity to re-express the evolution of a particle in flat spacetime in the presence of an external force, instead as an evolution of a free particle in a non-flat spacetime. In the presence of an external force, encoded in a potential term $V(x)$, the particle trajectory in flat space-time will deviate from the straight line. A standard procedure in the study of weak-field gravity is to reinterpret the change in the particle trajectory due to a potential, instead, as the effect of a non-flat space-time without a potential present (see e.g. chapter 8 of \cite{carlip2019general}). This reinterpretation is possible, as long as the values of the potential are smaller than the rest energy $(mc^2)$ of the point mass, a condition which is very well fulfilled for the non-relativistic systems we are interested in solving.

As we will see in the following, one can introduce the effects of a potential $V(x)$ on a point particle with mass $m$ in the weak-field limit of general relativity by modifying the temporal component $g_{00}$ of the diagonal metric tensor 
\begin{align}
    g_{00}=c^2+2V(x)/m, \label{eq:modmetr}
\end{align}
while keeping $g_{11}=-1$. I.e. one endows the metric with a non-trivial dependence on the spatial coordinate, trading the absence of an explicit external force for a non-flat spacetime.

Let us now show that such a modification of the metric indeed recovers the non-relativistic action of a particle in the presence of the potential $V(x)$. To this end we insert the modified metric \cref{eq:modmetr} into the geodesic action \cref{eq:geodact}:
\begin{align}
{\cal S}=& \int_{\gamma_i}^{\gamma_f}\,d\gamma\,(-mc)\sqrt{ g_{00} \Big(\frac{dt}{d\gamma}\Big)^2 - \Big(\frac{dx}{d\gamma}\Big)^2  }\label{eq:eq2}\\
\overset{g_{00}>0}{=}& \int_{\gamma_i}^{\gamma_f}\,d\gamma\,(-mc) \sqrt{ g_{00} \Big(\frac{dt}{d\gamma}\Big)^2 } \sqrt{ 1  - \frac{1}{g_{00}} \underbracket{\Big(\frac{dx}{d\gamma}\Big)^2 \Big(\frac{dt}{d\gamma}\Big)^{-2}}_{(dx/dt)^2} }\label{eq:eq3}\\
\nonumber \overset{\frac{dx}{dt}^2 \ll g_{00}\sim c^2}{=}& \int_{\gamma_i}^{\gamma_f}\,d\gamma\,\Big|\frac{dt}{d\gamma}\Big|\,(-mc) \sqrt{ g_{00} } \Big( 1 - \frac{1}{2} \frac{1}{g_{00}} \Big(\frac{dx}{d\gamma}\Big)^2 \Big(\frac{dt}{d\gamma}\Big)^{-2} \\
& \qquad\qquad\qquad\qquad \qquad\qquad\qquad\qquad\qquad + {\cal O}( \frac{1}{g_{00}^2}\Big(\frac{dx}{dt}\Big)^4 )\Big)\\
\nonumber \overset{V/m \ll c^2}{=}& \int_{\gamma_i}^{\gamma_f}\,d\gamma\, \frac{dt}{d\gamma} \Big( -mc^2   + \frac{1}{2} m \Big(\frac{dx}{d\gamma}\Big)^2 \Big(\frac{dt}{d\gamma}\Big)^{-2} - V(x)\\ 
& \qquad \qquad \qquad\qquad \qquad\qquad +{\cal O}(\Big(\frac{V}{mc^2}\Big)^2) + {\cal O}( \frac{1}{c^2}\Big(\frac{dx}{dt}\Big)^4 )  \Big) \label{eq:eq5}\\
=&\int_{t_i}^{t_f}\,dt \Big( -mc^2   + \frac{1}{2} m \Big(\frac{dx}{dt}\Big)^2 - V(x) \Big)\label{eq:eq6}.
\end{align}
In the third line we have expanded the rightmost square root in \cref{eq:eq3}, assuming that the square of the physical velocity $(dx/dt)^2$ is much smaller than $g_{00}$, which is to say that the particle velocity $dx/dt$ itself is much smaller than the speed of light $c$. To go from the third to the fourth line, we have in addition assumed that the potential is much smaller than the rest energy of the point particle, which allows us to expand the term $\sqrt{g_{00}}=\sqrt{c^2+2V(x)/m}$ in terms of $V(x)/mc^2$. We will look for solutions where time flows forward and thus have dropped the absolute value around $dt/d\gamma$ at the beginning of the second to last line. Note that \cref{eq:eq6} is nothing but the standard non-relativistic action \cite{goldstein1980classical} for a point particle in the presence of an arbitrary potential term with the rest energy $mc^2$ included.

We have thus successfully related the (artificially constructed) fully geometric description of the particle in a non-flat spacetime in \cref{eq:eq3} with the standard description of a particle propagating in flat spacetime in the presence of an external potential in \cref{eq:eq6} in the non-relativistic limit.

We see in \cref{eq:eq6} that time emerges naturally as the independent variable in which the action integral is formulated. Of course, choosing time as independent variable hides the inherent reparametrization invariance, which persists even in the non-relativistic limit in \cref{eq:eq5}. Interestingly it turns out that \cref{eq:eq5} is a generalization of the ad-hoc construction of a reparametrization invariant non-relativistic action, discussed in standard textbooks on the calculus of variations (see e.g. \cite{jost1998calculus}). \Cref{eq:eq5} includes the rest mass term $-mc^2 (dt/d\gamma)$, which is missing in the standard derivation and which in the absence of a potential contributes a dependence on $(dt/d\gamma)$ that plays a role in obtaining a well-defined critical point for the time degree of freedom.

The reward for our efforts lies in the fact that \cref{eq:eq2} is manifestly invariant under the space-time symmetries of our $(1+1)$ dimensional system. If $V(x)=0$ only the derivatives $dt/d\gamma$ and $dx/d\gamma$ but not $t$ and $x$ itself appear in the action functional \cref{eq:geodact}. In turn adding a constant shift to either $t$ or $x$ as in ${\bf x}\to {\bf x}+{\bf s}$ leaves the action invariant. In the presence of a spatially dependent potential $V(x)$, $g_{00}(x)$ too becomes dependent on space $x$ and only time translation invariance remains (as the force induced by $V(x)$ changes the momentum of the point particle). 

Proving time translation invariance in the conventional action \cref{eq:eq6} is much more involved, as one needs to consider how $x$ as a function of $t$ changes under such translations and in addition the boundaries of the action integral themselves are affected by the shift. None of these complications arise in \cref{eq:eq2}\footnote{That the derivatives of space and time occur in eq.~(\ref{eq:eq2}) as squares under the square root with a relative minus sign (hiding in $g_{11}$) also entails that the action is manifestly invariant under so called Lorentz boosts. These transformations mix space and time components and are related to changes between inertial coordinate systems.}.

In the calculus of variations it is known that the critical point of the action ${\cal S}$ can be obtained by solving certain differential equations, the so called geodesic equations \cite{jost1998calculus}. It follows from considering the variation of the action in all of its dependent variables $t$, $\dot t=dt/d\gamma$, $x$ and $\dot x=dx/d\gamma$
\begin{align}
\delta {\cal S}[t,\dot t, x,\dot x]=&\int_{\gamma_i}^{\gamma_f} d\gamma \Big\{ 
\frac{\partial {\cal L}}{\partial t}\delta t + \frac{\partial {\cal L}}{\partial \dot t}\delta \dot t + \frac{\partial {\cal L}}{\partial x}\delta x + \frac{\partial {\cal L}}{\partial \dot x}\delta \dot x \Big\}\label{eq:varS}\\
=&\int_{\gamma_i}^{\gamma_f} d\gamma \Big\{ \Big( \frac{\partial {\cal L}}{\partial t} - \frac{d}{d\gamma}\frac{\partial {\cal L}}{\partial \dot t}\Big)\delta t + \Big( \frac{\partial {\cal L}}{\partial x} - \frac{d}{d\gamma}\frac{\partial {\cal L}}{\partial \dot x}\Big)\delta x\Big\}\label{eq:geodesicEL}\\ &+ \left.\Big[ \frac{\partial {\cal L}}{\partial \dot t} \delta t \Big]\right|_{\gamma_i}^{\gamma_f}+ \left.\Big[ \frac{\partial {\cal L}}{\partial \dot x} \delta x \Big]\right|_{\gamma_i}^{\gamma_f}.
\end{align}
where in the second line we have integrated by parts. As we are considering the variational problem as boundary value problem with the coordinates $t$ and $x$ fixed at the start and end points of the trajectory ${\bf x}(\gamma_i)={\bf x}_i, {\bf x}(\gamma_f)={\bf x}_f$, also the variations $\delta t$ and $\delta x$ on the boundary vanish and so do the two boundary terms above. Note that we consider $t$ and $x$ as distinct degrees of freedom, so that the terms in the parentheses, multiplying the arbitrary variations $\delta x$ and $\delta t$, must vanish each independently at the stationary point $\delta {\cal S}=0$. 

By deriving the Euler-Lagrange equations of the system in the spirit of the standard BVP treatment of classical mechanics, the above derivation tells us that we may locate the classical trajectory of a non-relativistic particle under the influence of a potential, by finding the critical point of the action \cref{eq:geodact} with modified $g_{00}$ component of the metric, while keeping the start and end coordinates ${\bf x}(\gamma_i)$ and ${\bf x}(\gamma_f)$ fixed.

Note that there exist infinitely many different parameterizations of the trajectory described by $\delta {\cal S}=0$, which all differ by the velocity in $\gamma$, in which this trajectory is traversed. In practice these different stationary points of ${\cal S}$ lead to difficulties in numerical optimization 
and we therefore follow the standard practice (see e.g. discussion in \cite{rizzuti_square_2019} or \cite{carroll2019spacetime}) of selecting a particular parameterization by choosing instead of $\cal S$ the variations of the functional
\begin{align}
    {\cal E}_{\rm BVP}=\int_{\gamma_i}^{\gamma_f}d\gamma E_{\rm BVP}[t,\dot t,x,\dot x]=\int_{\gamma_i}^{\gamma_f}d\gamma \frac{1}{2}\Big( g_{00} \Big(\frac{dt}{d\gamma}\Big)^2 + g_{11}\Big(\frac{dx}{d\gamma}\Big)^2 \Big)\label{eq:actnri}.
\end{align}
It differs from $\cal S$ via squaring the integrand and replacing the pre-factor $-mc$ by $1/2$. These are both irrelevant changes with respect to the classical equation of motion. Since $\cal E_{\rm BVP}$ and ${\cal S}$ differ by a monotonous function applied to their integrands, formally the same critical point ensues. I.e. the variation of $\cal E_{\rm BVP}$ is given by $\delta {\cal L}=\delta \sqrt{E_{\rm BVP}}=\frac{\delta E_{\rm BVP}}{2\sqrt{E_{\rm BVP}}}=0$, so that the trajectory that extremizes ${\cal E}_{\rm BVP}$ agrees with that for ${\cal S}$ at the critical point. Note that the functional ${\cal E}_{\rm BVP}$ is not reparametrization invariant anymore. The derivative terms enter quadratically, and produce a conversion factor $(d\gamma^\prime/d\gamma)^2$, which cannot be absorbed by the measure $d\gamma$ alone. 

Let us compute the Euler-Lagrange equations (the geodesic equations) for time $t$ and space $x$ following from the variation of \cref{eq:actnri} 
\begin{align}
\nonumber&\delta {\cal E}_{\rm BVP}[t,\dot t, x,\dot x]\\
=&\int_{\gamma_i}^{\gamma_f} d\gamma \Big\{ 
\frac{\partial E_{\rm BVP}}{\partial t}\delta t + \frac{\partial E_{\rm BVP}}{\partial \dot t}\delta \dot t + \frac{\partial E_{\rm BVP}}{\partial x}\delta x + \frac{\partial E_{\rm BVP}}{\partial \dot x}\delta \dot x \Big\}\label{eq:varE}\\
=&\int_{\gamma_i}^{\gamma_f} d\gamma \Big\{ \Big( \frac{\partial E_{\rm BVP}}{\partial t} - \frac{d}{d\gamma}\frac{\partial E_{\rm BVP}}{\partial \dot t}\Big)\delta t + \Big( \frac{\partial E_{\rm BVP}}{\partial x} - \frac{d}{d\gamma}\frac{\partial E_{\rm BVP}}{\partial \dot x}\Big)\delta x\Big\}\label{eq:geodesicEE}\\ &+ \left.\Big[ \frac{\partial E_{\rm BVP}}{\partial \dot t} \delta t \Big]\right|_{\gamma_i}^{\gamma_f}+ \left.\Big[ \frac{\partial E_{\rm BVP}}{\partial \dot x} \delta x \Big]\right|_{\gamma_i}^{\gamma_f}.
\end{align}
As the above boundary terms vanish, we are left with evaluating the individual expressions appearing in the parentheses of \cref{eq:geodesicEE}. Below we evaluate each of these terms individually 
\begin{align}
    &\frac{\partial E_{\rm BVP}}{\partial t}=0, \quad \frac{\partial E_{\rm BVP}}{\partial \dot t}=g_{00}(x)\frac{dt}{d\gamma},\\
    &\frac{\partial E_{\rm BVP}}{\partial x}=\frac{1}{2} \frac{\partial g_{00}(x)}{\partial x}\Big(\frac{dt}{d\gamma}\Big)^2, \quad \frac{\partial E_{\rm BVP}}{\partial \dot x}=g_{11}\frac{dx}{d\gamma}=-\frac{dx}{d\gamma},
\end{align}
making explicit the ingredients to the geodesic equations for the temporal and spatial degrees of freedom
\begin{align}
 &\frac{d}{d\gamma}\Big(g_{00}\frac{dt}{d\gamma}\Big)=0, \label{eq:geodt}\\
 &\frac{d}{d\gamma}\Big(\frac{dx}{d\gamma}\Big) + \frac{1}{2}\frac{\partial g_{00}}{\partial x}\Big( \frac{dt}{d\gamma} \Big)^2=0\label{eq:geodx}.
\end{align}
The attentive reader will have recognized that \cref{eq:geodt} constitutes a conservation equation for the expression inside the parenthesis. In the next chapter we will show that this quantity indeed is the conserved charge associated with the time translation symmetry of our system. In general the geodesic equations do not single out the conserved quantities in such a simple fashion. There however exists an systematic procedure to identify the space-time symmetries of the system in the form of different so-called Killing vectors, each of which leads to one conserved quantity (see \cref{sec:consquant}).

Note that the geodesic equations \cref{eq:geodt,eq:geodx} are often written in a more concise fashion in the general relativity literature (see e.g. \cite{stephani2004relativity}). They are expressed for a general metric using the so-called Christoffel symbols $\Gamma^\alpha_{\mu\nu}=\frac{1}{2}g^{\alpha\beta}\big( \partial g_{\beta \mu}/\partial x_\nu + \partial g_{\beta \nu}/\partial x_\mu - \partial g_{\mu\nu}/\partial x_\beta)$, where $g^{\alpha\beta}$ refers to the components of the inverse of the metric $g_{\alpha\beta}$. One obtains in short hand notation with Einstein summation implied
\begin{align}
    \frac{d^2 x^\alpha}{d \gamma^2}+\Gamma^\alpha_{\mu\nu} \frac{d x^\mu}{d \gamma}\frac{d x^\nu}{d \gamma} = 0 \label{eq:geodeqconsc}.
\end{align}
It is important to note that the derivation of the above expression involves application of the product rule, which in the discrete setting is not valid. Therefore even though in the continuum  \cref{eq:geodt,eq:geodx} and \cref{eq:geodeqconsc} are equivalent, we will work solely with the former, as only integration by parts (which is exactly mimicked by summation by parts) has been used in their derivation.

\subsection{Conserved quantities, Noether's theorem and stability}
\label{sec:consquant}

Conservation of momentum and energy in general relativity is conceptually more involved compared to flat space-time, since the comparison of two quantities at different space-time points becomes a non-trivial operation due to the effects of a non-flat metric. However there may exist a vector field $K^\mu(x)$ along which transported quantities remain constant. These vector fields are known as Killing\footnote{For completeness we note that a Killing vector field $K_\mu$ is defined as solution to the Killing equation
$\Big(\frac{\partial K_\mu}{\partial x^\nu}-\Gamma^\alpha_{\mu\nu}K_\alpha\Big) + \Big(\frac{\partial K_\nu}{\partial x^\mu}-\Gamma^\alpha_{\nu\mu}K_\alpha\Big) =0.$} vector fields $K^\mu(x)$.
The Killing vector fields are generators of infinitesimal isometries of the space-time manifold. Moving all points of the manifold in the direction of the Killing field leaves the manifold unchanged.

As discussed in standard literature on general relativity (see e.g. chapter 3.8 of \cite{carroll2019spacetime}), each Killing vector field $K^\mu$ can be used to define a conserved quantity $Q_K$ via the expression
\begin{align}
    Q_K=g_{\alpha\beta}K^{\alpha}\dot x^{\beta}\label{eq:killingcons}.
\end{align}
Computing the change of $Q_K$ along a geodesic, parameterized by $\gamma$, one finds from combining \cref{eq:geodeqconsc,eq:killingcons} and the equation that defines the Killing vector that $dQ_K/d\gamma=0$, i.e. it vanishes. We will give an explicit example of such a conserved quantity below.

More intuitively, one can think of the role of $K^\mu$ as pointing out directions along which the metric $g$ of spacetime in our system remains constant. In the spirit of Noether's theorem, assume that the integrand $E_{\rm BVP}$ of our action functional ${\cal E}_{\rm BVP}$ in \cref{eq:actnri} remains unchanged under infinitesimal translations with magnitude $\epsilon$ in the direction of $K^\mu$. The change in coordinates under such a shift is $\delta x^\mu=\epsilon K^\mu$. Noether's theorem tells us that the conserved quantity corresponding to $\delta x^\mu$ is given by $J=\delta x^\mu \frac{\partial E}{\partial \dot x^{\mu}}$, which, when written explicitly as $\epsilon K^\alpha g_{\alpha\beta} \frac{dx^\beta}{d\gamma}$, turns out to just be $\epsilon Q_K$.

In case of our geometrized problem of determining the dynamics of a point particle under the influence of a potential $V(x)$, the metric remains independent of time $t$. Thus the vector $K_t=(1,0)$ constitutes a Killing vector associated with time translation symmetry. The conservation of the associated conserved quantity $Q_t=K_t^\mu g_{\mu\nu}\dot x^\nu= g_{00}\dot t$ follows straight forwardly from the geodesic equation for $t$
\begin{align}
   \frac{d}{d\gamma} Q_t \overset{K_t=(1,0)}{\underset{\cref{eq:killingcons}}{=}} \frac{d}{d\gamma} \Big( g_{00}\dot t \Big) \overset{\cref{eq:geodt}}{=}0 \label{eq:Qt},
\end{align}
i.e. the quantity $Q_t$ remains constant along the geodesic. Note that this quantity is different from the usual energy considered in the non-relativistic formalism.

Turning to the question of stability, let us show next that as a consequence of the presence of a conserved quantity together with the form of the geodesic equations and the reasonable assumption that the potential of the system is bounded from below, it is possible to provide an upper bound on the derivatives of the trajectories obtained as critical point of the functional \cref{eq:actnri}.

In an analogy to the construction of a Hamiltonian from a Lagrangian, we define the following
\begin{align}
    {\cal H}_{\rm BVP}&=\int_{\gamma_i}^{\gamma_f}d\gamma\underbracket{\frac{1}{2}\Big( g_{00}(x)\Big(\frac{dt}{d\gamma}\Big)^2 - g_{11}\Big(\frac{dx}{d\gamma}\Big)^2 \Big)}_{H_{\rm BVP}},\\
    &=\int_{\gamma_i}^{\gamma_f}d\gamma\frac{1}{2}\Big( \big(c^2+2\frac{V(x)}{m}\big)\Big(\frac{dt}{d\gamma}\Big)^2 + \Big(\frac{dx}{d\gamma}\Big)^2 \Big).
\end{align}
Due to the flipped sign in front of $g_{11}$, compared to the action \cref{eq:actnri}, this quantity is actually positive definite, as long as $V(x)$ is bounded from below\footnote{Since physical forces arise from the derivative of the potential, we may always add a constant to a bounded potential that will make $g_{00}$ positive.}. ${\cal H}_{\rm BVP}$ thus provides a norm on the function space in which $t(\gamma)$ and $x(\gamma)$ reside. Now let us inspect the evolution of the integrand $H_{\rm BVP}$ 
\begin{align}
    \nonumber \frac{d H_{\rm BVP}}{d\gamma}&=\frac{1}{2}\frac{dg_{00}}{d\gamma}\Big(\frac{dt}{d\gamma}\Big)^2+ g_{00}\frac{dt}{d\gamma} \frac{d^2t}{d\gamma^2}+\frac{dx}{d\gamma} \frac{d^2x}{d\gamma^2},\\
    \nonumber &= \frac{dx}{d\gamma}\Big[ \frac{1}{2}\frac{\partial g_{00}}{\partial x}\Big(\frac{dt}{d\gamma}\Big)^2 +\frac{d^2x}{d\gamma^2} \Big] + g_{00}\frac{dt}{d\gamma} \frac{d^2t}{d\gamma^2},\\
    &\overset{\cref{eq:geodx}}{=}\underbracket{g_{00}\frac{dt}{d\gamma}}_{Q_t\,\rm const.} \frac{d^2t}{d\gamma^2}\label{eq:evolHBVP}.
\end{align}
To arrive at the final expression in \cref{eq:evolHBVP}, we use the fact that one can rewrite $dg_{00}/d\gamma=(\partial g_{00}(x)/\partial x) \dot x$ and combine the first and third term to apply \cref{eq:geodx}. This simplification tells us that the change in $H_{\rm BVP}$ is given solely by the second derivative of time with respect to the world-line parameter. Now we can integrate up twice ${\cal H}_{\rm BVP}= \int_{\gamma_i}^{\gamma_f}d\gamma \int_{\gamma_i}^{\gamma}d\gamma^\prime (dH_{\rm BVP}(\gamma^\prime)/d\gamma^\prime)$ to get 
\begin{align}
    \nonumber {\cal H}_{\rm BVP}&= m g_{00}(x_i)\dot t(\gamma_i)\;\int_{\gamma_i}^{\gamma_f} d\gamma \big( \dot t(\gamma)-\dot t(\gamma_i)\big),\\
    \nonumber &= g_{00}(x_i)\dot t(\gamma_i)\; \Big( -\dot t(\gamma_i)(\gamma_f-\gamma_i) + ( t(\gamma_f)-t(\gamma_i) ) \Big),\\
    &\leq  g_{00}(x_i)\dot t(\gamma_i)\big( t(\gamma_f)-t(\gamma_i) \big)  \label{eq:boundonsol}.
\end{align}
For the last inequality we use the fact that the world-line is parameterized by an increasing $\gamma$ and correspondingly time moves forward along the world-line.

In the BVP setting, where both $t(\gamma_i)$ and $t(\gamma_f)$ are given apriori, \cref{eq:boundonsol} constitutes a proof that the norm ${\cal H}_{\rm BVP}$ defined on the derivatives of the solution $t$ and $x$ grows at most linearly with time, precluding the occurrence of exponentially increasing behavior that would signal an instability, in turn establishing stability of the geometric approach.

\subsection{Initial value formulation}
\label{sec:initvalform}
So far we have shown how the geodesic equations \cref{eq:geodt,eq:geodx} can be obtained from a variational principle formulated as a boundary value problem in time. However for a causal description as an initial value problem, we must be able to determine the dynamics of the particle without knowledge of the final point of the trajectory. If one wishes to prescribe only initial values, i.e. positions and derivatives at $\gamma_i$, then the variations $\delta x^\mu$ in \cref{eq:varS} do not vanish at the end of the particle world line, i.e. at $\gamma_f$. In turn the equivalence between the critical point of ${\cal S}$ and the Euler-Lagrange equations in \cref{eq:geodesicEL} does not hold. As discussed by \cite{galley_classical_2013} and put into practice in our previous publication \cite{Rothkopf:2022zfb} one can overcome this issue by constructing an action with doubled degrees of freedom, living on a closed contour with a forward and backward branch in $\gamma$.

Since both time and space constitute dependent degrees of freedom in our approach, we need to introduce both forward and backward variants of each of them $x_1(\gamma), x_2(\gamma)$ and $t_1(\gamma), t_2(\gamma)$. The degrees of freedom on the forward contour enter the action functional with the usual Lagrangian, while those on the backward contour are assigned the negative Lagrangian. Choosing to build the doubled formalism based on the action ${\cal E}_{\rm BVP}$ we obtain
\begin{align}
    {\cal E}_{\rm IVP}&=\int_{\gamma_i}^{\gamma_f} d\gamma \, E_{\rm IVP}[t_1,\dot t_1,x_1,\dot x_1,t_2,\dot t_2,x_2,\dot x_2],\\
    &=\int_{\gamma_i}^{\gamma_f} d\gamma \Big\{ E_{\rm BVP}[t_1,\dot t_1,x_1,\dot x_1] - E_{\rm BVP}[t_2,\dot t_2,x_2,\dot x_2] \Big\}\label{eq:EIVPcont}.
\end{align}
As discussed in detail in \cite{Rothkopf:2022zfb}, the inner workings of the doubled formalism become more transparent, once we go over to expressing the action $\cal E_{\rm IVP}$ in terms of the central and difference coordinates $x_+=\frac{1}{2}(x_1+x_2)$ and $x_-=x_1-x_2$ and $t_+=\frac{1}{2}(t_1+t_2)$ and $t_-=t_1-t_2$ respectively.  The variation now proceeds in the independent degrees of freedom $x_{\pm}$ and $t_{\pm}$ and yields
\begin{align}
\nonumber &\delta {\cal E}_{\rm IVP}[t_{\pm},\dot t_{\pm}, x_{\pm},\dot x_{\pm}]=\\
&\int_{\gamma_i}^{\gamma_f} d\gamma \Big\{ 
\frac{\partial E_{\rm IVP}}{\partial t_+}\delta t_+ + \frac{\partial E_{\rm IVP}}{\partial \dot t_+}\delta \dot t_+ + \frac{\partial E_{\rm IVP}}{\partial t_-}\delta t_- + \frac{\partial E_{\rm IVP}}{\partial \dot t_-}\delta \dot t_- \label{eq:varEIVP}\\
\nonumber &\qquad + \frac{\partial E_{\rm IVP}}{\partial x_+}\delta x_+ + \frac{\partial E_{\rm IVP}}{\partial \dot x_+}\delta \dot x_+ + \frac{\partial E_{\rm IVP}}{\partial x_-}\delta x_- + \frac{\partial E_{\rm IVP}}{\partial \dot x_-}\delta \dot x_- \Big\}\\
\nonumber =&\int_{\gamma_i}^{\gamma_f} d\gamma \Big\{ \Big( \frac{\partial E_{\rm IVP}}{\partial t_+} - \frac{d}{d\gamma}\frac{\partial E_{\rm IVP}}{\partial \dot t_+}\Big)\delta t_+ +\Big( \frac{\partial E_{\rm IVP}}{\partial t_-} - \frac{d}{d\gamma}\frac{\partial E_{\rm IVP}}{\partial \dot t_-}\Big)\delta t_-\\
&+ \Big( \frac{\partial E_{\rm IVP}}{\partial x_+} - \frac{d}{d\gamma}\frac{\partial E_{\rm IVP}}{\partial \dot x_+}\Big)\delta x_+ + \Big( \frac{\partial E_{\rm IVP}}{\partial x_-} - \frac{d}{d\gamma}\frac{\partial E_{\rm IVP}}{\partial \dot x_-}\Big)\delta x_-\Big\}\label{eq:geodesicELIVP}\\ \nonumber&+ \left.\Big[ \frac{\partial E_{\rm IVP}}{\partial \dot t_+} \delta t_+ \Big]\right|_{\gamma_i}^{\gamma_f} + \left.\Big[ \frac{\partial E_{\rm IVP}}{\partial \dot t_-} \delta t_- \Big]\right|_{\gamma_i}^{\gamma_f} + \left.\Big[ \frac{\partial E_{\rm IVP}}{\partial \dot x_+} \delta x_+ \Big]\right|_{\gamma_i}^{\gamma_f} + \left.\Big[ \frac{\partial E_{\rm IVP}}{\partial \dot x_-} \delta x_- \Big]\right|_{\gamma_i}^{\gamma_f}.
\end{align}
To arrive at \cref{eq:geodesicELIVP} we have carried out four integrations by parts. As the next step, we consider under which conditions the boundary terms in the above expression vanish. Since we prescribe fixed initial values for both time and space, the variations $\delta t_{\pm}(\gamma_i)=0$ and $\delta x_{\pm}(\gamma_i)=0$ vanish. What about the variations at the end of the forward and backward world-line? As long as we require that 
\begin{align}
    x_2(\gamma_f)=x_1(\gamma_f), \quad t_2(\gamma_f)=t_1(\gamma_f)\label{eq:conconI},
\end{align}
it follows that $\delta x_{-}(\gamma_f)$ and $\delta t_{-}(\gamma_f)$ vanish and with it the corresponding boundary terms. The only remaining terms are those at $\gamma_f$ which feature $\delta x_+$ and $\delta t_+$. As these variations do not vanish, we instead inspect the terms multiplying them, i.e. $\partial E_{\rm IVP}/\partial\dot t_+$ and $\partial E_{\rm IVP}/\partial\dot x_+$. Using the definition $x_1=x_+ + \frac{1}{2} x_-$ and $x_2=x_+ - \frac{1}{2} x_-$ and correspondingly for $t_{1,2}$, we find from the defining equation for $E_{\rm IVP}$ \cref{eq:actnri}
\begin{align}
    \nonumber\frac{d E_{\rm IVP}}{d \dot x_+} &=  \frac{\partial E_{\rm IVP}[t_{1,2},\dot t_{1,2},x_{1,2},\dot x_{1,2}]}{\partial \dot x_1}\frac{d \dot x_1}{d \dot x_+} + \frac{\partial E_{\rm IVP}[t_{1,2},\dot t_{1,2},x_{1,2},\dot x_{1,2}]}{\partial \dot x_2}\frac{d \dot x_2}{d\dot x_+},\\
    &= \frac{\partial E_{\rm BVP}[t_1,\dot t_1,x_1,\dot x_1]}{\partial \dot x_1}\frac{d \dot x_1}{d \dot x_+} - \frac{\partial E_{\rm BVP}[t_2,\dot t_2,x_2,\dot x_2]}{\partial \dot x_2}\frac{d \dot x_2}{d\dot x_+}, \\
    &= g_{11}(x_1)\dot x_1-g_{11}(x_2)\dot x_2=-\dot x_1+\dot x_2.
\end{align}
Similarly one obtains 
\begin{align}
    \frac{d E_{\rm IVP}}{d \dot t_+} =  g_{00}(x_1)\dot t_1-g_{00}(x_2)\dot t_2.
\end{align}
Together with condition \cref{eq:conconI} that the values of $x_{1,2}$ and $t_{1,2}$ must agree at $\gamma_f$, this result tells us that in order for the two remaining boundary terms to vanish, we need to also identify the derivatives of $x_{1,2}$ and $t_{1,2}$ at the point $\gamma_f$
\begin{align}
    \dot x_2(\gamma_f)=\dot x_1(\gamma_f), \quad \dot t_2(\gamma_f)=\dot t_1(\gamma_f)\label{eq:conconII}.
\end{align}

Note that we have now managed to remove the boundary terms without the need for specifying the concrete value of t's and x's at the final point $\gamma_f$. This is the central contribution of the forward-backward construction.

The last remaining step is to undo the proliferation of degrees of freedom that occurred when introducing the forward-backward construction. It has been shown \cite{Berges:2007ym,galley_classical_2013} that taking the so-called physical limit achieves this goal, where the constraints $x_1(\gamma)-x_2(\gamma)=x_-(\gamma)=0$ and $t_1(\gamma)-t_2(\gamma)=t_-(\gamma)=0$ are enforced. The remaining $x_+$ and $t_+$ are identified with the true classical geodesics.

In terms of the Euler-Lagrange equations in parentheses in \cref{eq:geodesicELIVP} 
\begin{align}
    \frac{\partial E_{\rm IVP}}{\partial x_\pm} - \frac{d}{d\gamma}\frac{\partial E_{\rm IVP}}{\partial \dot x_\pm}=0, \qquad \frac{\partial E_{\rm IVP}}{\partial t_\pm} - \frac{d}{d\gamma}\frac{\partial E_{\rm IVP}}{\partial \dot t_\pm}=0,
\end{align}
the physical limit entails that only those equations independent of $x_-$ and $t_-$ survive. With the construction of the action $E_{\rm IVP}= E_{\rm BVP}[x_1,\dot x_1,t_1,\dot t_1] - E_{\rm BVP}[x_2,\dot x_2,t_2,\dot t_2] $ from a difference of the $E_{\rm BVP}$ functionals, there will appear at least a linear dependence on the minus degrees of freedom. Hence in the physical limit only those Euler-Lagrange equations linear in $x_-$ and $t_-$ will survive, where the minus degrees of freedom have been removed by taking the derivative with respect to $x_-$ or $t_-$. 


Note that we have decided to not only specify the value and derivative of $x$ at initial $\gamma_i$ but also those of $t$. As we wish to determine the dynamics of a point particle in the presence of a potential with given $x(t_i)$ and $dx/dt(t_i)$, there remains a freedom in choosing $\dot x(\gamma_i)$ and $\dot t(\gamma_i)$, since only their ratio needs to be fixed $dx/dt(t=t_0)=\dot x(\gamma_i)/\dot t (\gamma_i)$. The end of the time interval traversed by the world line parameter $\gamma$, will consequently depend on the value prescribed to $\dot t(\gamma_i)$ and emerges dynamically from the combined evolution of $x$ and $t$.

At this point we have formulated a manifest time translation symmetric variational principle that encodes the dynamics of a point particle evolving in the presence of a non-relativistic potential as initial value problem. Our next goal is to discretize the action functional ${\cal E}_{\rm IVP}$ in \cref{sec:disc} using SBP finite difference operators. Since all derivations of the Euler-Lagrange equations, as well as that of the conserved quantity $Q_t$ have made ample reference to integration by parts, it is paramount to use such a discretization technique, which faithfully mimics this continuum property on a finite mesh.

\section{Discretized formalism for IVPs}
\label{sec:disc}

The central novelty we introduce in this section is related to the fact that the discretization of the action functional takes place in the world-line parameter $\gamma$ and not in the time variable $t$, as in conventional discretization prescriptions. I.e. the values of both time $t(\gamma)$ and position $x(\gamma)$ remain continuous and in turn we achieve \textit{preservation of the continuum space-time symmetries} even after discretization.

In the presence of a potential that depends on $x$ but not on $t$, the invariance under infinitesimal constant shifts in time is hence retained. This comes about, since the metric remains invariant under changes in $t$, which in turn leads to a simple form of the corresponding Killing equation, which shows that $K_t=(1,0)$ indeed is a Killing vector. The symmetry of the metric under time translation is intimately related to energy conservation via $Q_t$ and thus the stability of the simulation. In the absence of a potential, when the metric does not depend on neither $t$ nor $x$, our discretized approach, in addition to $K_t=(1,0)$, retains the continuum invariance under shifts in $x$ via the Killing vector $K_x=(0,1)$, as well as the invariance under boosts via the Killing vector $K_\eta=(x,t)$.

We will give numerical evidence that we achieve exact conservation of $Q_t$ in the interior of the simulated domain, even in the case of highly non-harmonic motion. In contrast to other formally energy preserving schemes, such as the leap-frog, our approach, using SBP operators, is consistent with the continuum formulation, in that it only requires the actual initial conditions of the system at hand, avoiding the need to stagger the degrees of freedom (also known as insertion of dummy points).

After introducing the discretization on the level of the underlying action functional, we will obtain the classical trajectory by numerically finding the critical point of that functional without the need to derive the corresponding equations of motion. To make sure that the solution of the discretized variational principle mimics as accurately as possible the continuum theory, we deploy summation-by-parts finite difference operators \cite{svard2014review,fernandez2014review,lundquist2014sbp}.

Note that we are discretizing the world-line parameter $\gamma$ with equidistant steps, whereas both the values of $t$ and $x$ arise dynamically from the evolution of the simulation along $\gamma$. I.e. a not necessarily equidistant discretization of the time coordinate emerges dynamically in our approach. As we will see in \cref{sec:num} this dynamical time discretization realizes a one-dimensional form of automatic adaptive mesh refinement, guided by the symmetries of the system. I.e. the non-equidistant discretization in $t$ plays a crucial role in guaranteeing that the Noether charge $Q_t$ remains conserved. 

Another non-standard feature of our technique is the departure from the conventional notion of carrying out a simulation on a predefined time interval. We instead provide the initial time and its velocity with respect to $\gamma$, so that the end-point of the simulation too emerges dynamically. 

In the following we will consider the trajectory of a point particle propagating under the influence of an arbitrary $x$ but not $t$ dependent potential $V(x)$. We begin by discretizing the action functional ${\cal E}_{\rm IVP}$ of \cref{eq:EIVPcont} along the world-line parameter $\gamma$ between $\gamma_i$ and $\gamma_f$ with $N_\gamma$ steps, leading to a step-size of $d\gamma=(\gamma_f-\gamma_i)/(N_\gamma-1)$. We will add to $E_{\rm IVP}$ Lagrange multipliers to explicitly account for both the initial conditions and the connecting conditions required by doubling of the degrees of freedom. The forward and backward paths $x_{1,2}$ and times $t_{1,2}$ are described by ${\bf x}_{1,2}=(x_{1,2}(0),x_{1,2}(\Delta \gamma),x_{1,2}(2\Delta \gamma),\ldots,x_{1,2}((N_\gamma-1)\Delta\Gamma))^{\rm T}$ and ${\bf t}_{1,2}=(t_{1,2}(0),t_{1,2}(\Delta \gamma),t_{1,2}(2\Delta \gamma),\ldots,t_{1,2}((N_\gamma-1)\Delta\gamma))^{\rm T}$ respectively.

The integral in ${\cal E}_{\rm IVP}$ is approximated with a quadrature rule, consistent with our choice of finite difference operator, in the form of a diagonal positive definite matrix $\mathds{H}$. The inner product on discretized paths and times thus reads $({\bf x},{\bf x}^\prime)={\bf x}^{\rm T} \mathds{H} {\bf x}^\prime$. 

With integration by parts being a central element in establishing both equations of motion and the existence of conserved quantities, we must use a discretization that mimics IBP exactly, which is achieved by deploying summation-by-parts (SBP) operators $\mathds{D}$ with the defining properties
\begin{align}
    \mathds{D}=\mathds{H}^{-1}\mathds{Q}, \qquad \mathds{Q}^{\rm T}+\mathds{Q}=\mathds{E}_N-\mathds{E}_0={\rm diag}[-1,0,\ldots,0,1].
\end{align}

In this study we consider both the lowest order SBP discretization scheme, referred to as \texttt{SBP21} and the next higher order scheme \texttt{SBP42}. The former is second order in the interior and exhibits one order less on the boundary. Using the trapezoidal rule for integration one has 
\begin{equation}
 \mathds{H}^{[2,1]}=\Delta \gamma \left[ \begin{array}{ccccc} 1/2 & & & & \\ &1 & & &\\ & &\ddots && \\ &&&1&\\ &&&&1/2 \end{array} \right],
\quad 
\mathds{D}^{[2,1]}=
\frac{1}{2 \Delta \gamma}
\left[ \begin{array}{ccccc} -2 &2 & & &\\ -1& 0& 1& &\\ & &\ddots && \\ &&-1&0&1\\ &&&-2&2 \end{array} \right].\label{eq:SBP21}
\end{equation}
The \texttt{SBP42} scheme achieves fourth order accuracy in the interior, which reduces to second order on the boundary
\begin{align}
\nonumber &\mathds{H}^{[4,2]}=\Delta \gamma \left[ \begin{array}{ccccccc} 
\frac{17}{48} & & & & & & \\
 & \frac{59}{48} & & & & & \\
 & &\frac{43}{48} & & & & \\
 & & &\frac{49}{48} & & & \\
 & & & & &1 & \\
 & & & & & &\ddots \\
 \end{array} \right],\end{align}

 \begin{align}
&\mathds{D}^{[4,2]}=
\frac{1}{\Delta \gamma}
\left[ \begin{array}{ccccccccc} 
-\frac{24}{17}&\frac{59}{34} & -\frac{4}{17} & -\frac{3}{34} & && & & \\
-\frac{1}{2}& 0 & \frac{1}{2} & 0 & && & & \\
 \frac{4}{43}& -\frac{59}{86} & 0 & \frac{59}{86}&-\frac{4}{43} && & & \\
\frac{3}{98}& 0& -\frac{59}{86}  & 0&\frac{32}{49}&-\frac{4}{49}& & & \\
& &\frac{1}{12}  & -\frac{2}{3}&0&\frac{2}{3}& -\frac{1}{12}& & \\
&&&&&&&\ddots&
 \end{array} \right].\label{eq:SBP42}
\end{align}

The SBP operators defined above are not yet ready for duty in our variational approach, as they allow for non-physical zero modes. As discussed in detail in \cite{Rothkopf:2022zfb}, we can construct null-space consistent\footnote{Note that in the context of PDE's, SBP operators are considered null-space consistent by construction, as only their right eigenvectors play a role in the equation of motion. Here due to the presence of $\mathds{D}^{\rm T}$ in the action functional, also the left eigenvectors contribute, among which a highly oscillating null-mode (the so-called $\pi$-mode) can be identified (see ref.~\cite{Rothkopf:2022zfb})} SBP operators $\bar{\mathds{D}}$ from the conventional $\mathds{D}$ by deploying affine coordinates and by absorbing penalty terms, inspired by the simultaneous-approximation terms (SAT) technique \cite{carpenter1994time}, used to regularize SBP operators. A brief overview of this regularization is given in \cref{sec:appregSBP}.

The idea behind the penalty term construction is that we are assigning a penalty to all functions that do not fulfill the initial conditions in $t$ and $x$, which includes the non-physical zero mode of $\mathds{D}$. In turn, when we will be searching for the critical point of the discretized action functional ${\cal E}_{\rm IVP}$ the minimizer will approach the correct solution globally and the presence of the penalty term effectively prevents contamination of the correct solution by the non-constant zero mode. 

Explicitly our regularized and null-space consistent \texttt{SBP21} operators read
\begin{align}
&{\bar{\mathds{D}}}^{{\rm R},[2,1]}_t=
\left[ \begin{array}{cccccc} -\frac{1}{\Delta \gamma} + \frac{2}{\Delta \gamma}&\frac{1}{\Delta \gamma} & & &&-\frac{2}{\Delta \gamma} t_i\\ -\frac{1}{2\Delta \gamma}& 0& \frac{1}{2\Delta \gamma}& &&0\\ & &\ddots &&&\vdots\\ &&-\frac{1}{2\Delta \gamma}&0&\frac{1}{2\Delta \gamma}&0\\ &&&-\frac{1}{\Delta \gamma}&\frac{1}{\Delta \gamma}&0\\  0&&\ldots &&0&1\\  \end{array} \right]\label{eq:SBPregt}\\
&{\bar{\mathds{D}}}^{{\rm R},[2,1]}_x=
\left[ \begin{array}{cccccc} -\frac{1}{\Delta \gamma} + \frac{2}{\Delta \gamma}&\frac{1}{\Delta \gamma} & & &&- \frac{2}{\Delta \gamma} x_i\\ -\frac{1}{2\Delta \gamma}& 0& \frac{1}{2\Delta \gamma}& &&0\\ & &\ddots &&&\vdots\\ &&-\frac{1}{2\Delta \gamma}&0&\frac{1}{2\Delta \gamma}&0\\ &&&-\frac{1}{\Delta \gamma}&\frac{1}{\Delta \gamma}&0\\  0&&\ldots &&0&1\\  \end{array} \right].\label{eq:SBPregx}
\end{align}

Using the operators defined above, we can now write the discretized action functional in the following fashion
\begin{align}
\nonumber \mathds{E}_{\rm IVP}= &\frac{1}{2} \left\{ ({\bar{\mathds{D}}^{\rm R}}_t{\bf t}_1)^{\rm T} \mathbb{d}\left[c^2+\frac{2 {\bf V}({\bf x}_1)}{m}\right] {\bar{\mathds{H}}} ({\bar{\mathds{D}}^{\rm R}}_t{\bf t}_1) -  ({\bar{\mathds{D}}^{\rm R}}_x{\bf x}_1)^{\rm T} {\bar{\mathds{H}}} ({\bar{\mathds{D}}^{\rm R}}_x{\bf x}_1)\right\}\\
\nonumber-&\frac{1}{2} \left\{({\bar{\mathds{D}}^{\rm R}}_t{\bf t}_2)^{\rm T} \mathbb{d}\left[c^2+\frac{2 {\bf V}({\bf x}_2)}{m}\right] {\bar{\mathds{H}}} ({\bar{\mathds{D}}^{\rm R}}_t{\bf t}_2) -  ({\bar{\mathds{D}}^{\rm R}}_x{\bf x}_2)^{\rm T} {\bar{\mathds{H}}} ({\bar{\mathds{D}}^{\rm R}}_x{\bf x}_2)\right\}\\
\nonumber+&\lambda_1\big({\bf t}_1[1]-t_i\big)+\lambda_2\big((\mathds{D}{\bf t}_1)[1]-\dot t_i\big)+\lambda_3\big({\bf x}_1[1]-x_i\big)\\
\nonumber+&\lambda_4\big((\mathds{D}{\bf x}_1)[1]-\dot x_i\big)\\
\nonumber+&\lambda_5\big({\bf t}_1[N_\gamma]-{\bf t}_2[N_\gamma]\big) + \lambda_6\big({\bf x}_1[N_\gamma]-{\bf x}_2[N_\gamma]\big)\\
+&\lambda_7\big( (\mathds{D}{\bf t}_1)[N_\gamma]- (\mathds{D}{\bf t}_2)[N_\gamma]\big)+\lambda_8\big( (\mathds{D}{\bf x}_1)[N_\gamma]- (\mathds{D}{\bf x}_2)[N_\gamma]\big)\label{eq:discEIVP}.
\end{align}
Conventional matrix vector multiplication is implied in the above expression, whenever a matrix quantity such as $\bar{\mathds{H}}$ or $\bar{\mathds{D}}$ acts on a vector ${\bf x}_{1,2}$ or ${\bf t}_{1,2}$. The matrix denoted by $\mathbb{d}[f({\bf x})]$ contains on its diagonal the values $\mathbb{d}_{kk}=f({\bf x}(\gamma_k))$ and zero otherwise. We deploy an appropriately modified matrix ${\bar{\mathds{H}}}$ for the inner product in the presence of the affine-coordinate regularized SBP operators (see \cref{sec:appregSBP}).

The initial conditions we supply are the values of the spatial and temporal coordinate $x_i, t_i$, as well as the initial velocities with respect to the world line parameter $\gamma$, i.e. $\dot x_i$ and $\dot t_i$. Since our physical problem is formulated as an initial value problem, given $t_i$, $x_i$ and the physical velocity $v_i=dx/dt$, there exists a freedom to choose $\dot t_i$ and $\dot x_i$, as only their ratio is fixed $v_i=\dot x_i/\dot t_i$. We have added eight Lagrange multipliers, whose role is to explicitly implement the initial conditions $(\lambda_{1-4})$ and the connecting conditions at the end of the forward and backward branches of our doubled degree of freedom construction $(\lambda_{5-8})$.

Once the action functional has been formulated in its discrete form, changing from  \texttt{SBP21} to \texttt{SBP42} only requires replacement of the corresponding difference operator $\mathds{D}$ and quadrature matrix $\mathds{H}$ but no further changes to the functional itself.

This concludes the description of our novel variational approach and we proceed to evaluate its properties and performance based on two concrete numerical examples.

\section{Numerical results}
\label{sec:num}

In this section we will present explicit results for the numerically obtained classical trajectory of a point particle in the presence of two different potentials, $V_1(x)=\alpha x$ 
 and $V_4(x)=\kappa x^4$. These two choices correspond to a model of a point mass falling in a constant gravitational field 
and carrying out highly-nonlinear anharmonic motion. We set the mass of the particle to unity, as well as adopt without loss of generality the convention that the speed of light $c=1$, which simply amounts to a particular choice of units for length and time.

Let us stress again that while standard numerical methods exist to solve the equations of motion for each of these systems, the novelty of the approach presented here lies in the fact that we retain the continuum time shift invariance of the system and thus achieve \textit{exact conservation} of $Q_t$ in the interior of the simulated time domain. In addition we determine the classical trajectory directly from the action functional of the geometrized problem, without the need to derive the equation of motion.

We implement the action functional \cref{eq:discEIVP} in the Mathematica language\footnote{The code using both the \texttt{SBP21} or \texttt{SBP42} operator is available under open access on the Zenodo repository \cite{zenodoIVP2023}.}. As the critical point of the action may be a saddle point, instead of an actual minimum, we must be careful in deploying established numerical optimization algorithms in the dynamical degrees of freedom ${\bf d}=\{{\bf t}_{1,2},{\bf x}_{1,2},\lambda_{1-8}\}$. Instead of minimizing $\mathds{E}_{\rm IVP}$ directly, we will minimize the Euclidean norm of the gradient $|\nabla_{\bf d}\mathds{E}_{\rm IVP}|^2$. Via this detour, a saddle point is converted into a minimum. In practice we deploy a chain of minimization algorithms. We start with a preconditioning based on the LBFGS quasi-Newton algorithm, which features cost efficient iteration steps, when far away from the true critical point. It is followed by further iterations based on the full Newton method, which exhibits a faster convergence rate than the LBFGS algorithm when close to the critical point. Once the critical point has been approached to at least floating point precision we switch to the interior point optimization, which showed reliable performance in identifying the critical point to any desired tolerance. For our numerical tests in Mathematica, we used \texttt{WorkingPrecision} of $40$ and \texttt{PrecisionGoal} of $40$.

The figures shown in the following are based on results from the \texttt{SBP21} operator and include the outcomes from the \texttt{SBP42} operators when indicated in the text.

\subsection{Linear potential case}

We discretize the continuous action functional 
\begin{align}
    \nonumber{\cal E}^{\rm lin}_{\rm IVP}=&\int_{\gamma_i}^{\gamma_f}d\gamma \frac{1}{2}\left\{ \Big( 1+2\alpha x_1(\gamma) \Big) \Big(\frac{d t_1}{d\gamma}\Big)^2-\Big(\frac{d x_1}{d\gamma}\Big)^2\right\}\\
   \nonumber -&\int_{\gamma_i}^{\gamma_f}d\gamma \frac{1}{2}\left\{ \Big( 1+2\alpha x_2(\gamma) \Big) \Big(\frac{d t_2}{d\gamma}\Big)^2-\Big(\frac{d x_2}{d\gamma}\Big)^2\right\}\\
   \nonumber+&\lambda_1\big(t(\gamma_i)-t_i\big)+\lambda_2\big(\dot t_1(\gamma_i)-\dot t_i\big)+\lambda_3\big(x_1(\gamma_i)-x_i\big)+\lambda_4\big(\dot x_1(\gamma_i)-\dot x_i\big)\\
\nonumber   +&\lambda_5\big(t_1(\gamma_f)-t_2(\gamma_2)\big) +\lambda_6\big(\dot t_1(\gamma_f)-\dot t_2(\gamma_2)\big)\\
   +&\lambda_7\big(x_1(\gamma_f)-x_2(\gamma_2)\big) +\lambda_8\big(\dot x_1(\gamma_f)-\dot x_2(\gamma_2)\big)
\end{align}
along the world-line of the particle motion between $\gamma_i=0$ and $\gamma_f=1$ with $N_\gamma=32$ points. Without loss of generality, we arbitrarily set the starting time to $t_i=0$ and the starting position to $x_i=1$. To obtain an initial velocity $v_i=1/10$ we choose $\dot t=1$ and $\dot x=v_i$. Note that we do not fix the value of $t_f$ but only the initial velocity of time with respect to $\gamma$. The choice of $\dot t=1$ will lead to dynamics, such that $t_f$ will be of the order of one. (In the next subsection we will also provide results for different choices of $\dot t_i$.) As strength for the linear potential we choose $\alpha=1/4$. The corresponding discrete action functional reads explicitly
\begin{align}
\nonumber \mathds{E}^{\rm lin}_{\rm IVP}= &\frac{1}{2} \left\{ ({\bar{\mathds{D}}^{\rm R}}_t{\bf t}_1)^{\rm T} \mathbb{d}\left[1+2 \alpha {\bf x}_1\right]{\bar{\mathds{H}}} ({\bar{\mathds{D}}^{\rm R}}_t{\bf t}_1) -  ({\bar{\mathds{D}}^{\rm R}}_x{\bf x}_1)^{\rm T} {\bar{\mathds{H}}} ({\bar{\mathds{D}}^{\rm R}}_x{\bf x}_1)\right\}\\
\nonumber-&\frac{1}{2} \left\{ ({\bar{\mathds{D}}^{\rm R}}_t{\bf t}_2)^{\rm T} \mathbb{d}\left[1+2 \alpha {\bf x}_2\right] {\bar{\mathds{H}}} ({\bar{\mathds{D}}^{\rm R}}_t{\bf t}_2) -  ({\bar{\mathds{D}}^{\rm R}}_x{\bf x}_2)^{\rm T} {\bar{\mathds{H}}} ({\bar{\mathds{D}}^{\rm R}}_x{\bf x}_2)\right\}\\
\nonumber+&\lambda_1\big({\bf t}_1[1]-t_i\big)+\lambda_2\big((\mathds{D}{\bf t}_1)[1]-\dot t_i\big)\\
\nonumber+&\lambda_3\big({\bf x}_1[1]-x_i\big)+\lambda_4\big((\mathds{D}{\bf x}_1)[1]-\dot x_i\big)\\
\nonumber+&\lambda_5\big({\bf t}_1[N_\gamma]-{\bf t}_2[N_\gamma]\big) + \lambda_6\big({\bf x}_1[N_\gamma]-{\bf x}_2[N_\gamma]\big)\\
+&\lambda_7\big( (\mathds{D}{\bf t}_1)[N_\gamma]- (\mathds{D}{\bf t}_2)[N_\gamma]\big)+\lambda_8\big( (\mathds{D}{\bf x}_1)[N_\gamma]- (\mathds{D}{\bf x}_2)[N_\gamma]\big)\label{eq:discEIVPlin}.
\end{align}

\begin{figure}
    \includegraphics[scale=0.3]{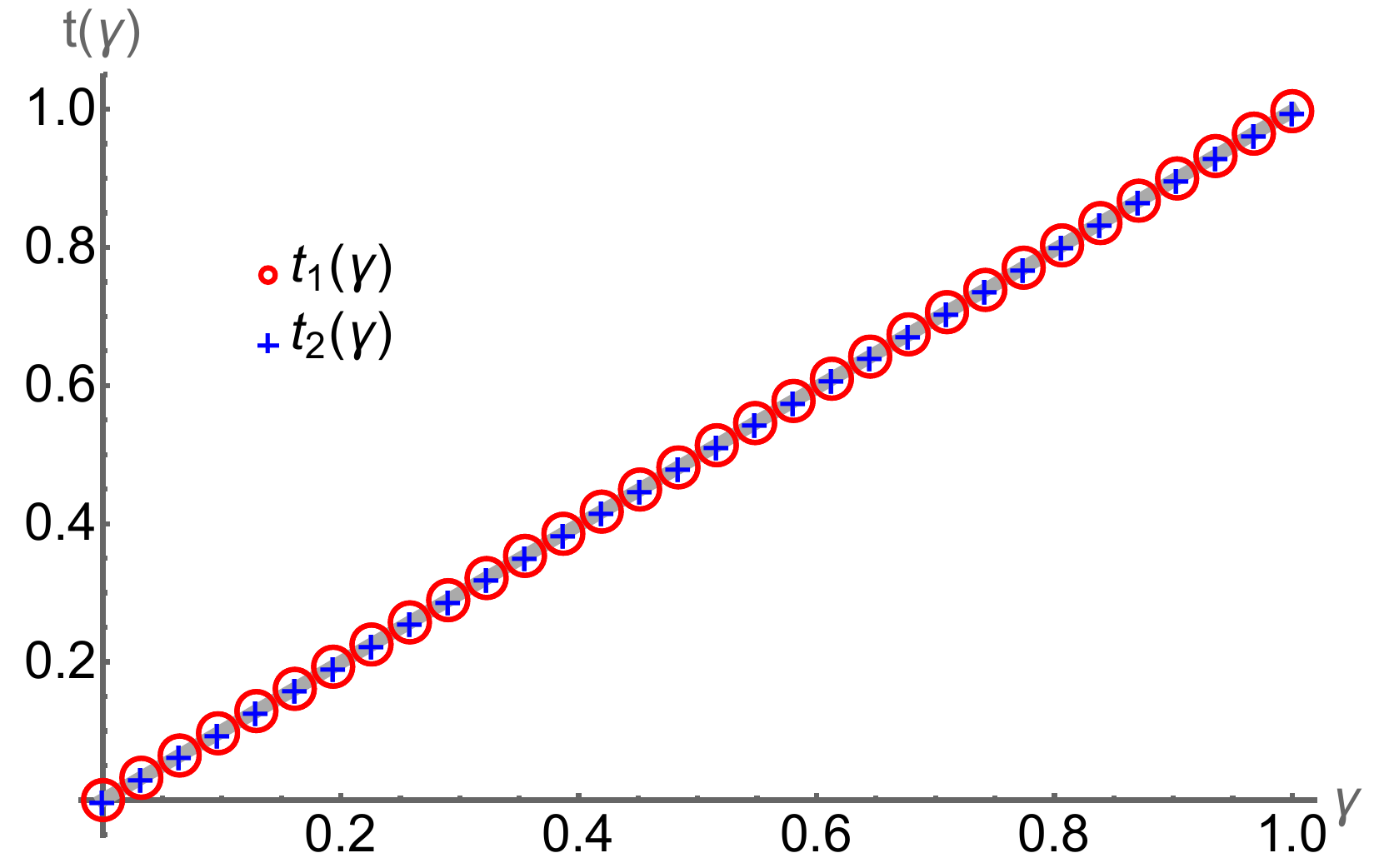}
    \includegraphics[scale=0.3]{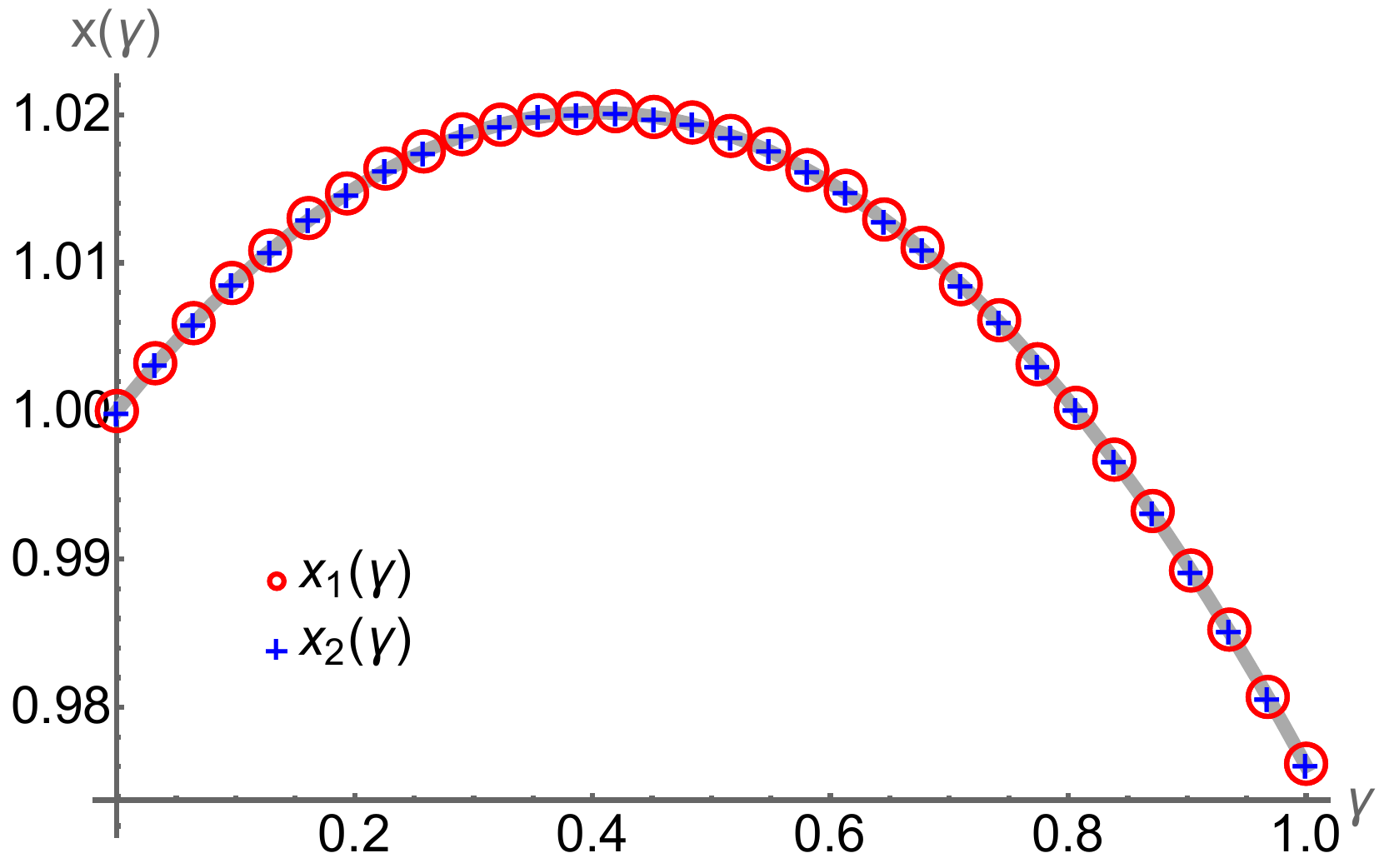}
    \caption{The values of the (top) time coordinates ${\bf t}_1(\gamma_i)$ (red circles) and ${\bf t}_2(\gamma_i)$ (blue crosses) and that of the (bottom) spatial coordinates ${\bf x}_1(\gamma_i)$ (red circles) and ${\bf x}_2(\gamma_i)$ (blue crosses) along the world-line parameter $\gamma$, as obtained from the critical point of $\mathds{E}^{\rm lin}_{\rm IVP}$ with $V(x)=\alpha x$, discretized with $N_\gamma=32$ points and the \texttt{SBP21} operators. The solution of the corresponding geodesic equations via Mathematica's \texttt{NDSolve} is shown as light gray solid line.}
    \label{fig:Gravdofofgamma}
\end{figure}

Let us take a look in \cref{fig:Gravdofofgamma} at the raw results for the forward and backward time and spatial coordinates, as obtained from the critical point of $\mathds{E}^{\rm lin}_{\rm IVP}$ with $V(x)=\alpha x$. In the top panel, we show ${\bf t}_1(\gamma_i)$ as red circles and ${\bf t}_2(\gamma_i)$ as blue crosses, while in the bottom panel these symbols denote the spatial coordinate of the point particle trajectory ${\bf x}_1(\gamma_i)$ and ${\bf x}_2(\gamma_i)$ respectively. As required by the physical limit (discussed in \cref{sec:initvalform}), we find that the values of the doubled degrees of freedom coincide at the critical point. The solution of the corresponding continuum geodesic equations, obtained via the \texttt{LSODA} algorithm of Mathematica's \texttt{NDSolve} command is shown as gray solid line and excellent agreement is observed.
Note that due to our choice of $\dot t_i=1$ the maximum time traversed by the simulation is close to one. 

\begin{figure}
    \includegraphics[scale=0.3]{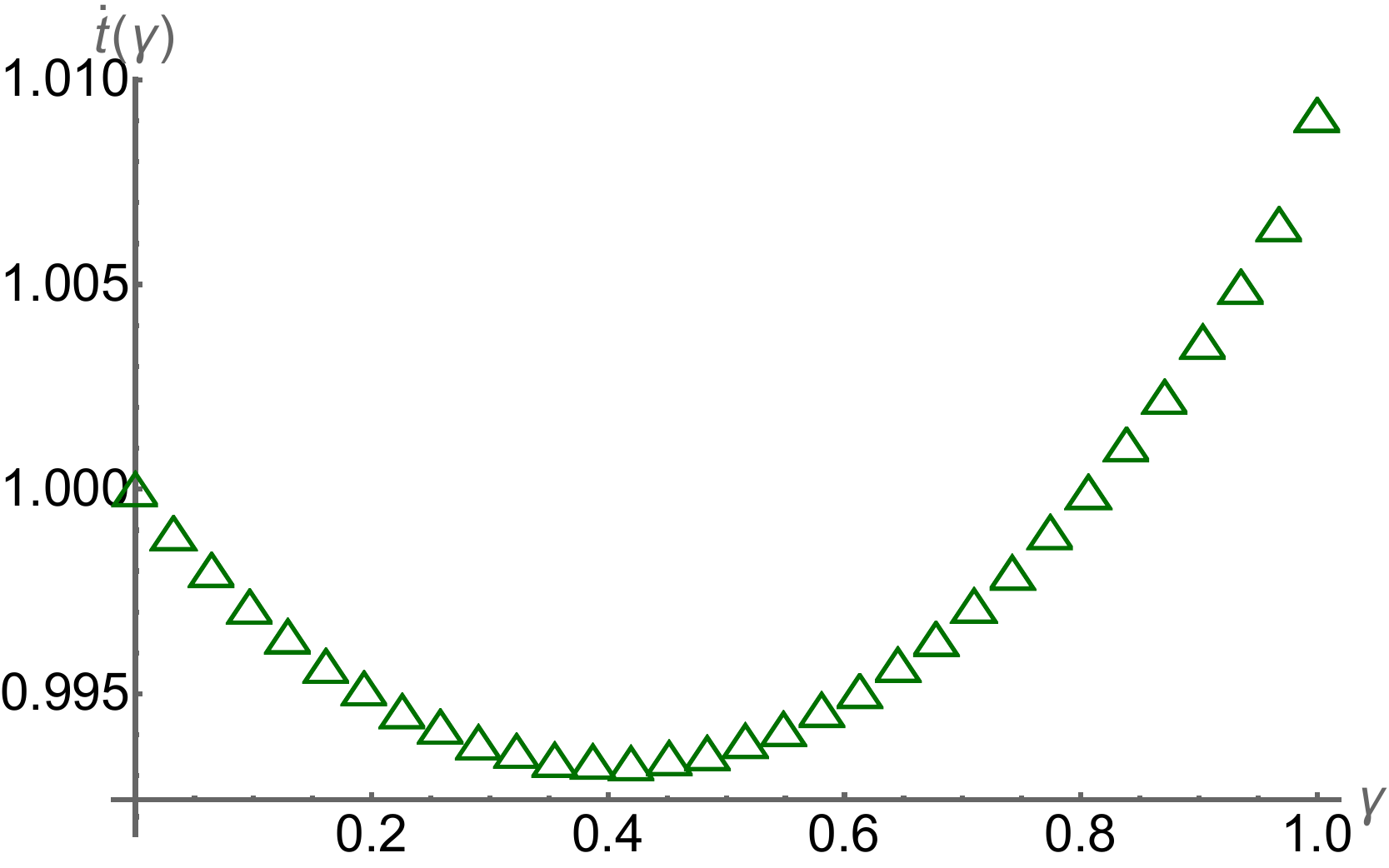}
    \caption{The change of time with respect to the world-line coordinate $\gamma$, as obtained from the critical point of $\mathds{E}_{\rm IVP}$ with $V(x)=\alpha x$.}
    \label{fig:Gravtdotofgamma}
\end{figure}

At first sight it appears that an equidistant discretization of time in $\gamma$ emerges, but an inspection of the velocity of time with respect to $\gamma$ in \cref{fig:Gravtdotofgamma} reveals that the time spacing dynamically adapts to the behavior observed in the spatial coordinate $x$. Close to the maximum of $x(\gamma)$ at around $\gamma=0.4$ the temporal spacing e.g. has a minimum. This dynamically emerging time discretization constitutes an automatically generated non-trivial mesh for the time coordinate and arises naturally in our formalism. In fact an automatic AMR procedure results.

\begin{figure}
    \includegraphics[scale=0.3]{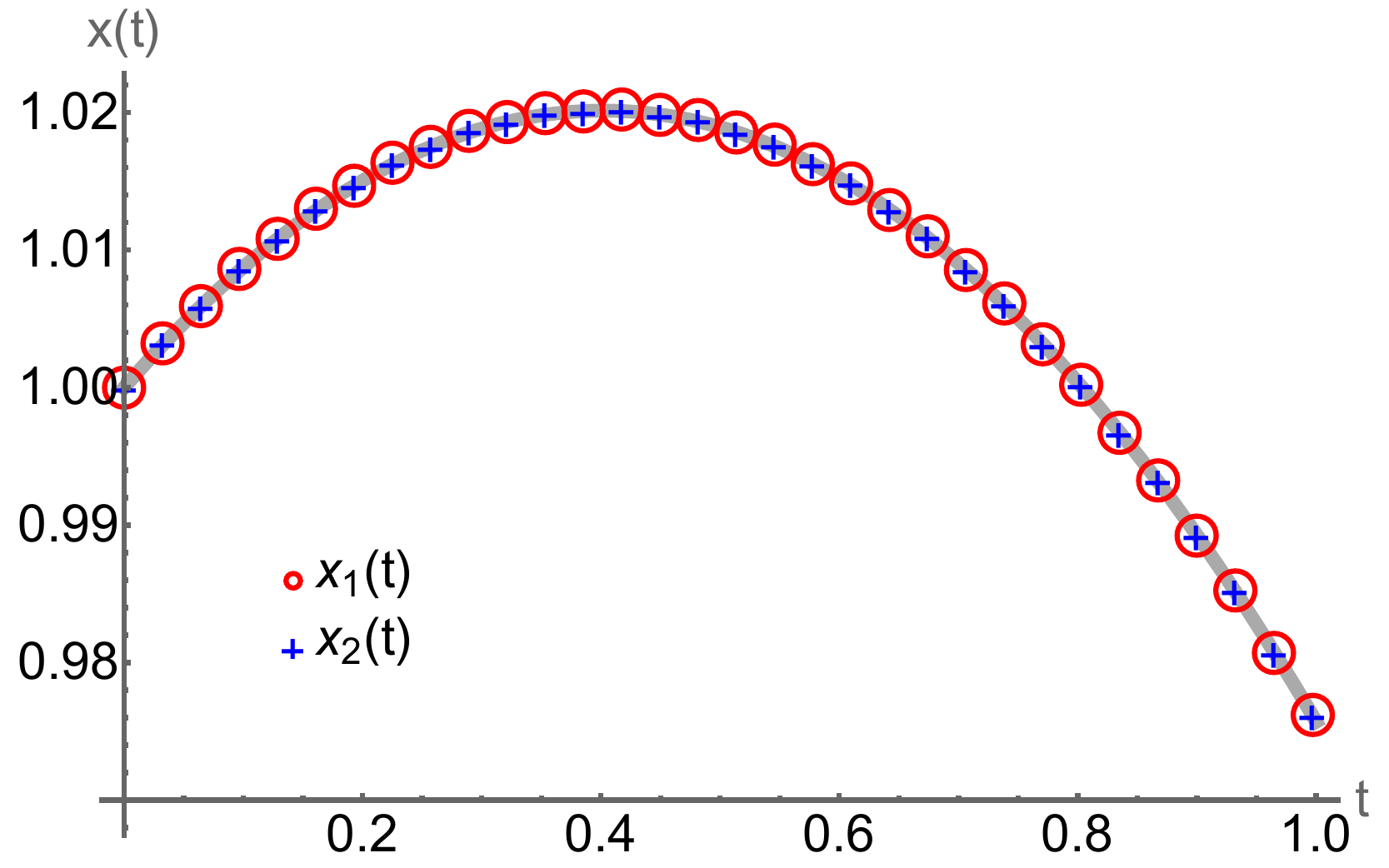}
    \caption{The forward ${\bf x}_1(t_1)$ (red circles) and backward ${\bf x}_2(t_2)$ (blue crosses) degrees of freedom obtained from the critical point of $\mathds{E}^{\rm lin}_{\rm IVP}$ with $V(x)=\alpha x$. We discretize with $N_\gamma=32$ points and the \texttt{SBP21} operators. The trajectory from a \texttt{LSODA} solver of the corresponding equation of motion is given as gray solid line.}
    \label{fig:Gravxoft}
\end{figure}

Let us plot next in \cref{fig:Gravxoft}, the results from our geometrized formalism as physical trajectory, i.e. as ${\bf x}_{1,2}(t_{1,2})$ (red circles and blue crosses). This allows us to compare the outcome to the solution one would obtain by following the conventional approach in the literature (see e.g. chapter 7.9 in \cite{goldstein1980classical}). There one considers time as independent variable and simply adds a potential term to the free relativistic action \cref{eq:flataction} before deriving the corresponding Euler-Lagrange equation, which for the linear potential reads $d^2x/dt^2 = -(\alpha)(1-(dx/dt)^2)^{(3/2)}$. Using the \texttt{LSODA} algorithm of Mathematica's \texttt{NDSolve} command, we compute the solution of this equation of motion and plot it as gray solid line. Excellent agreement with the solution from our variational approach is observed, indicating that the geometrization strategy indeed reproduces the solution of the physical problem at hand.

Note that the change in the velocity of the time coordinate manifests itself here as a slightly denser time grid around the maximum of the trajectory.

After this qualitative visual inspection, let us take a closer look at the properties of the obtained solution. The first question we may ask is how well quantitatively the solution follows the naively discretized geodesic equations for time \cref{eq:geodt} and space \cref{eq:geodx} respectively. The \textit{continuum} geodesic equations for the system at hand read
\begin{align}
 &\frac{d}{d\gamma}\Big(g_{00}\frac{dt}{d\gamma}\Big)=\frac{d}{d\gamma}\Big( \big(1+2 \alpha x\big) \frac{dt}{d\gamma}\Big) 
 =0,\\
 &\frac{d}{d\gamma}\Big(\frac{dx}{d\gamma}\Big) + \frac{1}{2}\frac{\partial g_{00}}{\partial x}\Big( \frac{dt}{d\gamma} \Big)^2=\frac{d^2x}{d\gamma^2} + \alpha\Big( \frac{dt}{d\gamma} \Big)^2=0.
\end{align}
When deriving these equations of motion from the continuum action functional \cref{eq:EIVPcont} we have only used integration by parts. This motivates us to proceed, considering them naively discretized by replacing the derivatives with SBP finite difference operators
\begin{align}
    & \mathds{D}\big( (1+2\alpha {\bf x} )\circ(\mathds{D}{\bf t}) \big)=\Delta{\bf G}^t,\\
    &\mathds{D}\mathds{D}{\bf x} + \alpha (\mathds{D}{\bf t})\circ (\mathds{D}{\bf t})=\Delta{\bf G}^x.
\end{align}
Here element-wise multiplication of entries of vector quantities is explicitly denoted by the symbol $\circ$, which implements e.g. ${\bf x}_1\circ {\bf x}_2= ( x_1(0)x_2(0), x_1(\Delta \gamma)x_2(\Delta\gamma),\ldots, x_1(\Delta \gamma(N_\gamma-1))x_2(\Delta\gamma(N_\gamma-1)))^{\rm T}$. Note that we have introduced on the right of the above equations two quantities $\Delta {\bf G}^x$ and $\Delta {\bf G}^t$, which denote the deviation from the value zero, to which the equations of motion evaluate in the continuum. By inspecting $\Delta {\bf G}^x$ and $\Delta {\bf G}^t$ for the trajectories ${\bf x}_{1,2}$ and ${\bf t}_{1,2}$ obtained from the critical point of the discretized action functional $\mathds{E}^{\rm lin}_{\rm IVP}$, we can obtain first quantitative insight into the performance of our variational approach. 

\begin{figure}
    \includegraphics[scale=0.3]{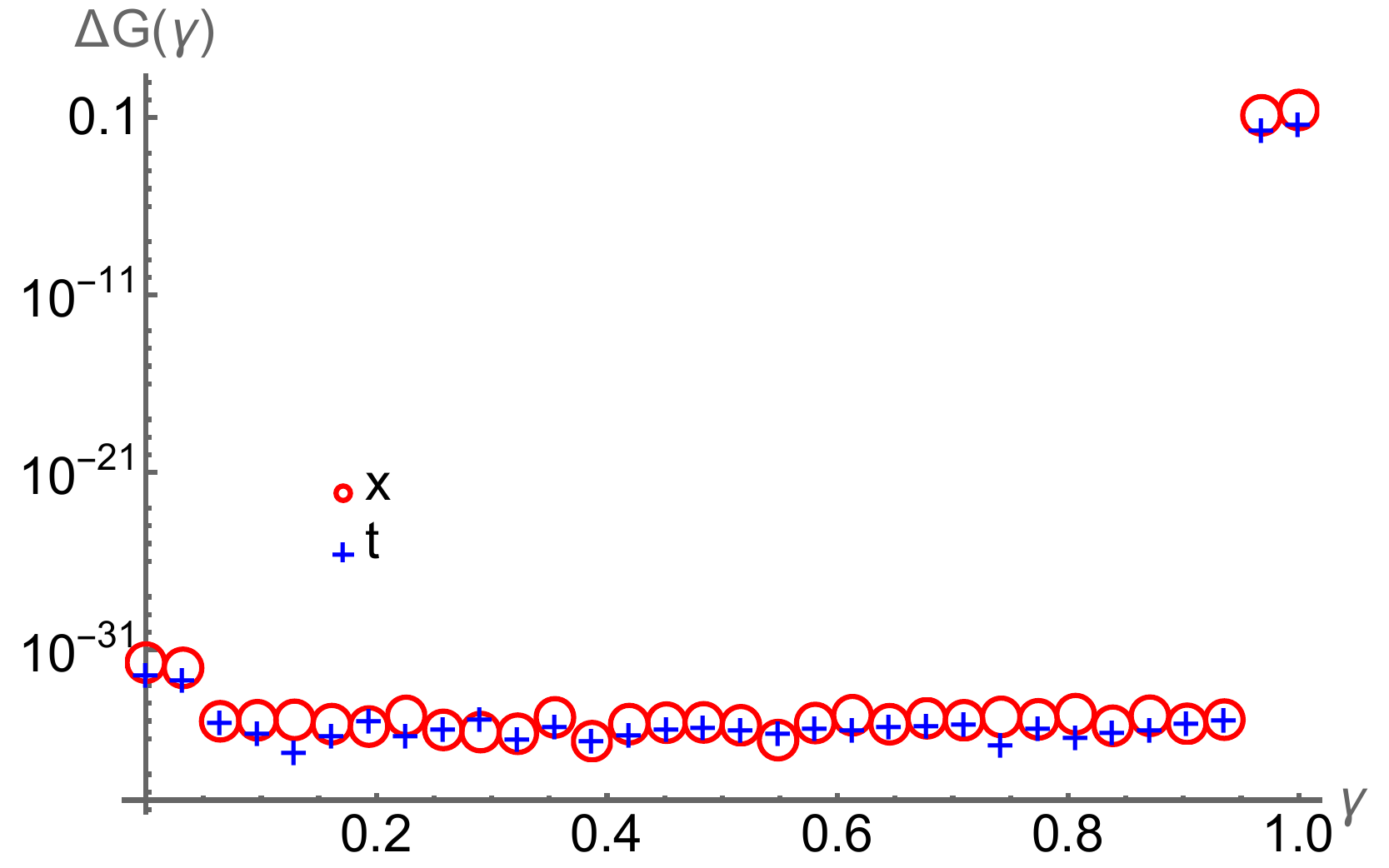}
    \includegraphics[scale=0.3]{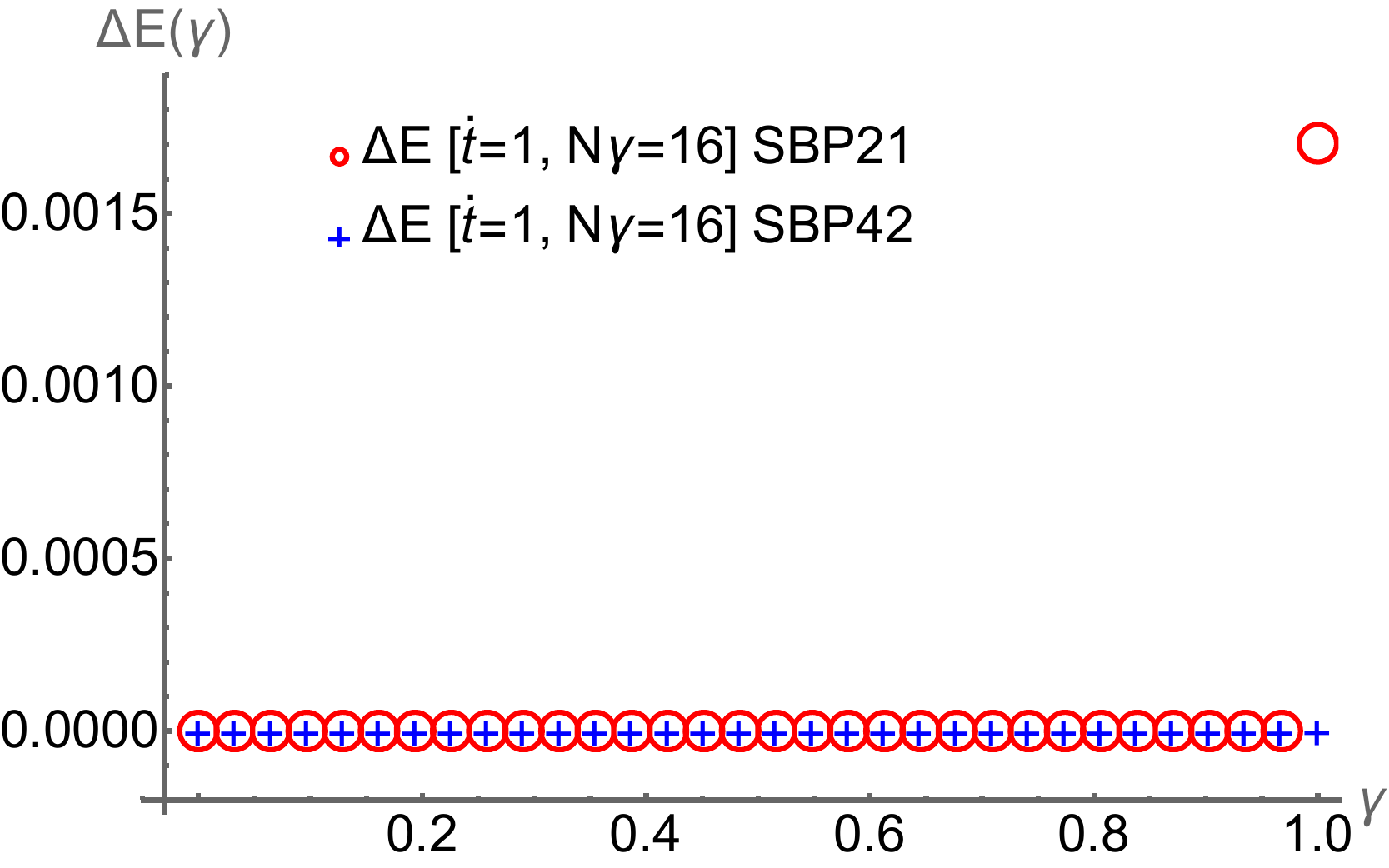}
    \caption{(top) Deviation of the spatial coordinate $\Delta{\bf G}^x$ (red circles) and time coordinate $\Delta{\bf G}^t$ (blue crosses) from the discretized geodesic equation, as obtained from the critical point of $\mathds{E}^{\rm lin}_{\rm IVP}$ with $V(x)=\alpha x$ at $N_\gamma=32$ discretized with the \texttt{SBP21} operator. (bottom) Deviation $\Delta {\bf E}$ of the quantity ${\mathds Q}_t$ from its continuum value as given by the initial conditions. Results for the \texttt{SBP21} operator are given as red circles and those for \texttt{SBP42} as blue crosses. Note that the deviation $\Delta E$ in the interior is exactly zero within machine precision.}
    \label{fig:GravDeltaGDeltaE}
\end{figure}

We plot the values of both quantities $\Delta{\bf G}^x$ and $\Delta{\bf G}^t$ in the top panel of \cref{fig:GravDeltaGDeltaE}. At first sight we find that deviations from the naively discretized geodesic equations are minute, except for the two last points. Note that the plot is given in logarithmic scale. 

Since we use a minimizer in Mathematica with \texttt{WorkingPrecision} set to 40, the values of $<10^{-30}$ reflect a true zero. It is apparent that both the naively discretized geodesic equation for $x$ and $t$ are fulfilled down to machine precision.

Let us proceed to the central quantity of interest in this study $Q_t$, defined in \cref{eq:Qt}, which in the continuum represents the conserved quantity associated with the time-translation symmetry of the system. We again consider its naively discretized form in the following
\begin{align}
    {\bf \mathds{Q}_t}=(\mathds{D}{\bf t})\circ({\bf 1} + 2\alpha {\bf x}).
\end{align}
With the discrete action functional $\mathds{E}^{\rm lin}_{\rm IVP}$ retaining manifest invariance under shifts in the time coordinates ${\bf t}_{1,2}$ we wish to investigate whether also the discretized $\mathds{Q}_t$ retains its role as conserved Noether charge. To this end let us focus here on the deviation $\Delta {\bf E}$ of $\mathds{Q}_t$ from its continuum value
\begin{align}
    \Delta {\bf E} = {\bf \mathds{Q}_t} -Q_t =  (\mathds{D}{\bf t})\circ({\bf 1} + 2\alpha {\bf x}) - \dot t_i (1+2\alpha x_i).
\end{align}
Note that $\mathds{Q}_t$ takes on the continuum value by construction at the first point in $\gamma$, as there it is defined by the initial conditions. The values obtained for $\Delta {\bf E}$ from the critical point of $\mathds{E}^{\rm lin}_{\rm IVP}$ using either the \texttt{SBP21} (red circles) or \texttt{SBP42} operator (blue crosses) are shown in the bottom panel of \cref{fig:GravDeltaGDeltaE}. There are two important observations to be made. 

First, the discretized quantity $\mathds{Q}_t$ is \textit{exactly conserved} in the discrete setting in the \textit{interior} of the simulated time domain and only at the final point $\gamma_f$ it deviates from that constant. While the deviation $\Delta {\bf E}(\gamma_f)$ in case of the \texttt{SBP21} operator is already smaller than two permille, it reduces even further to a value of $10^{-6}$ when deploying the \texttt{SBP42} operator.

We have investigated various potential reasons for the slight difference at the final point, such as a potential over-constraint from the connecting conditions in \cref{eq:conconI,eq:conconII}, but we have not identified the source as of yet. One venue to explore in the future is whether the exact enforcement of the connecting conditions plays a role, which however requires the development of a genuinely weak formulation of our approach without the use of Lagrange multipliers. It is important to point out that, as we will show explicitly below, the presence of this final differing point does not spoil the convergence to the correct continuum limit.

Secondly, the value of $\mathds{Q}_t$ that remains conserved in the interior agrees with the true continuum value, prescribed by the initial conditions, \textit{within machine precision}. This is a highly non-trivial result, as even in energy preserving schemes, such as the leap-frog, the conserved quantities do not necessarily agree with the continuum ones.

\begin{figure}
    \includegraphics[scale=0.3]{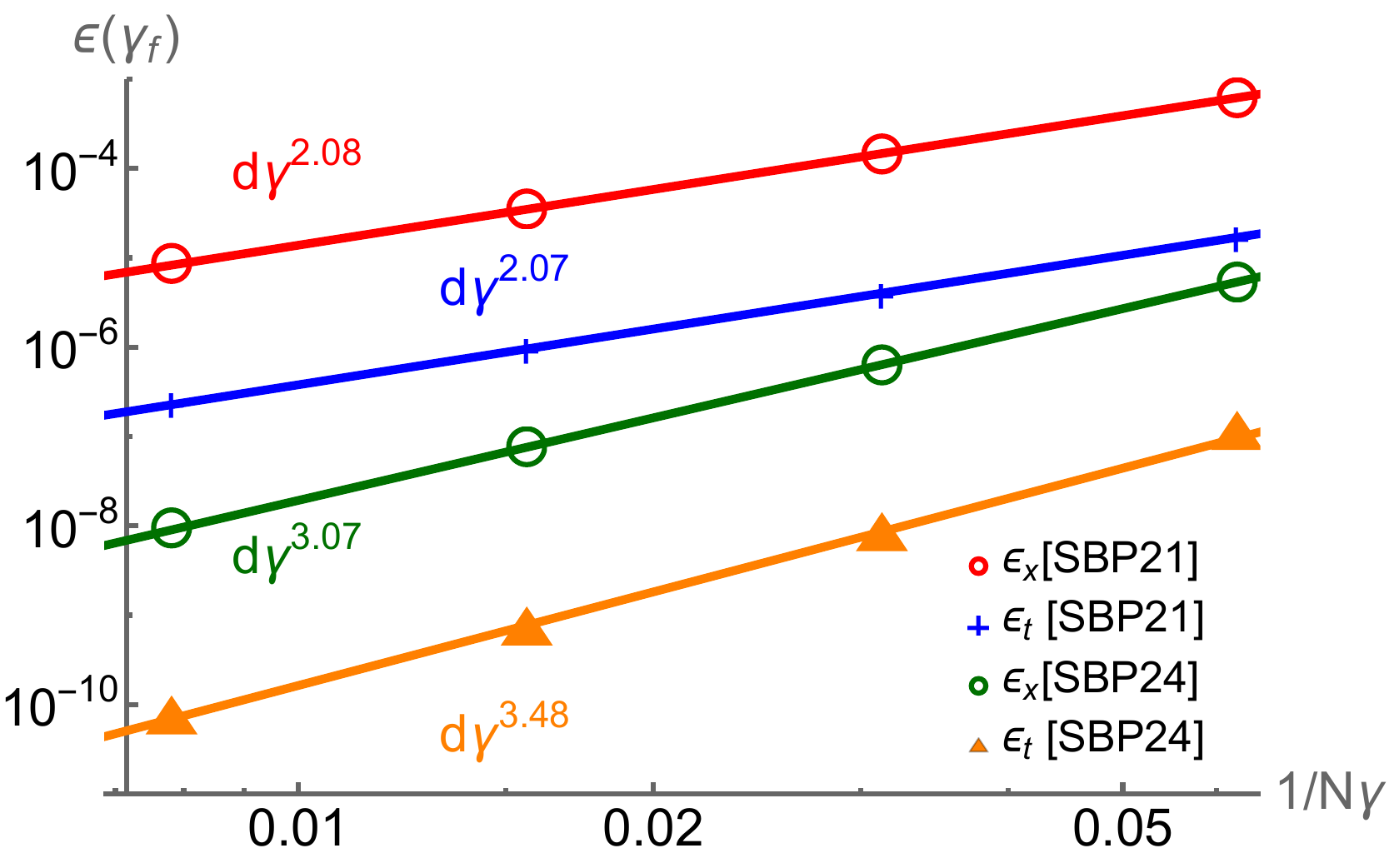}
    \includegraphics[scale=0.3]{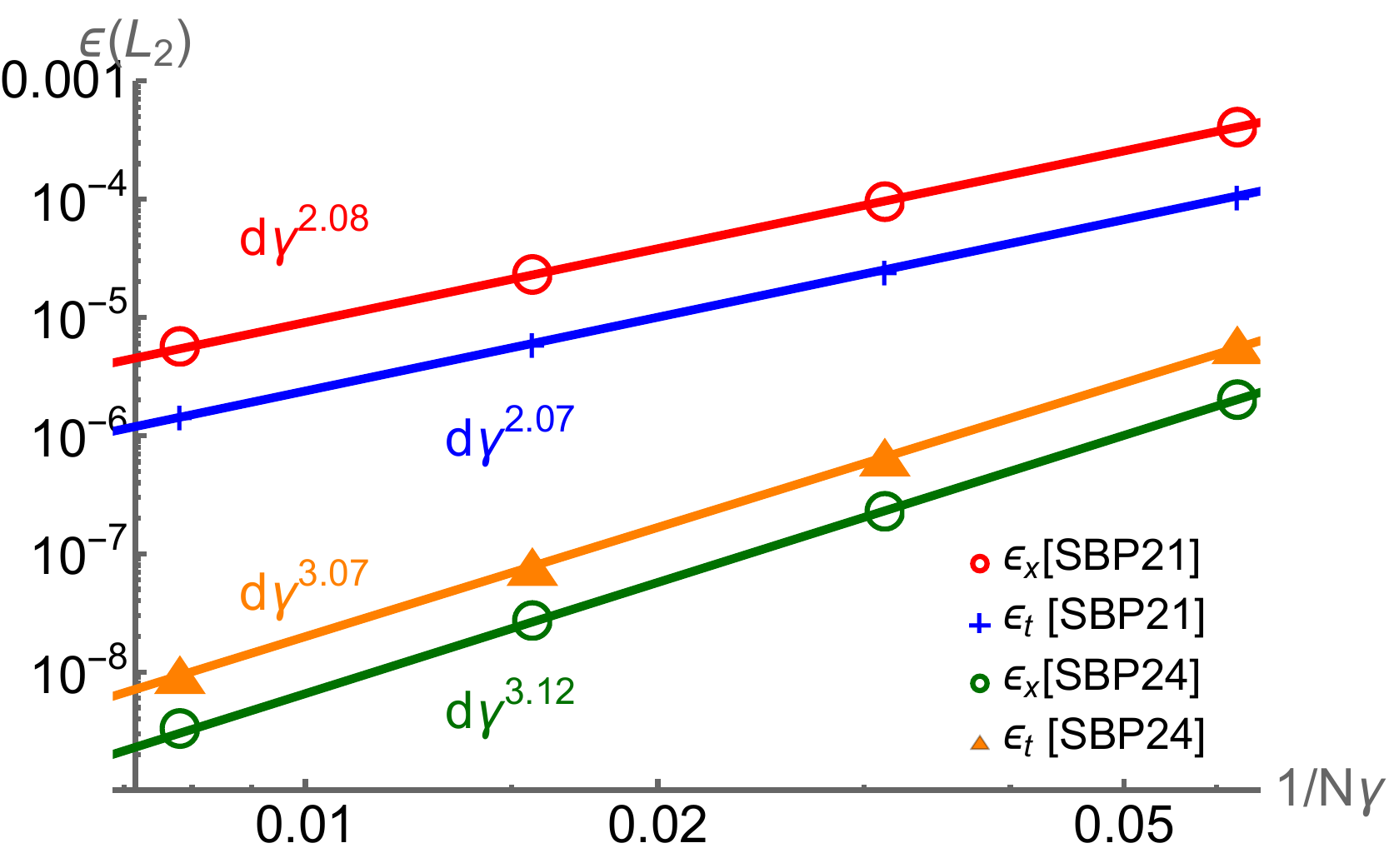}
    \caption{Visualization of the convergence of our novel variational approach under the refinement of the $\gamma$ grid, based on the final point via $\epsilon(\gamma_f)$ (top panel) and globally via $\epsilon(L_2)$ (bottom panel). In each plot we show the results from the \texttt{SBP21} (top two lines) and \texttt{SBP42} (bottom two lines) operators. Red circles and blue crosses denote the absolute deviation in the spatial coordinate $\epsilon_x$ and temporal coordinate $\epsilon_t$ between the continuum solution and the values obtained from the the critical point of $\mathds{E}^{\rm lin}_{\rm IVP}$ with $V(x)=\alpha x$. The corresponding deviations for the \texttt{SBP42} operator are given by the green circles and orange triangles respectively.}
    \label{fig:GravConv}
\end{figure}

We surmise that it is the interplay of a manifest time-translation invariant formulation of the action functional, together with the resulting dynamically emerging time discretization, which achieves the conservation of the discrete $\mathds{Q}_t$ at its continuum value in the interior of the simulation domain.

The presence of two points that deviate from the naively discretized continuum geodesic equations may appear troublesome. However as we show in \cref{fig:GravConv} these points do not spoil the convergence to the correct continuum limit under grid refinement. 

In the top panel of \cref{fig:GravConv}, we select the apparently most disadvantageous points for our convergence study, i.e. we compare the deviation from the continuum geodesic equations $\epsilon(\gamma_f)_x=|{\bf x}[N_\gamma]-x_{\rm true}(1)|$ and $\epsilon(\gamma_f)_t=|{\bf t}[N_\gamma]-t_{\rm true}(1)|$
at $\gamma_f$, exactly where the deviation from the continuum result was maximal in the top panel of \cref{fig:GravDeltaGDeltaE}. Grid refinement is carried out and we provide the results for both the lowest order \texttt{SBP21} operator and the next higher \texttt{SBP42} operator.

Even in this disadvantaged scenario, we find that under grid refinement, the discrete solution approaches the true continuum values as expected from a scheme that is second order in the interior. Taking the \texttt{SBP21} results, the best fit to $\epsilon_x$ reveals a scaling with $\Delta\gamma^{2.08}$, while for $\epsilon_t$ an virtually identical $\Delta\gamma^{2.07}$ ensues. Going over to the \texttt{SBP42} results we find that the convergence is in line with expectations for an SBP operator of 4th order in the interior with $\epsilon_x$ exhibiting a scaling of $\Delta\gamma^{3.07}$ and $\epsilon_t$ a somewhat better value of $\Delta\gamma^{3.48}$.

In the bottom panel of \cref{fig:GravConv} we instead investigate the global convergence of our approach using the $L_2$ norm $\epsilon(L_2)_x=\sqrt{({\bf x}-{\bf x}_{\rm true})^{\rm T}.\mathds{H}.({\bf x}-{\bf x}_{\rm true})}$ and $\epsilon(L_2)_t=\sqrt{({\bf t}-{\bf t}_{\rm true})^{\rm T}.\mathds{H}.({\bf t}-{\bf t}_{\rm true})}$, where ${\bf x}_{\rm true}$ and ${\bf t}_{\rm true}$ are taken from the numerical solution of the geodesic equations, used for comparison in \cref{fig:Gravxoft}. We find that similar convergence rates ensue, where \texttt{SBP21} shows scaling $\Delta\gamma^\beta$ with exponent $\beta\ge2$ and \texttt{SBP42} shows scaling with exponent $\beta\ge 3$.

These convergence result agrees with the findings of our previous study \cite{Rothkopf:2022zfb}, where the standard action functional was discretized with time as independent parameter. 


\subsection{Quartic potential}

After considering the simplest possible non-trivial scenario with a linear potential, we now turn to a system with a quartic potential and the following continuum action functional
\begin{align}
    \nonumber{\cal E}^{\rm qrt}_{\rm IVP}=&\int_{\gamma_i}^{\gamma_f}d\gamma \frac{1}{2}\left\{ \Big( 1+2\kappa x_1^4(\gamma) \Big) \Big(\frac{d t_1}{d\gamma}\Big)^2-\Big(\frac{d x_1}{d\gamma}\Big)^2\right\}\\
   \nonumber -&\int_{\gamma_i}^{\gamma_f}d\gamma \frac{1}{2}\left\{ \Big( 1+2\alpha x_2^4(\gamma) \Big) \Big(\frac{d t_2}{d\gamma}\Big)^2-\Big(\frac{d x_2}{d\gamma}\Big)^2\right\}\\
   \nonumber+&\lambda_1\big(t(\gamma_i)-t_i\big)+\lambda_2\big(\dot t_1(\gamma_i)-\dot t_i\big)+\lambda_3\big(x_1(\gamma_i)-x_i\big)+\lambda_4\big(\dot x_1(\gamma_i)-\dot x_i\big)\\
\nonumber   +&\lambda_5\big(t_1(\gamma_f)-t_2(\gamma_2)\big) +\lambda_6\big(\dot t_1(\gamma_f)-\dot t_2(\gamma_2)\big)\\
   +&\lambda_7\big(x_1(\gamma_f)-x_2(\gamma_2)\big) +\lambda_8\big(\dot x_1(\gamma_f)-\dot x_2(\gamma_2)\big).
\end{align}
Again we discretize along $N_\gamma=32$ in the world-line parameter $\gamma$. Using $\kappa=1/2$ in the potential $V(x)=\kappa x^4$ leads to dynamics that already in the small time regime considered here are distinctly anharmonic.

\begin{figure}
    \includegraphics[scale=0.3]{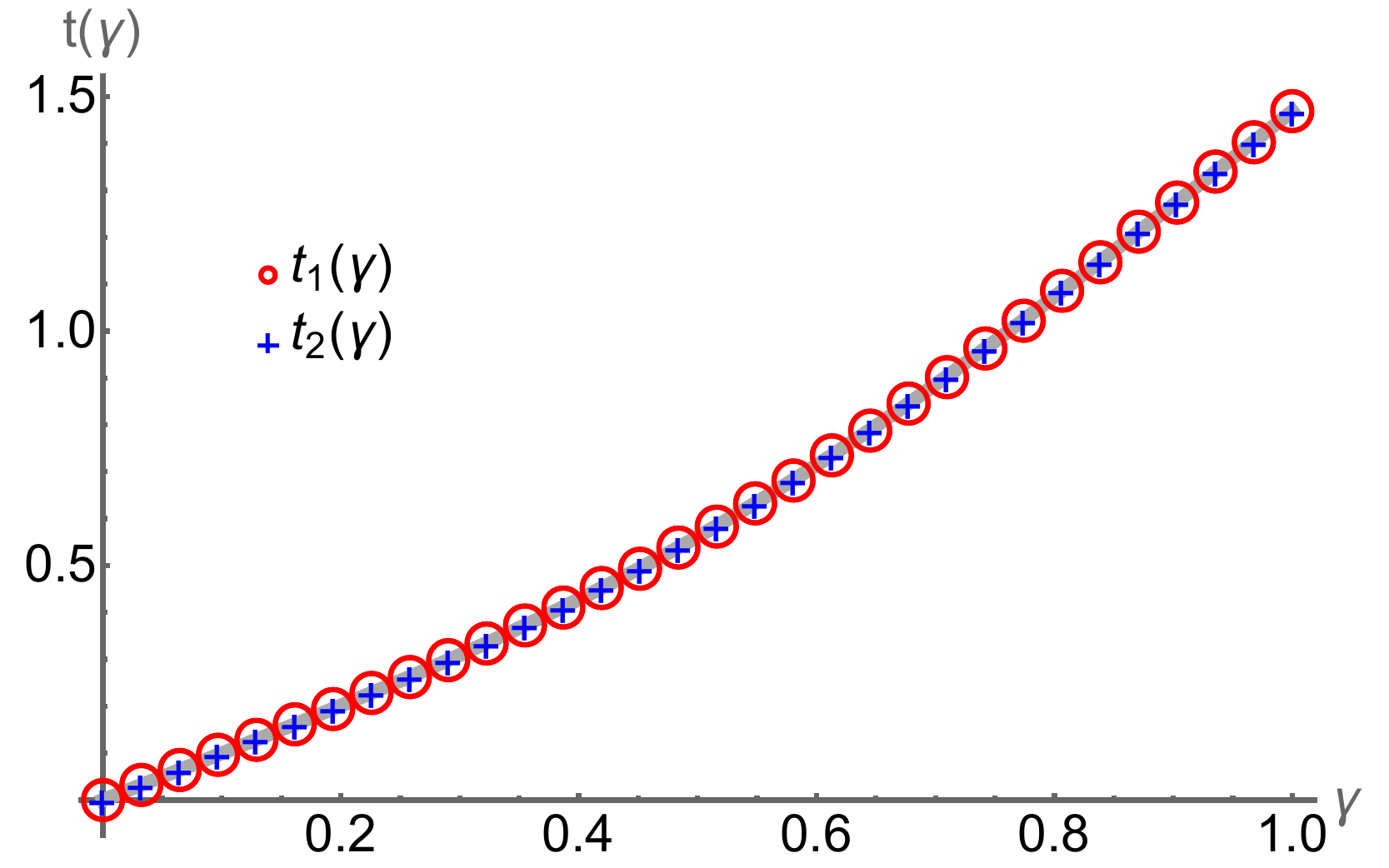}
    \includegraphics[scale=0.3]{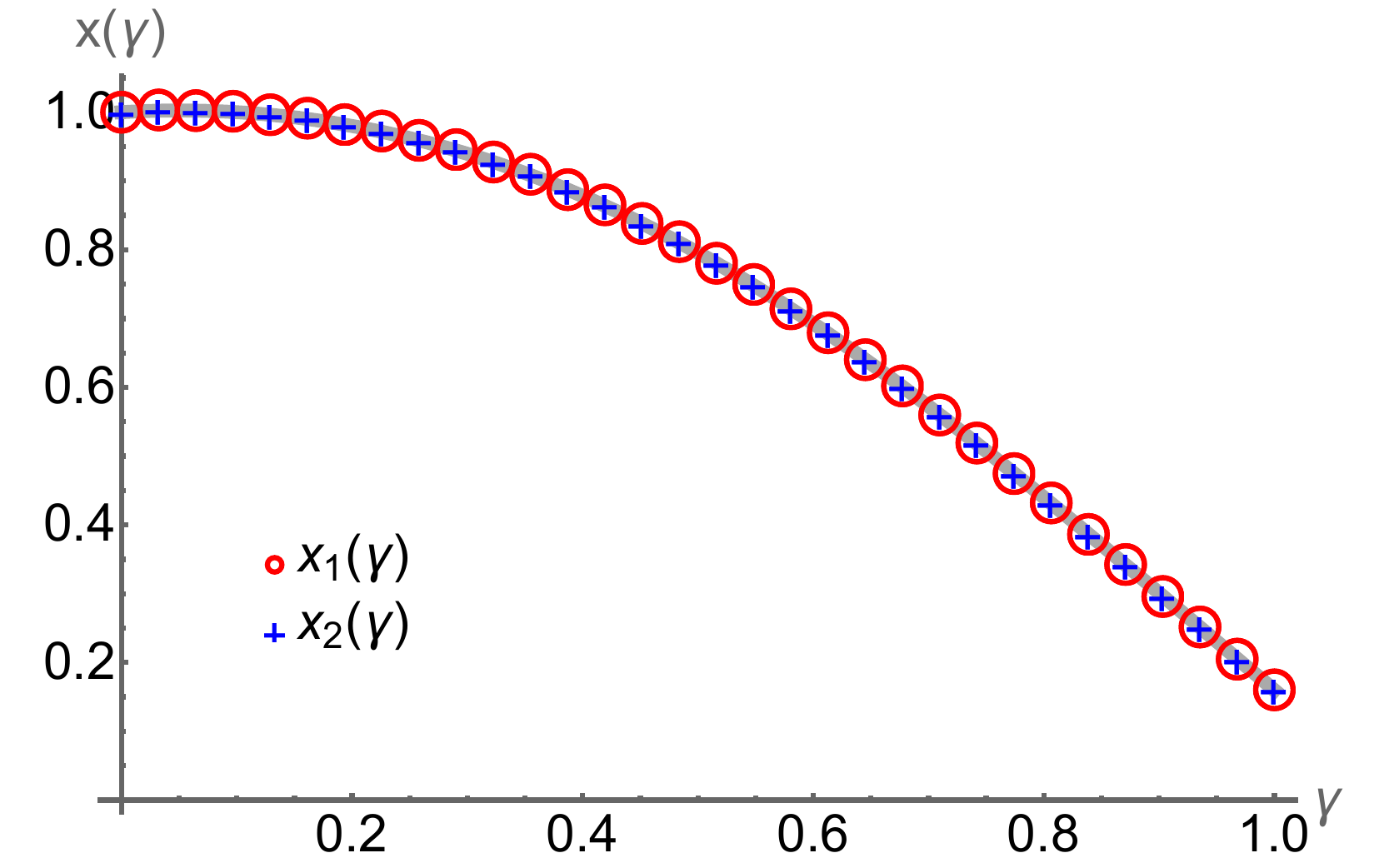}
    \caption{The values of the (top) time coordinates ${\bf t}_1(\gamma_i)$ (red circles) and ${\bf t}_2(\gamma_i)$ (blue crosses) and that of the (bottom) spatial coordinates ${\bf x}_1(\gamma_i)$ (red circles) and ${\bf x}_2(\gamma_i)$ (blue crosses) along the world-line parameter $\gamma$, as obtained from the critical point of $\mathds{E}^{\rm qrt}_{\rm IVP}$ with $V(x)=\kappa x^4$, discretized with $N_\gamma=32$ points and the \texttt{SBP21} operators. The solution of the corresponding geodesic equations via Mathematica's \texttt{NDSolve} is shown as light gray solid line.}
    \label{fig:NLdofofgamma}
\end{figure}

As in the previous subsection we discretize the world-line of the particle motion between $\gamma_i=0$ and $\gamma_f=1$, set the starting time to $t_i=0$ and the starting position to $x_i=1$. For our choice of $v_i=1/10$ we again decide on $\dot t=1$ and $\dot x=v_i$. The discretized action functional thus reads
\begin{align}
\nonumber \mathds{E}^{\rm qrt}_{\rm IVP}= &\frac{1}{2} \left\{ ({\bar{\mathds{D}}^{\rm R}}_t{\bf t}_1)^{\rm T}\mathbb{d}\left[1+2 \kappa {\bf x}^4_1\right] {\bar{\mathds{H}}} ({\bar{\mathds{D}}^{\rm R}}_t{\bf t}_1) -  ({\bar{\mathds{D}}^{\rm R}}_x{\bf x}_1)^{\rm T} {\bar{\mathds{H}}} ({\bar{\mathds{D}}^{\rm R}}_x{\bf x}_1)\right\}\\
\nonumber-&\frac{1}{2} \left\{ ({\bar{\mathds{D}}^{\rm R}}_t{\bf t}_2)^{\rm T} \mathbb{d}\left[1+2 \kappa {\bf x}^4_2\right] {\bar{\mathds{H}}} ({\bar{\mathds{D}}^{\rm R}}_t{\bf t}_2) -  ({\bar{\mathds{D}}^{\rm R}}_x{\bf x}_2)^{\rm T} {\bar{\mathds{H}}} ({\bar{\mathds{D}}^{\rm R}}_x{\bf x}_2)\right\}\\
\nonumber+&\lambda_1\big({\bf t}_1[1]-t_i\big)+\lambda_2\big((\mathds{D}{\bf t}_1)[1]-\dot t_i\big)\\
\nonumber+&\lambda_3\big({\bf x}_1[1]-x_i\big)+\lambda_4\big((\mathds{D}{\bf x}_1)[1]-\dot x_i\big)\\
\nonumber +&\lambda_5\big({\bf t}_1[N_\gamma]-{\bf t}_2[N_\gamma]\big) + \lambda_6\big({\bf x}_1[N_\gamma]-{\bf x}_2[N_\gamma]\big)\\
+&\lambda_7\big( (\mathds{D}{\bf t}_1)[N_\gamma]- (\mathds{D}{\bf t}_2)[N_\gamma]\big)+\lambda_8\big( (\mathds{D}{\bf x}_1)[N_\gamma]- (\mathds{D}{\bf x}_2)[N_\gamma]\big)\label{eq:discEIVPlin}
\end{align}
and taking the fourth power of the ${\bf x}_{1,2}$ vector is to be understood in an element wise fashion.

While for the linear potential, the time geodesic appeared to depend almost linearly on $\gamma$, we find that here a distinct curvature along $\gamma$ emerges, as shown in the top panel of \cref{fig:NLdofofgamma}. We plot the values of ${\bf t}_1(\gamma_i)$ as red circles and ${\bf t}_2(\gamma_i)$ as blue crosses and show as gray solid line the solution of the corresponding geodesic equation, obtained from the \texttt{LSODA} algorithm of Mathematica's \texttt{NDSolve} command. Again the physical limit of equal values ${\bf t}_1(\gamma)={\bf t}_2(\gamma)$ is realized.

The values of the spatial coordinate ${\bf x}_1(\gamma_i)$ and ${\bf x}_2(\gamma_i)$ as obtained from the critical point of $\mathds{E}^{\rm qrt}_{\rm IVP}$ with $V(x)=\kappa x^4$ are plotted in the bottom panel of \cref{fig:NLdofofgamma} with the direct numerical solution of the geodesic equation added as gray solid line.

\begin{figure}
    \includegraphics[scale=0.3]{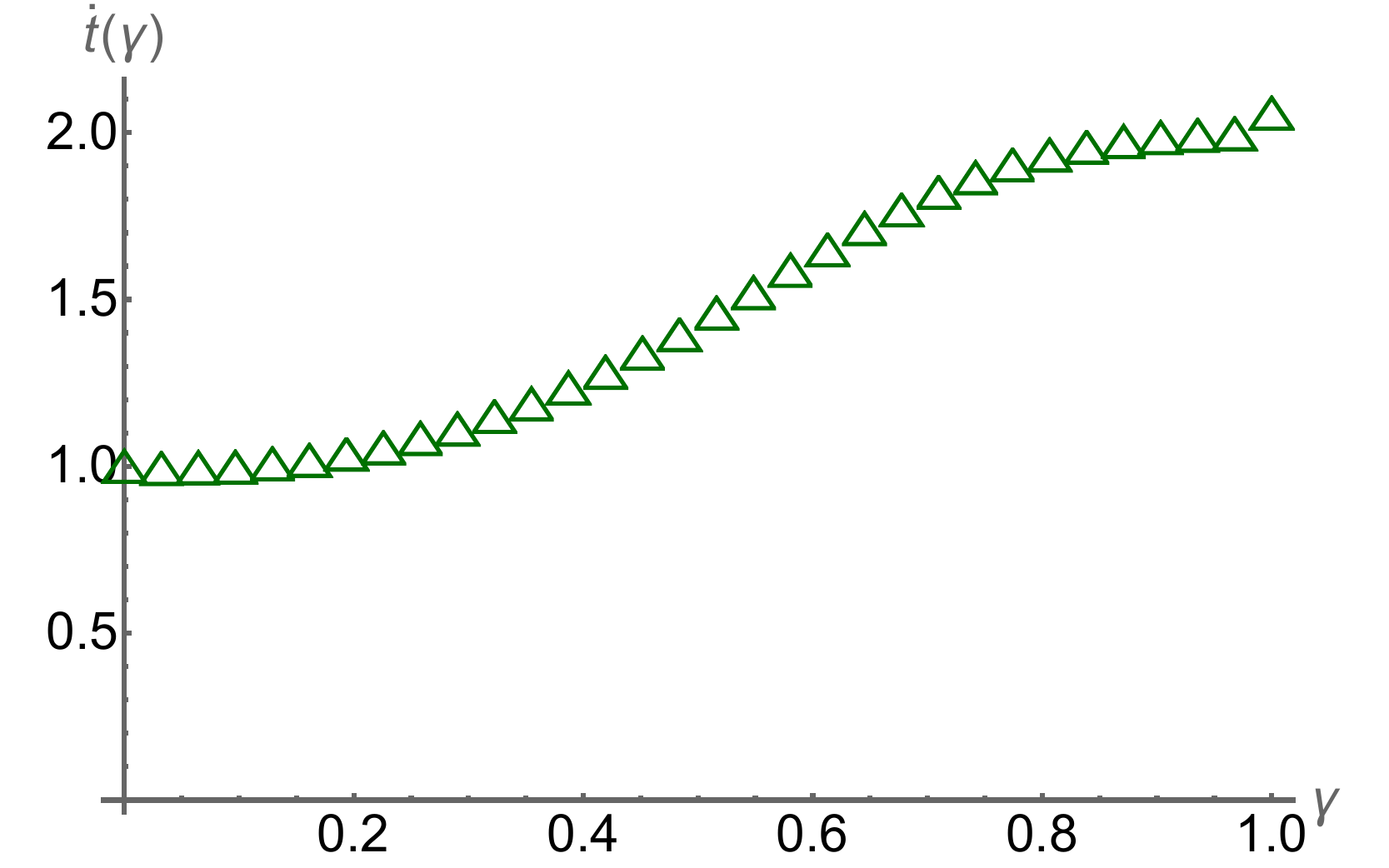}
    \caption{The change of time with respect to the world-line coordinate $\gamma$, as obtained from the critical point of $\mathds{E}^{\rm qrt}_{\rm IVP}$ with $V(x)=\kappa x^4$.}
    \label{fig:NLtdotofgamma}
\end{figure}

Note that even though we have provided an initial velocity of the time along $\gamma$ again with value $\dot t_i=1$, the final time reached by the simulation now lies at ${\bf t}[N_\gamma]=1.47$. Similarly one finds that that a dynamical discretization in $t$ emerges, which, as shown in \cref{fig:NLtdotofgamma}, varies from the initial values $\dot t_i=1$ to $(\mathds{D}{\bf t})[N_\gamma]=2.06$. This behavior can be understood when realizing that the trajectory $x(t)$ in the non-linear case shows a stronger curvature close to $t=0$ than at later times. I.e. we find again that the automatically generated non-trivial mesh (through automatic AMR) for the time coordinate adapts to the dynamics, by exhibiting a finer spacing at initial times.

\begin{figure}
    \includegraphics[scale=0.3]{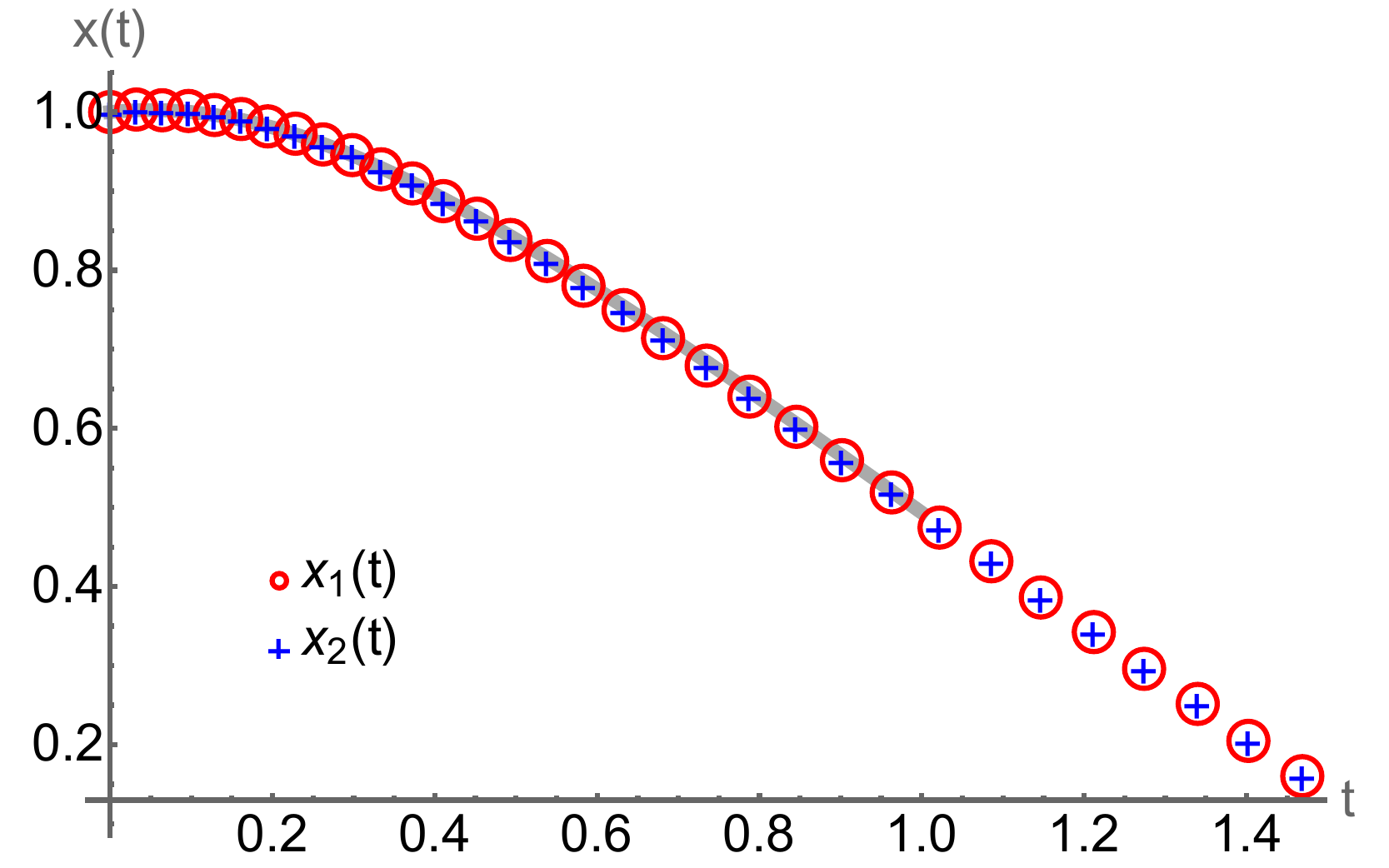}
    \caption{The forward ${\bf x}_1(t_1)$ (red circles) and backward ${\bf x}_2(t_2)$ (blue crosses) degrees of freedom obtained from the critical point of $\mathds{E}^{\rm qrt}_{\rm IVP}$ with $V(x)=\kappa x^4$. We discretize with $N_\gamma=32$ points and the \texttt{SBP21} operators. The trajectory from a \texttt{LSODA} solver of the corresponding equation of motion is given as gray solid line.}
    \label{fig:NLxoft}
\end{figure}

Let us take a look at the results from our geometrized formalism as physical trajectory in \cref{fig:NLxoft}, i.e. plotted as ${\bf x}_{1,2}(t_{1,2})$ (red circles and blue crosses). They are compared to the solution of the conventional equation of motion, obtained from treating time as independent variable $d^2x/dt^2 = -(4\kappa x^3 )(1-(dx/dt)^2)^{(3/2)}$, computed via the \texttt{LSODA} algorithm of Mathematica's \texttt{NDSolve} command (gray solid line) in the range $t\in[0,1]$. We find that within this range the solution from our geometrized discrete approach shows excellent agreement. Note that due to the non-equidistant emergent time discretization, the physical trajectory $x(t)$, shown in \cref{fig:NLxoft} extends beyond the point $t=1$.


As for the linear potential, let us investigate quantitatively the properties of the trajectories ${\bf t}(\gamma_i)$ and ${\bf x}(\gamma_i)$ by inserting them into the naively discretized geodesic equations. For the quartic potential, the continuum geodesic equations for the temporal and spatial coordinate read
\begin{align}
 &\frac{d}{d\gamma}\Big(g_{00}\frac{dt}{d\gamma}\Big)=\frac{d}{d\gamma}\Big( \big(1+2 \kappa x^4\big) \frac{dt}{d\gamma}\Big) 
 =0,\\
 &\frac{d}{d\gamma}\Big(\frac{dx}{d\gamma}\Big) + \frac{1}{2}\frac{\partial g_{00}}{\partial x}\Big( \frac{dt}{d\gamma} \Big)^2=\frac{d^2x}{d\gamma^2} + 4\kappa x^3\Big( \frac{dt}{d\gamma} \Big)^2=0.
\end{align}
Naively discretizing these equations by replacing derivatives with SBP operators leads to the following discrete geodesic equations
\begin{align}
    &\mathds{D}\big( (1+2\kappa {\bf x}^4 )\circ\mathds{D}{\bf t} \big)=\Delta{\bf G}^t,\\
    &\mathds{D}\mathds{D}{\bf x} + (4\kappa {\bf x}^3) \circ (\mathds{D}{\bf t})\circ (\mathds{D}{\bf t})=\Delta{\bf G}^x.
\end{align}
where again taking a power of the ${\bf x}_{1,2}$ vector is to be understood in an element wise fashion. To evaluate how well the solution obtained from the critical point of $ \mathds{E}^{\rm qrt}_{\rm IVP}$ fulfills the naive discretized geodesic equations we have again introduced the quantities $\Delta {\bf G}^t$ and $\Delta {\bf G}^x$ above.

\begin{figure}
    \includegraphics[scale=0.3]{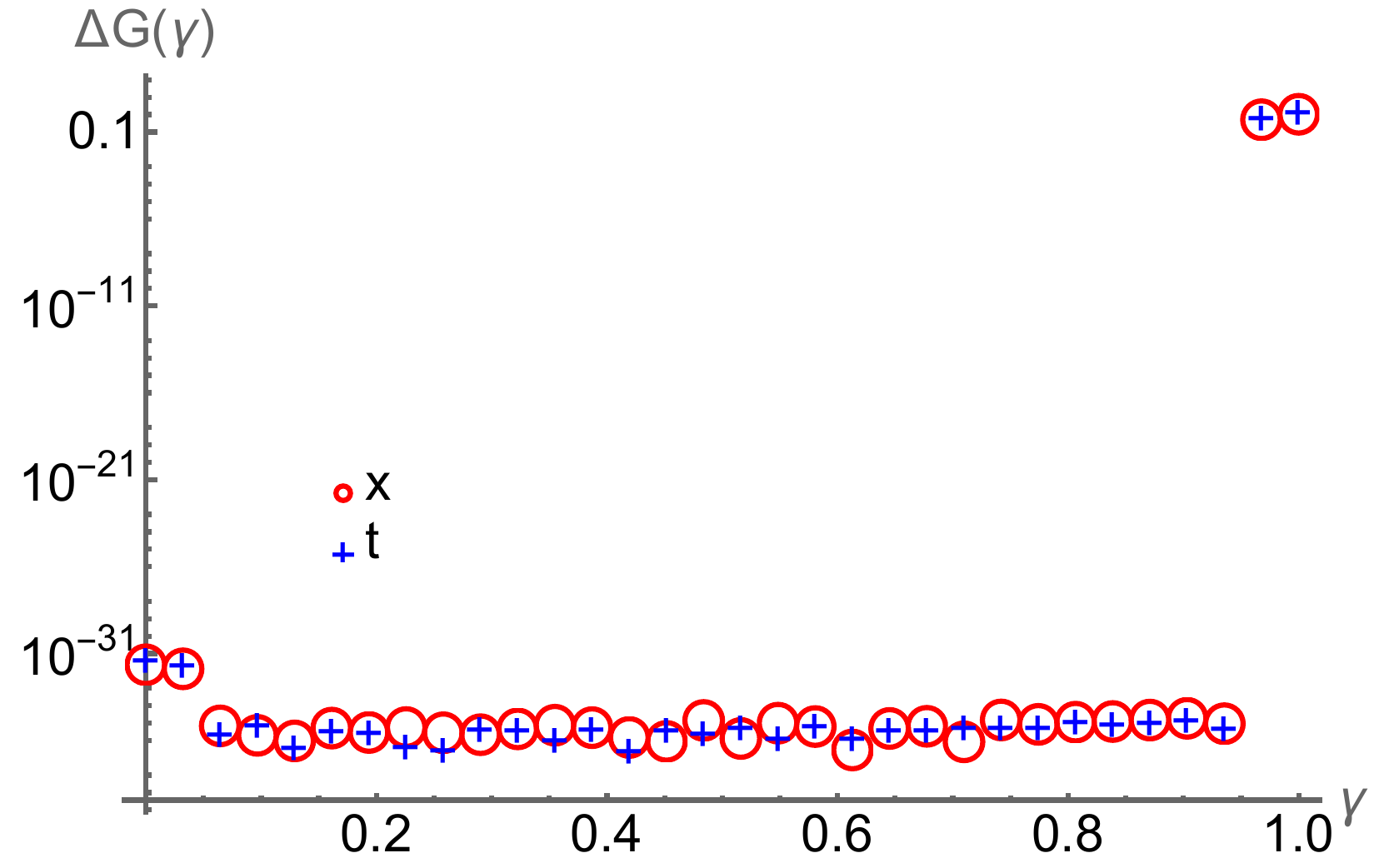}
    \includegraphics[scale=0.3]{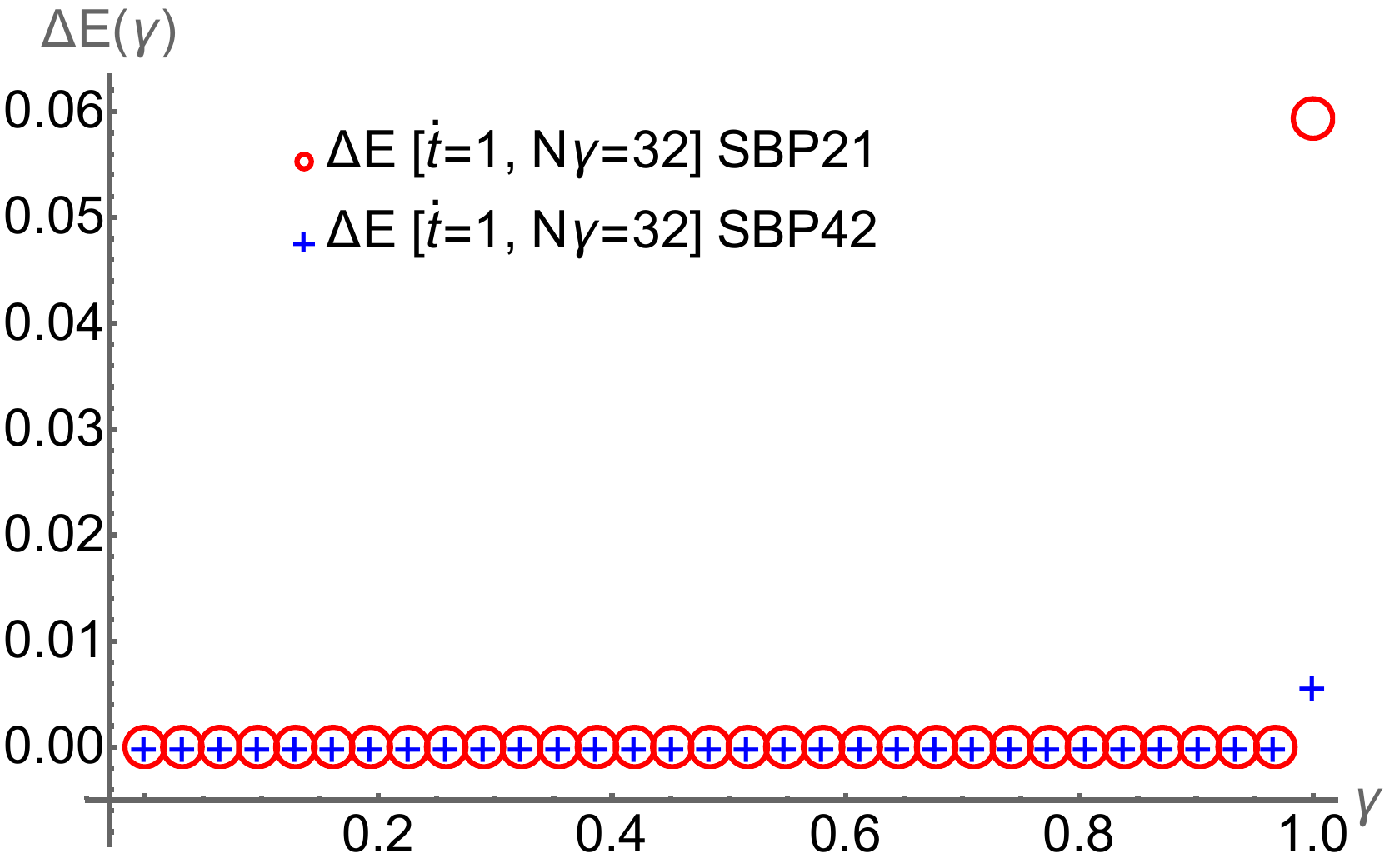}
    \caption{(top) Deviation of the spatial coordinate $\Delta{\bf G}^x$ (red circles) and time coordinate $\Delta{\bf G}^t$ (blue crosses) from the discretized geodesic equation, as obtained from the critical point of $\mathds{E}^{\rm qrt}_{\rm IVP}$ with $V(x)=\kappa x^4$ at $N_\gamma=32$ discretized with the \texttt{SBP21} operator. (bottom) Deviation $\Delta {\bf E}$ of the quantity ${\mathds Q}_t$ from its continuum value as given by the initial conditions. Results for the \texttt{SBP21} operator are given as red circles and those for \texttt{SBP42} as blue crosses. Note that the deviation $\Delta{\bf E}$ in the interior is exactly zero within machine precision.}
    \label{fig:NLDeltaGDeltaE}
\end{figure}

As shown in \cref{fig:NLDeltaGDeltaE} also here in the highly non-linear scenario, we find that the values of both ${\bf x}$ (red circles) and ${\bf t}$ (blue crosses) follow the discretized geodesic equations excellently, except for the last two points.

The most important question however remains whether in the non-linear discretized system, the continuum quantity $Q_t$ from \cref{eq:Qt} also remains conserved. Its naively discretized counterpart here reads
\begin{align}
    {\bf \mathds{Q}_t}=(\mathds{D}{\bf t})\circ({\bf 1} + 2\kappa {\bf x}^4),
\end{align}
and we define its deviation from the continuum result via the difference
\begin{align}
    \Delta {\bf E} = {\bf \mathds{Q}_t} -Q_t = (\mathds{D}{\bf t})\circ({\bf 1} + 2 \kappa {\bf x}^4) - \dot t_i (1+2\kappa x_i^4),
\end{align}
which we plot in the bottom panel of \cref{fig:NLDeltaGDeltaE} using the \texttt{SBP21} operator (red circles) and the \texttt{SBP42} operator (blue crosses).

We find also in the case of a non-linear potential that ${\bf \mathds{Q}_t}$ is preserved exactly in the interior of the simulation time domain. Up to machine precision its values in the interior also take on the correct continuum value. Similar to what we saw in the linear case, the last point deviates from the continuum value. It is reassuring to see that the absolute deviation at $\gamma_f$ reduces already by an order of magnitude when going from a \texttt{SBP21} to an \texttt{SBP42} operator.

\begin{figure}
    \includegraphics[scale=0.3]{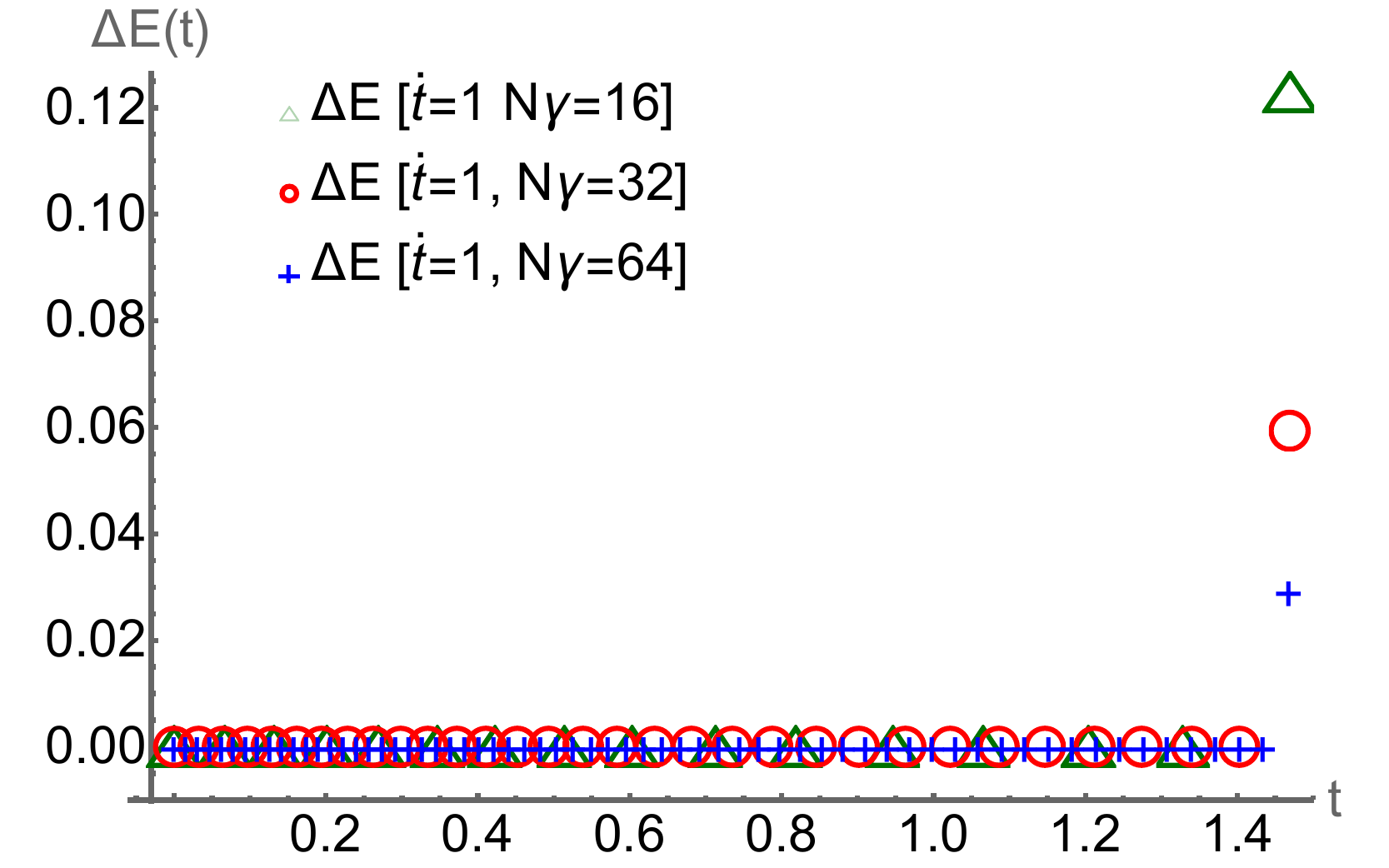}
    \caption{Deviation of the quantity $\mathds{Q}_t$ from its continuum value $Q_t$ in case of a fixed final time and grid refinement in $N_\gamma$. Note that the magnitude of the deviation at the final point monotonously decreases with grid refinement.}
    \label{fig:NLDeltaGDeltaE2}
\end{figure}

One may now ask whether the deviation of $\Delta {\bf E}$ from its continuum value at $\gamma_f$ is in some way related to the fact that we use $N_\gamma=32$ points to discretize the world-line parameter. The answer is negative, as demonstrated in \cref{fig:NLDeltaGDeltaE2}.
Three different datasets are shown in \cref{fig:NLDeltaGDeltaE2}, where for fixed $\dot t_i$ the grid spacing in $\gamma$ is changed. The green triangles denote the results for $\Delta E$ when using $N_\gamma=16$, the red circles $N_\gamma=32$ and the blue crosses $N_\gamma=64$. We have confirmed explicitly that in all cases the values of $Q_t$ are preserved up to machine precision in the interior of the simulated time domain. It is indeed only the last point that shows a deviation and we see that the absolute magnitude of the deviation reduces as the grid is refined. 

\begin{figure}
    \includegraphics[scale=0.3]{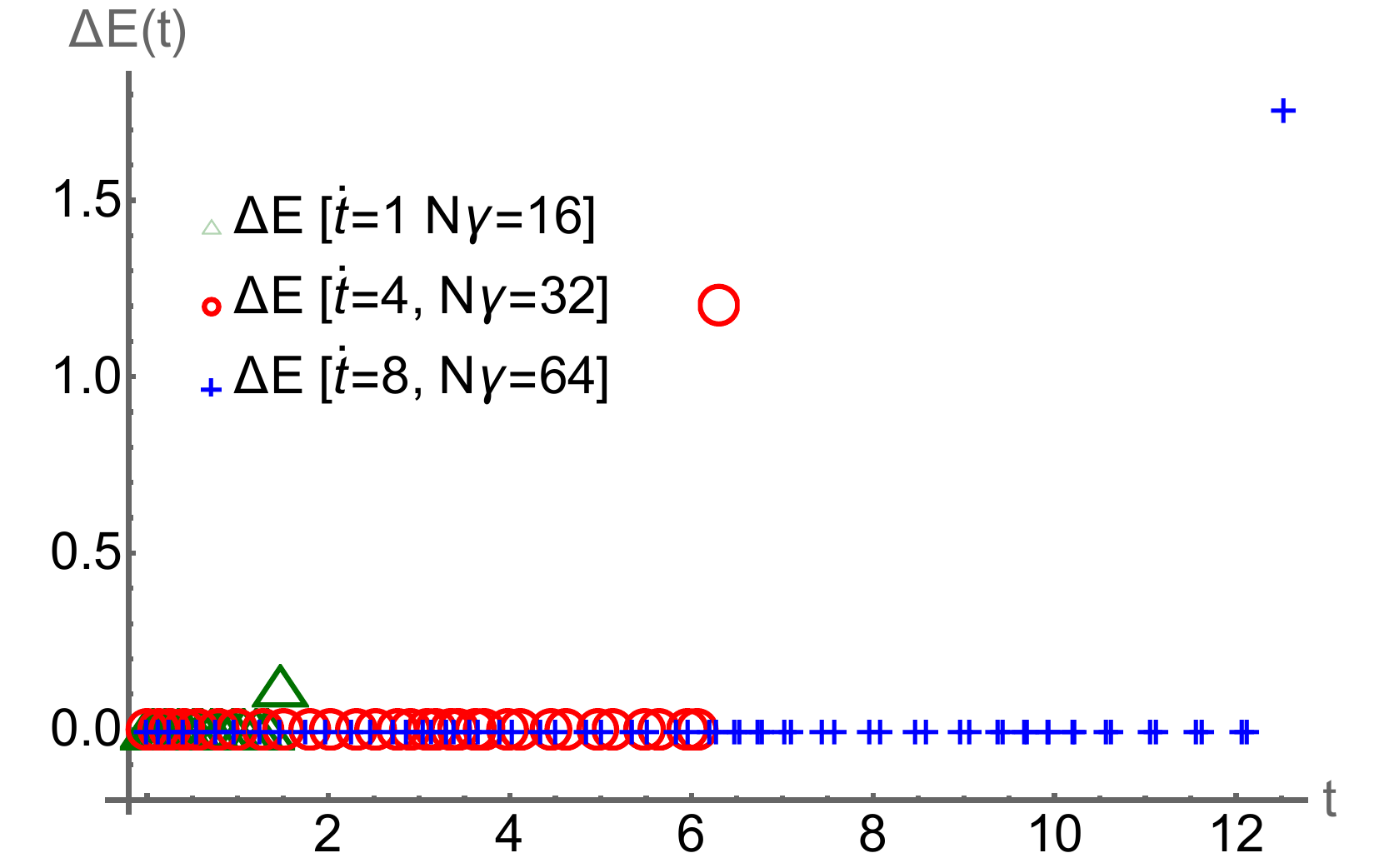}
    \includegraphics[scale=0.3]{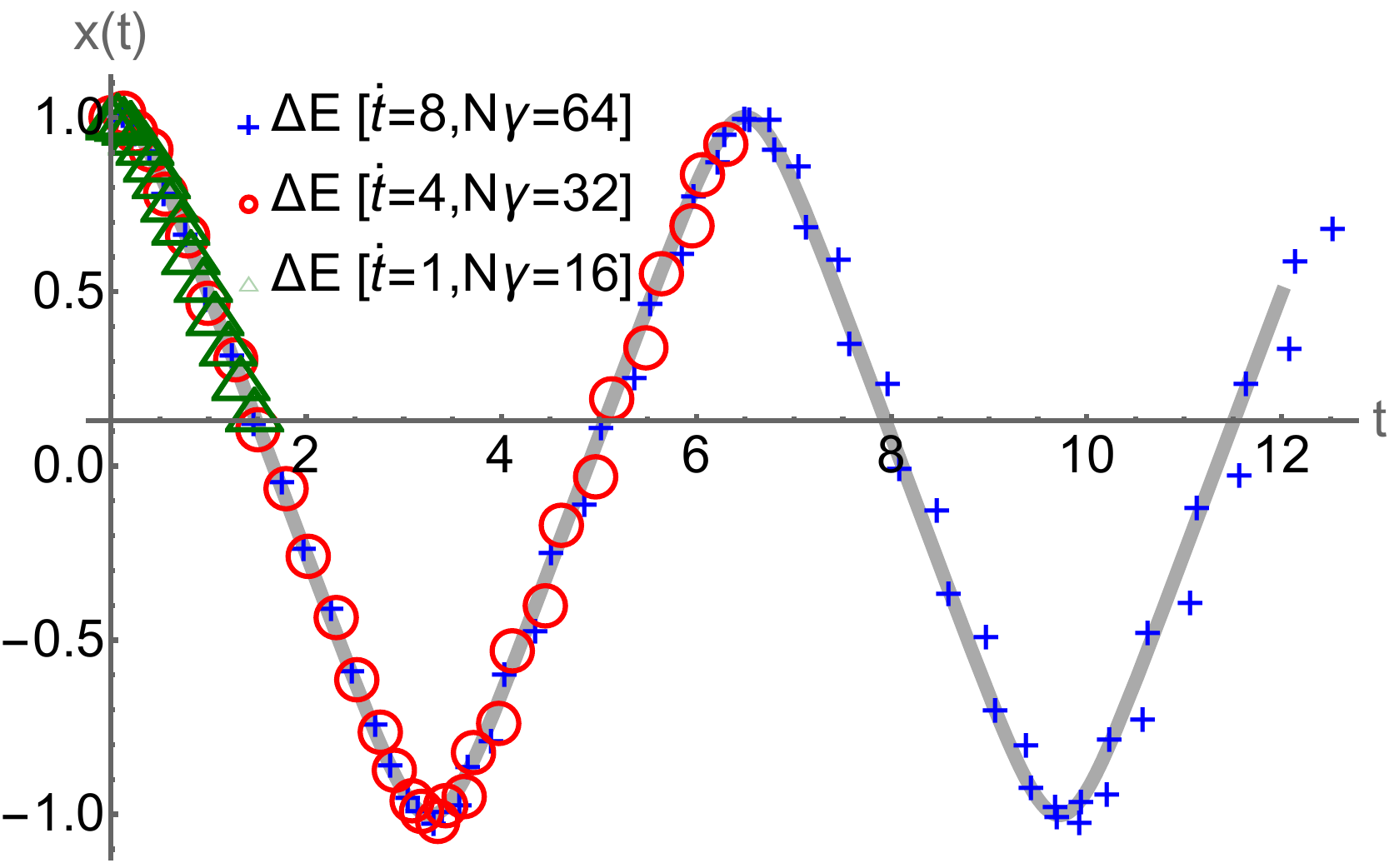}
    \caption{(top) Deviation of the quantity $\mathds{Q}_t$ from its continuum value $Q_t$ when moving the latest time to larger values by specifying different initial $\dot t_i$. (bottom) The corresponding simulated trajectories ${\bf x}(t)$.}
    \label{fig:NLDeltaGDeltaE3}
\end{figure}

For the next test, we instead increase $N_\gamma$ together with $\dot t_i$ to let the simulation proceed to larger values of time $t$. In the top panel of \cref{fig:NLDeltaGDeltaE3} we plot the deviation of $Q_t$ from its continuum value for three choices $\dot t_i=1, N_\gamma=16$ (green triangles), $\dot t_i=4, N_\gamma=32$ (red circles) and $\dot t_i=8,N_\gamma=64$  (blue crosses). As seen before in the interior of the simulated time domain, the values of $\mathds{Q}_t$ remain exactly preserved and only the last point deviates. We find that the magnitude of the deviation in the last point changes only marginally with the length of the simulated trajectory. For completeness the corresponding trajectories $x(t)$ are plotted in the bottom panel of \cref{fig:NLDeltaGDeltaE3}. Again let us emphasize that, as we will show below, the presence of this single deviating point does not spoil the convergence to the correct solution under grid refinement.

The exact conservation of the quantity ${\mathds Q}_t$ in the interior is remarkable, as e.g. the trajectory in the bottom panel of \cref{fig:NLDeltaGDeltaE3} for $\dot t_i=8,N_\gamma=64$ shows sizable discretization artifacts (which disappear under grid refinement). We believe that it is due to the manifest time-translation invariance of the underlying action functional that the combined dynamics of $x(\gamma)$ and $t(\gamma)$, including the automatically generated non-equidistant time mesh, achieve conservation of the continuum quantity. 

The fact that the solutions we obtain fulfill the naively discretized geodesic equations and provide exact conservation of the continuum conserved charge in the interior of the simulated domain (see \cref{fig:NLDeltaGDeltaE}) bodes well for establishing its stability. Since in the IVP setting $t(\gamma_f)$ is not given but emerges dynamically we cannot directly apply \cref{eq:boundonsol} as proof of stability. However, as long as we can assume that the simulated time range (given a certain $\dot t(\gamma_i)$ is finite, the linear bound of \cref{eq:boundonsol} on the norm ${\cal H}_{\rm BVP}$ holds in the discrete setting. In turn we deduce that the solution cannot exhibit stronger than linear rise of the derivatives of either $t(\gamma)$ or $x(\gamma)$, implying stability of the approach.

\begin{figure}
    \includegraphics[scale=0.3]{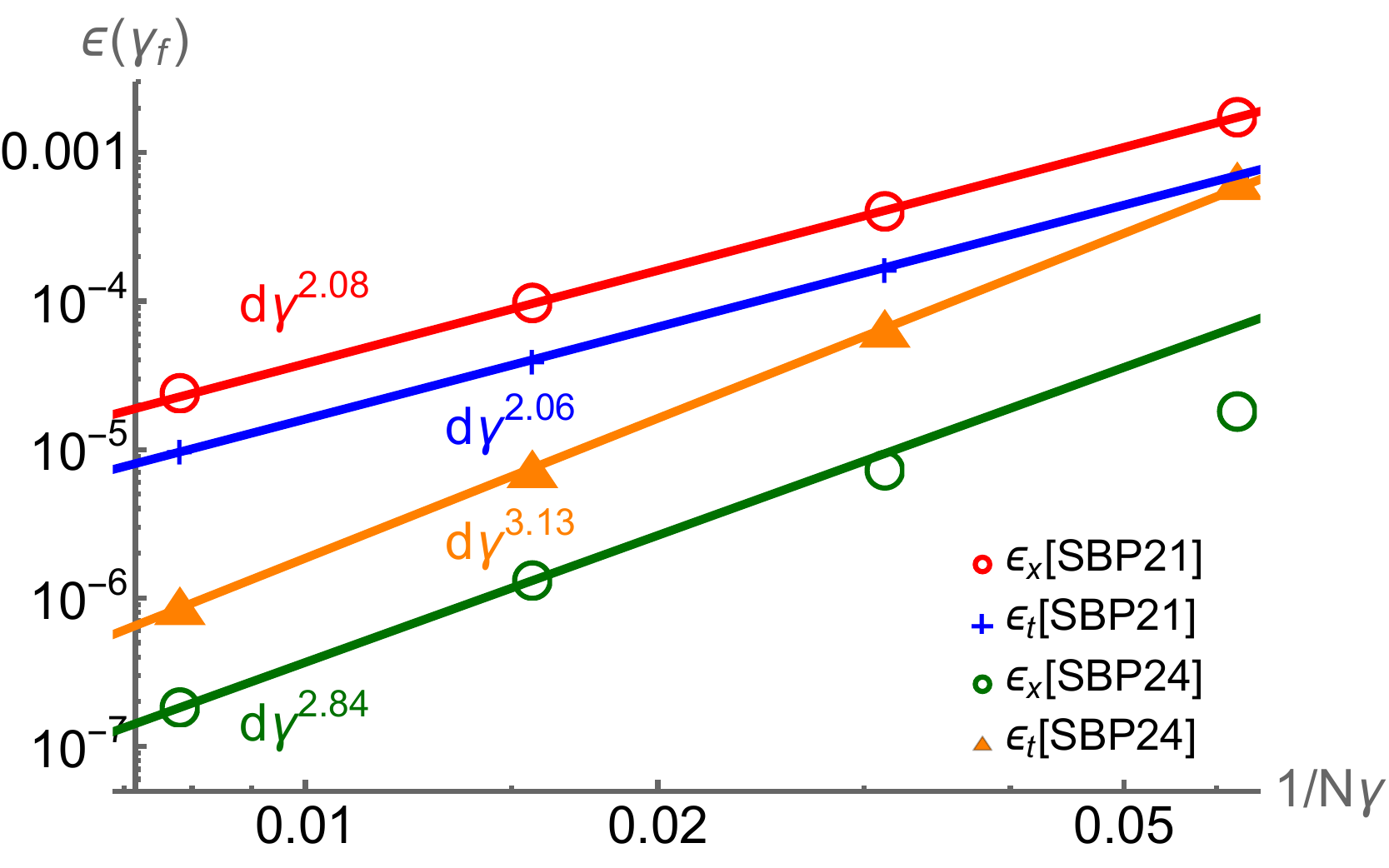}
    \includegraphics[scale=0.3]{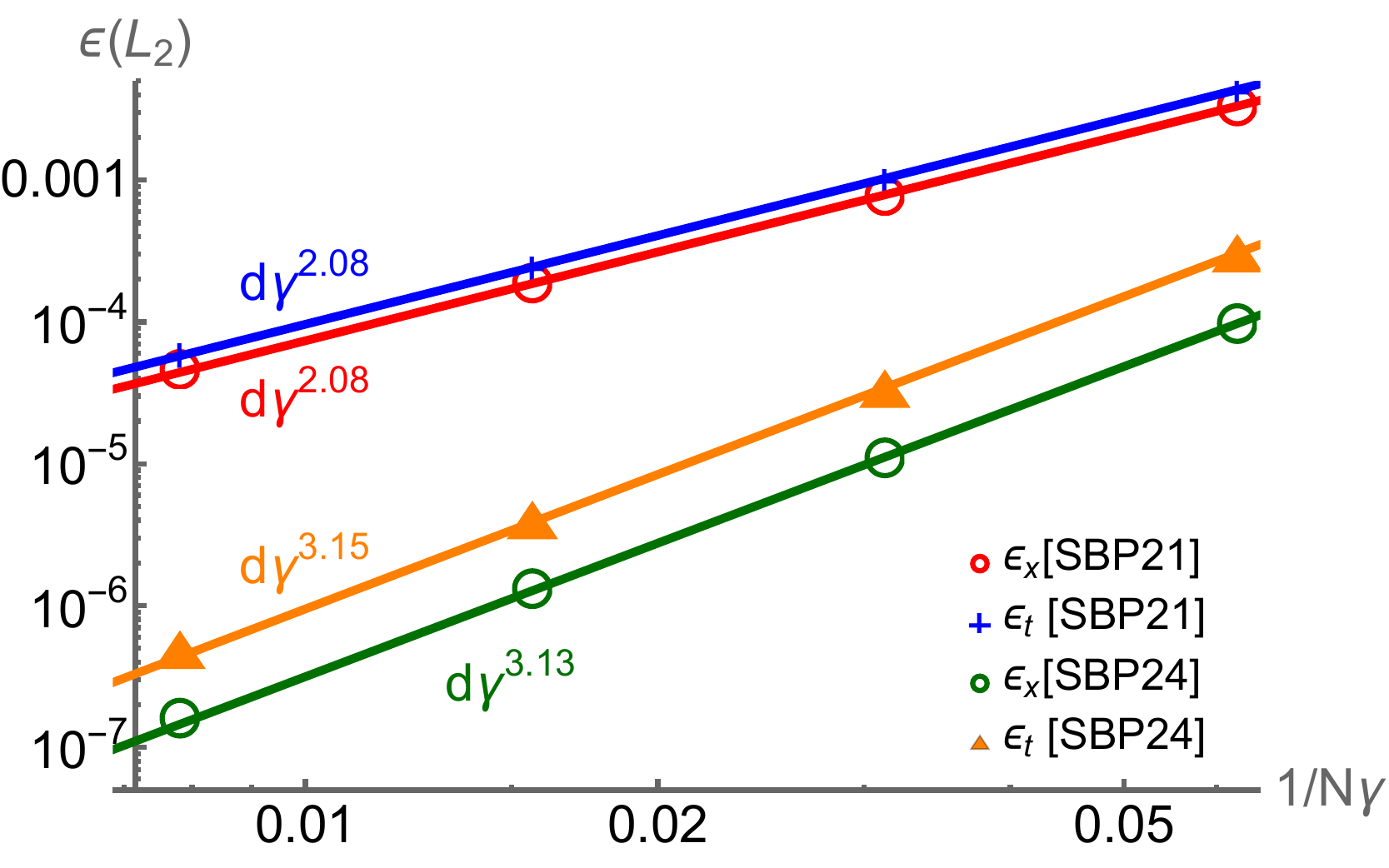}
 \caption{Visualization of the convergence of our novel variational approach under the refinement of the $\gamma$ grid, based on the final point via $\epsilon(\gamma_f)$ (top panel) and globally via $\epsilon(L_2)$ (bottom panel). In each plot we show the results from the \texttt{SBP21} (top two lines) and \texttt{SBP42} (bottom two lines) operators. Red circles and blue crosses denote the absolute deviation in the spatial coordinate $\epsilon_x$ and temporal coordinate $\epsilon_t$ between the continuum solution and the values obtained from the the critical point of $\mathds{E}^{\rm qrt}_{\rm IVP}$ with $V(x)=\kappa x^4$. The corresponding deviations for the \texttt{SBP42} operator are given by the green circles and orange triangles respectively.}
    \label{fig:NLConv}
\end{figure}

Let us now quantify the convergence properties of our variational approach using the results from the lowest order \texttt{SBP21} operator and those coming from the \texttt{SBP42} operator in \cref{fig:NLConv}.

As in the linear potential case, in the top panel of \cref{fig:NLConv}, we select the most disadvantageous points for our convergence study, i.e. we compare the deviation from the continuum geodesic equations $\epsilon(\gamma_f)_x=|{\bf x}[N_\gamma]-x_{\rm true}(1)|$ and $\epsilon(\gamma_f)_t=|{\bf t}[N_\gamma]-t_{\rm true}(1)|$
at $\gamma_f$, exactly where the deviation from the continuum result was maximal in the top panel of \cref{fig:NLDeltaGDeltaE}. Also in the non-linear scenario we find that under refinement of the $\gamma$ grid, the discrete solution monotonously approaches the true continuum values.

Taking the \texttt{SBP21} results, the best fit to $\epsilon(\gamma_f)_x$ reveals a scaling with $\Delta\gamma^{2.08}$, while for $\epsilon(\gamma_f)_t$ an virtually identical $\Delta\gamma^{2.06}$ is obtained. 

For \texttt{SBP42}, we find that the convergence is slightly worse than in the linear potential case. As seen in the green circles plotted in \cref{fig:NLConv}, the asymptotic convergence regime is reached for $32 <N_\gamma <64$. Once we are in that regime, we find that $\epsilon(\gamma_f)_x$ exhibits a scaling of $\Delta\gamma^{2.84}$, close to the expected value of three. On the other hand $\epsilon(\gamma_f)_t$ shows a consistent performance with a scaling of $\Delta\gamma^{3.13}$ already at $N_\gamma=32$.

Let us now investigate the global convergence in the bottom panel of \cref{fig:NLConv} using the $L_2$ norm $\epsilon(L_2)_x=\sqrt{({\bf x}-{\bf x}_{\rm true})^{\rm T}.\mathds{H}.({\bf x}-{\bf x}_{\rm true})}$ and correspondingly $\epsilon(L_2)_t=\sqrt{({\bf t}-{\bf t}_{\rm true})^{\rm T}.\mathds{H}.({\bf t}-{\bf t}_{\rm true})}$, where ${\bf x}_{\rm true}$ and ${\bf t}_{\rm true}$ are taken from the numerical solution of the geodesic equations, used for comparison in \cref{fig:NLdofofgamma}.

Reassuringly we find that the global convergence properties of our approach are better than indicated by those of the most disadvantaged point in the top panel of \cref{fig:NLConv}. Indeed we find that for the \texttt{SBP41} operators, the global scaling regime is reached already at $N_\gamma=32$, similarly to the \texttt{SBP21} case. In addition, the global convergence rate $\Delta\gamma^\beta$ for \texttt{SBP42} operators lies consistently above $\beta\ge3$ for both the $x$ and $t$ degrees of freedom.

Again, these convergence result are in good agreement with those of our previous study \cite{Rothkopf:2022zfb}, where the standard action functional was discretized with time as independent parameter. 

\section{Summary and Outlook}
\label{sec:conc}

In this study we have put forward a novel geometric variational approach for solving a large class of initial value problems, associated with the dynamics of point particles evolving under a generic $x$ dependent potential $V(x)$. Taking inspiration from the general theory of relativity, we consider both time and spatial coordinates of the point particle as dependent variables of a world-line parameter $\gamma$. We select a continuum action functional, which in the non-relativistic limit reduces to the standard action of point mechanics and whose critical point encodes a set of geodesic equations for $x(\gamma)$ and $t(\gamma)$. After doubling the degrees of freedom $t_{1,2}$ and $x_{1,2}$ we can relate the critical point of the corresponding doubled d.o.f. action with the classical trajectory. Using the concept of Killing vectors we identify conserved quantities, e.g. related to the continuum time translation invariance of the action.

Deploying the regularized SBP operators originally introduced in \cite{Rothkopf:2022zfb}, we discretize the continuum action and add Lagrange multipliers to enforce the initial and connecting conditions between the doubled ${\bf t}_{1,2}$ and ${\bf x}_{1,2}$. The main novelty of our approach is that the discretized action \textit{retains the continuum symmetries}, in particular the invariance under time translations. Exactly mimicking integration by part through the use of SBP finite difference operators entails that the derivation of the conserved charges associated with the Killing vectors of the system is also exactly mimicked in the discrete setting. I.e. the continuum conserved quantities $Q_K$ retain their role even after discretization.

The numerical results we obtain for both a linear and highly non-linear potential show that a discretization of time $t$ now indeed emerges dynamically, adapting to the behavior of the spatial coordinate $x$. This is a concrete realization of an automatically generated non-equidistant mesh for the time coordinate, guided by our action functional with manifest continuum translation symmetry, i.e. an automatic AMR procedure. We have shown that except for the last two points along the discrete $\gamma$, the solution we obtain follows the naively discretized geodesic equations excellently. 

Even more importantly, the naively discretized counterpart $\mathds{Q}_t$ of the continuum conserved quantity $Q_t$ \textit{remains exactly preserved} in the interior of the simulated time domain, where it even \textit{retains its continuum value} exactly within machine precision. A small deviation from the values in the interior for $\mathds{Q}_t$ is observed at the last step $\gamma_f$. This deviation however decreases both under grid refinement, as well as when increasing the order of the SBP operator.

Point-wise, as well as global scaling analyses under grid refinement show that even in the presence of two points deviating from the naively discretized geodesic equations at the last two $\gamma$ steps, the solution monotonously improves and manages to approach the true solution. When deploying the \texttt{SBP21} operator, we achieve consistent scaling in $\Delta \gamma^\beta$ with $\beta\gtrsim 2$ for both the linear and non-linear potential. For \texttt{SBP42} in case of a linear potential the dependence on the grid spacing follows the  expected power law $\Delta \gamma^\beta$ with  $\beta\gtrsim3$ for all values of $N_\gamma$ we inspected. For the non-linear potential, the scaling regime for point-wise convergence at the last point $\gamma_f$ is reached with \texttt{SBP42} for $32<N_\gamma<64$ with a slightly worse scaling of $2.84 \le \beta \le 3.13$. Global convergence on the other hand shows consistent scaling at all $N_\gamma$ we considered, with exponents $\beta\ge3$, in agreement with the findings in our previous paper \cite{Rothkopf:2022zfb}, where the standard action functional was discretized with time as independent variable.

This study presents a proof of principle that initial value problems can be discretized, while retaining continuum symmetries. Three future directions will be explored:  we may ask how we can capture systems of ordinary differential equations that e.g. contain a term that is proportional to a first derivative in $x$ with respect to time? To this end we must exploit the versatility of the doubled d.o.f. approach more thoroughly. Furthermore we will explore how the reparametrization invariant formulation can be applied to partial differential equations in higher dimensions, taking insight from how the non-relativistic action emerges from our relativistic starting point in \cref{eq:eq6}. In addition, to better understand the origin of the single deviating value in the otherwise exactly preserved $\mathds{Q}_t$, we will develop a genuinely weak formulation of our approach, devoid of Langrange multipliers for enforcing initial and connecting conditions.

We believe that the quest for retention of defining continuum properties in discretized systems is both conceptually and practically valuable. Not only does the preservation of symmetries place powerful physical constraints on the solution but in addition offers a mechanism for the automatic generation of optimal discrete spacetime grids to ensure conservation of the Noether charges associated with these symmetries. We hope that this study provides the community with a novel impulse in this direction. 

\section*{Acknowledgements}
A.~R. thanks Will Horowitz for inspiring and insightful discussions and Alex Nielsen for valuable insight on the general theory of relativity. A.~R.~ gladly acknowledges support by the Research Council of Norway under the FRIPRO Young Research Talent grant 286883. J.~N. was supported by the Swedish Research Council grant nr. 2021-05484. The study has benefited from computing resources provided by  
UNINETT Sigma2 - the National Infrastructure for High Performance Computing and Data Storage in Norway under project NN9578K-QCDrtX "Real-time dynamics of nuclear matter under extreme conditions" 


\appendix    
\section{Regularized SBP operators in affine coordinates}
\label{sec:appregSBP}

We here briefly review the idea and some technical aspects of constructing null-space consistent regularized SBP operators using affine coordinates, developed in our study \cite{Rothkopf:2022zfb}.

The goal of regularizing conventional SBP operators $\mathds{D}$, such as those defined e.g. in \cref{eq:SBP21} and \cref{eq:SBP42}, lies in removing their unphysical zero modes. These may appear as highly oscillatory eigenfunctions to $\mathds{D}^{\rm T}$ with zero eigenvalue. To this end we take inspiration from regularization techniques developed for partial differential equations. There the concept of null-space consistent SBP operators has been discussed in detail (see e.g. \cite{svard2019convergence,linders_properties_2020,svard_convergence_2021,ranocha2021new}).

For a differential equation, the boundary conditions may be enforced in the weak sense by adding a simultaneous approximation penalty term (SAT) \cite{carpenter1994time}, which can be partially absorbed into the finite difference operator, lifting its zero modes. Take for example a simple discretized first order differential equation
\begin{align}
\mathds{D} {\bf u} = \lambda {\bf u} + \sigma_0 \mathds{H}^{-1} \mathds{E}_0\big( {\bf u} - {\bf g}\big), \label{eq:simpleODE}
\end{align}
where the SAT penalty term has been added to the right-hand side. It features the matrix $\mathds{E}_0={\rm diag}[1,0,\ldots,0]$ that makes reference only to the first entry in the discretized functions ${\bf u}$ and ${\bf g}$, the latter of which contains the initial value in its first entry ${\bf g}=(u_0,0,\ldots,0)$. The SAT term also contains $\mathds{H}^{-1}$, i.e. ${\Delta t}^{-1}$, which increases the strength of the penalty as ${\Delta t}\to0$.  The parameter $\sigma_0$ in the SBP-SAT approach is tuned to satisfy stability properties and its optimal value is found to be $\sigma_0=-1$ (see e.g. ref.~\cite{aalund2016provably,ALUND2019209,glaubitz2023summation}), a choice we adopt in the following. In the differential equation context one conventionally absorbs the term proportional to ${\bf u}$ into a new $\tilde{\mathds{D}}=\mathds{D} - \sigma_0 \mathds{H}^{-1} E_0$. This new operator is devoid of zero modes \cite{ruggiu_eigenvalue_2020} and may be inverted to obtain the solution ${\bf u}$.

In the context of an action functional, such as \cref{eq:discEIVP}, we do not have an equal sign around which we can move the SAT term. Instead we must incorporate the whole of the penalty term directly in a modified SBP operator. Since the penalty term in our example \cref{eq:simpleODE} contains both a contribution that is proportional to the function ${\bf u}$ and a constant shift ${\bf g}$ it amounts to an affine transformation on ${\bf u}$, which can be captured efficiently using affine coordinates. To this end let us write $\bar A[{\bf b}] \bar{\bf x} = A{\bf x}+{\bf b}$, where $\bar A[{\bf b}]$ refers to a matrix $A$ extended by an additional  row and column with the value $1$ placed in the lower right corner. The new column available in $\bar A[{\bf b}]$ is populated with the entries of ${\bf b}$. The vector $\bar{\bf x}$ is nothing but ${\bf x}$ extended by one more entry with value unity. We will use this construction principle to define a regularized $\bar{\mathds{D}}$ from our conventional SBP operator $\mathds{D}$.

 Since we have both ${\bf x}$ and ${\bf t}$ as independent degrees of freedom each with independent initial conditions $x_i$ and $t_i$, we must define different shifts ${\bf b}^x$ and ${\bf b}^t$ respectively and thus end up with two different regularized SBP operators ${\bar{\mathds{D}}}_t$ and ${\bar{\mathds{D}}}_x$. The shift terms are nothing but the constant part of the corresponding SAT term, absorbed into the SBP operator 
\begin{align}
    {\bf b}^x= \sigma_0 \mathds{H}^{-1} E_0 {\bf g}^x, \quad {\bf b}^t= \sigma_0 \mathds{H}^{-1} E_0 {\bf g}^t.
\end{align}
Here ${\bf g}^x={\rm diag}[x_i,0,\cdots,0]$ and ${\bf g}^t={\rm diag}[t_i,0,\cdots,0]$ encode the initial values for $x$ and $t$ respectively. As mentioned before, we choose the parameter $\sigma_0=-1$, whenever a penalty term is incorporated in ${\bar {\mathds{D}}}$, motivated by the fact that in the conventional treatment of IVPs using the SBP-SAT approach, this value leads to a minimal discretization error (see e.g. ref.~\cite{aalund2016provably,ALUND2019209,glaubitz2023summation}). The resulting regularized SBP operators to be deployed on ${\bf t}_{1,2}$ or ${\bf x}_{1,2}$, are given explicitly in \cref{eq:SBPregt} and \cref{eq:SBPregx} respectively.

Consistent with the affine coordinates used in the newly defined ${\bar{\mathds{D}}}_t$ and ${\bar{\mathds{D}}}_x$, we also amend the discretized trajectories ${\bf t}_{1,2}$ and ${\bf x}_{1,2}$ by one more entry that is given the value one. 

In order to compute inner products in the space of discretized functions, we also have to modify the quadrature matrix $\mathds{H}\to\bar{\mathds{H}}$ by amending it by one row and column filled with zeros. We do not include the value one in the lower right corner in order to correctly account for the fact that the vectors appearing as arguments to the inner product contain an auxiliary final entry, which does not contribute to the value of the inner product and only facilitates the efficient implementation of shift operations. For more details on the affine coordinate regularization technique see \cite{Rothkopf:2022zfb}.

\FloatBarrier

\begin{backmatter}

\section*{Competing interests}
  The authors declare that they have no competing interests.

\section*{Author's contributions}
    \begin{itemize}
         \item A. Rothkopf: formulation of the geometric variational approach, literature review, numerical experiments, writing, editing
         \item J. Nordst\"om: guidance on the formulation and implementation of SBP based discretization schemes, literature review, editing
    \end{itemize}


\bibliographystyle{stavanger-mathphys}


\bibliography{references}


\begin{thebibliography}{45}
\ifx \bisbn   \undefined \def \bisbn  #1{ISBN #1}\fi
\ifx \binits  \undefined \def \binits#1{#1}\fi
\ifx \bauthor  \undefined \def \bauthor#1{#1}\fi
\ifx \batitle  \undefined \def \batitle#1{#1}\fi
\ifx \bjtitle  \undefined \def \bjtitle#1{#1}\fi
\ifx \bvolume  \undefined \def \bvolume#1{\textbf{#1}}\fi
\ifx \byear  \undefined \def \byear#1{#1}\fi
\ifx \bissue  \undefined \def \bissue#1{#1}\fi
\ifx \bfpage  \undefined \def \bfpage#1{#1}\fi
\ifx \blpage  \undefined \def \blpage #1{#1}\fi
\ifx \burl  \undefined \def \burl#1{\textsf{#1}}\fi
\ifx \doiurl  \undefined \def \doiurl#1{\textsf{#1}}\fi
\ifx \betal  \undefined \def \betal{\textit{et al.}}\fi
\ifx \binstitute  \undefined \def \binstitute#1{#1}\fi
\ifx \binstitutionaled  \undefined \def \binstitutionaled#1{#1}\fi
\ifx \bctitle  \undefined \def \bctitle#1{#1}\fi
\ifx \beditor  \undefined \def \beditor#1{#1}\fi
\ifx \bpublisher  \undefined \def \bpublisher#1{#1}\fi
\ifx \bbtitle  \undefined \def \bbtitle#1{#1}\fi
\ifx \bedition  \undefined \def \bedition#1{#1}\fi
\ifx \bseriesno  \undefined \def \bseriesno#1{#1}\fi
\ifx \blocation  \undefined \def \blocation#1{#1}\fi
\ifx \bsertitle  \undefined \def \bsertitle#1{#1}\fi
\ifx \bsnm \undefined \def \bsnm#1{#1}\fi
\ifx \bsuffix \undefined \def \bsuffix#1{#1}\fi
\ifx \bparticle \undefined \def \bparticle#1{#1}\fi
\ifx \barticle \undefined \def \barticle#1{#1}\fi
\ifx \bconfdate \undefined \def \bconfdate #1{#1}\fi
\ifx \botherref \undefined \def \botherref #1{#1}\fi
\ifx \url \undefined \def \url#1{\textsf{#1}}\fi
\ifx \bchapter \undefined \def \bchapter#1{#1}\fi
\ifx \bbook \undefined \def \bbook#1{#1}\fi
\ifx \bcomment \undefined \def \bcomment#1{#1}\fi
\ifx \oauthor \undefined \def \oauthor#1{#1}\fi
\ifx \citeauthoryear \undefined \def \citeauthoryear#1{#1}\fi
\ifx \endbibitem  \undefined \def \endbibitem {}\fi
\ifx \bconflocation  \undefined \def \bconflocation#1{#1}\fi
\ifx \arxivurl  \undefined \def \arxivurl#1{\textsf{#1}}\fi
\csname PreBibitemsHook\endcsname

\bibitem{goldstein1980classical}
\begin{bbook}
\bauthor{\bsnm{Goldstein}, \binits{H.}},
\bauthor{\bsnm{Poole}, \binits{C.P.}},
\bauthor{\bsnm{Safko}, \binits{J.L.}}:
\bbtitle{Classical Mechanics}.
\bpublisher{Addison Wesley}
(\byear{2002})
\end{bbook}
\endbibitem

\bibitem{arnold2013mathematical}
\begin{bbook}
\bauthor{\bsnm{Arnold}, \binits{V.I.}},
\bauthor{\bsnm{Vogtmann}, \binits{K.}},
\bauthor{\bsnm{Weinstein}, \binits{A.}}:
\bbtitle{Mathematical Methods of Classical Mechanics}.
\bsertitle{Graduate Texts in Mathematics}.
\bpublisher{Springer}
(\byear{2013})
\end{bbook}
\endbibitem

\bibitem{Coleman:1985rnk}
\begin{bbook}
\bauthor{\bsnm{Coleman}, \binits{S.}}:
\bbtitle{{Aspects of Symmetry}: {Selected Erice Lectures}}.
\bpublisher{Cambridge University Press},
\blocation{Cambridge, U.K.}
(\byear{1985}).
doi:\doiurl{10.1017/CBO9780511565045}
\end{bbook}
\endbibitem

\bibitem{noether1971invariant}
\begin{barticle}
\bauthor{\bsnm{Noether}, \binits{E.}}:
\batitle{Invariant variation problems}.
\bjtitle{Transport theory and statistical physics}
\bvolume{1}(\bissue{3}),
\bfpage{186}--\blpage{207}
(\byear{1971})
\end{barticle}
\endbibitem

\bibitem{landau2000classical}
\begin{bbook}
\bauthor{\bsnm{Landau}, \binits{L.D.}},
\bauthor{\bsnm{Lifshitz}, \binits{E.M.}}:
\bbtitle{The Classical Theory of Fields: Volume 2}.
\bsertitle{Course of theoretical physics}.
\bpublisher{Elsevier Science}
(\byear{2000})
\end{bbook}
\endbibitem

\bibitem{Yanagihara:2018qqg}
\begin{barticle}
\bauthor{\bsnm{Yanagihara}, \binits{R.}},
\bauthor{\bsnm{Iritani}, \binits{T.}},
\bauthor{\bsnm{Kitazawa}, \binits{M.}},
\bauthor{\bsnm{Asakawa}, \binits{M.}},
\bauthor{\bsnm{Hatsuda}, \binits{T.}}:
\batitle{{Distribution of Stress Tensor around Static Quark--Anti-Quark from
  Yang-Mills Gradient Flow}}.
\bjtitle{Phys. Lett. B}
\bvolume{789},
\bfpage{210}--\blpage{214}
(\byear{2019}).
doi:\doiurl{10.1016/j.physletb.2018.09.067}.
\arxivurl{1803.05656}
\end{barticle}
\endbibitem

\bibitem{cockburn2012discontinuous}
\begin{bbook}
\bauthor{\bsnm{Cockburn}, \binits{B.}},
\bauthor{\bsnm{Karniadakis}, \binits{G.E.}},
\bauthor{\bsnm{Shu}, \binits{C.-W.}}:
\bbtitle{Discontinuous Galerkin Methods: Theory, Computation and Applications}
vol. \bseriesno{11}.
\bpublisher{Springer}
(\byear{2012})
\end{bbook}
\endbibitem

\bibitem{svard2014review}
\begin{barticle}
\bauthor{\bsnm{Sv{\"a}rd}, \binits{M.}},
\bauthor{\bsnm{Nordstr{\"o}m}, \binits{J.}}:
\batitle{Review of summation-by-parts schemes for initial--boundary-value
  problems}.
\bjtitle{Journal of Computational Physics}
\bvolume{268},
\bfpage{17}--\blpage{38}
(\byear{2014})
\end{barticle}
\endbibitem

\bibitem{fernandez2014review}
\begin{barticle}
\bauthor{\bsnm{Fern{\'a}ndez}, \binits{D.C.D.R.}},
\bauthor{\bsnm{Hicken}, \binits{J.E.}},
\bauthor{\bsnm{Zingg}, \binits{D.W.}}:
\batitle{Review of summation-by-parts operators with simultaneous approximation
  terms for the numerical solution of partial differential equations}.
\bjtitle{Computers \& Fluids}
\bvolume{95},
\bfpage{171}--\blpage{196}
(\byear{2014})
\end{barticle}
\endbibitem

\bibitem{lundquist2014sbp}
\begin{barticle}
\bauthor{\bsnm{Lundquist}, \binits{T.}},
\bauthor{\bsnm{Nordstr{\"o}m}, \binits{J.}}:
\batitle{The {SBP-SAT} technique for initial value problems}.
\bjtitle{Journal of Computational Physics}
\bvolume{270},
\bfpage{86}--\blpage{104}
(\byear{2014})
\end{barticle}
\endbibitem

\bibitem{nordstrom2013summation}
\begin{barticle}
\bauthor{\bsnm{Nordstr{\"o}m}, \binits{J.}},
\bauthor{\bsnm{Lundquist}, \binits{T.}}:
\batitle{Summation-by-parts in time}.
\bjtitle{Journal of Computational Physics}
\bvolume{251},
\bfpage{487}--\blpage{499}
(\byear{2013})
\end{barticle}
\endbibitem

\bibitem{nordstrom2016summation}
\begin{barticle}
\bauthor{\bsnm{Nordstr{\"o}m}, \binits{J.}},
\bauthor{\bsnm{Lundquist}, \binits{T.}}:
\batitle{Summation-by-parts in time: the second derivative}.
\bjtitle{SIAM Journal on Scientific Computing}
\bvolume{38}(\bissue{3}),
\bfpage{1561}--\blpage{1586}
(\byear{2016})
\end{barticle}
\endbibitem

\bibitem{JOHNSON1982147}
\begin{barticle}
\bauthor{\bsnm{Johnson}, \binits{R.C.}}:
\batitle{Angular momentum on a lattice}.
\bjtitle{Physics Letters B}
\bvolume{114}(\bissue{2}),
\bfpage{147}--\blpage{151}
(\byear{1982}).
doi:\doiurl{10.1016/0370-2693(82)90134-4}
\end{barticle}
\endbibitem

\bibitem{e6656c89-b0e4-3df4-8651-29dbb3d55273}
\begin{barticle}
\bauthor{\bsnm{Regan}, \binits{H.M.}}:
\batitle{Von neumann stability analysis of symplectic integrators applied to
  hamiltonian pdes}.
\bjtitle{Journal of Computational Mathematics}
\bvolume{20}(\bissue{6}),
\bfpage{611}--\blpage{618}
(\byear{2002})
\end{barticle}
\endbibitem

\bibitem{nordstrom2023nonlinear}
\begin{botherref}
\oauthor{\bsnm{Nordström}, \binits{J.}}:
Nonlinear Boundary Conditions for Initial Boundary Value Problems with
  Applications in Computational Fluid Dynamics.
\arxivurl{2306.01297}
\end{botherref}
\endbibitem

\bibitem{verlet_computer_1967}
\begin{barticle}
\bauthor{\bsnm{Verlet}, \binits{L.}}:
\batitle{Computer "{Experiments}" on {Classical} {Fluids}. {I}.
  {Thermodynamical} {Properties} of {Lennard}-{Jones} {Molecules}}.
\bjtitle{Physical Review}
\bvolume{159}(\bissue{1}),
\bfpage{98}--\blpage{103}
(\byear{1967}).
doi:\doiurl{10.1103/PhysRev.159.98}.
\bcomment{Publisher: American Physical Society}
\end{barticle}
\endbibitem

\bibitem{dirac_generalized_1950}
\begin{barticle}
\bauthor{\bsnm{Dirac}, \binits{P.a.M.}}:
\batitle{Generalized {Hamiltonian} {Dynamics}}.
\bjtitle{Canadian Journal of Mathematics}
\bvolume{2},
\bfpage{129}--\blpage{148}
(\byear{1950}).
doi:\doiurl{10.4153/CJM-1950-012-1}.
\bcomment{Publisher: Cambridge University Press}
\end{barticle}
\endbibitem

\bibitem{anerot_noethers-type_2020}
\begin{barticle}
\bauthor{\bsnm{Anerot}, \binits{B.}},
\bauthor{\bsnm{Cresson}, \binits{J.}},
\bauthor{\bsnm{Hariz~Belgacem}, \binits{K.}},
\bauthor{\bsnm{Pierret}, \binits{F.}}:
\batitle{Noether’s-type theorems on time scales}.
\bjtitle{Journal of Mathematical Physics}
\bvolume{61}(\bissue{11}),
\bfpage{113502}
(\byear{2020}).
doi:\doiurl{10.1063/1.5140201}.
\bcomment{Number: 11}
\end{barticle}
\endbibitem

\bibitem{stephani2004relativity}
\begin{bbook}
\bauthor{\bsnm{Stephani}, \binits{H.}}:
\bbtitle{Relativity: An Introduction to Special and General Relativity}.
\bpublisher{Cambridge University Press}
(\byear{2004})
\end{bbook}
\endbibitem

\bibitem{Rothkopf:2022zfb}
\begin{barticle}
\bauthor{\bsnm{Rothkopf}, \binits{A.}},
\bauthor{\bsnm{Nordstr\"om}, \binits{J.}}:
\batitle{{A new variational discretization technique for initial value problems
  bypassing governing equations}}.
\bjtitle{J. Comput. Phys.}
\bvolume{477},
\bfpage{111942}
(\byear{2023}).
doi:\doiurl{10.1016/j.jcp.2023.111942}.
\arxivurl{2205.14028}
\end{barticle}
\endbibitem

\bibitem{berger1984adaptive}
\begin{barticle}
\bauthor{\bsnm{Berger}, \binits{M.J.}},
\bauthor{\bsnm{Oliger}, \binits{J.}}:
\batitle{Adaptive mesh refinement for hyperbolic partial differential
  equations}.
\bjtitle{Journal of computational Physics}
\bvolume{53}(\bissue{3}),
\bfpage{484}--\blpage{512}
(\byear{1984})
\end{barticle}
\endbibitem

\bibitem{lohner1987adaptive}
\begin{barticle}
\bauthor{\bsnm{L{\"o}hner}, \binits{R.}}:
\batitle{An adaptive finite element scheme for transient problems in {CFD}}.
\bjtitle{Computer methods in applied mechanics and engineering}
\bvolume{61}(\bissue{3}),
\bfpage{323}--\blpage{338}
(\byear{1987})
\end{barticle}
\endbibitem

\bibitem{berger1989local}
\begin{barticle}
\bauthor{\bsnm{Berger}, \binits{M.J.}},
\bauthor{\bsnm{Colella}, \binits{P.}}:
\batitle{Local adaptive mesh refinement for shock hydrodynamics}.
\bjtitle{Journal of computational Physics}
\bvolume{82}(\bissue{1}),
\bfpage{64}--\blpage{84}
(\byear{1989})
\end{barticle}
\endbibitem

\bibitem{persson2006sub}
\begin{bchapter}
\bauthor{\bsnm{Persson}, \binits{P.-O.}},
\bauthor{\bsnm{Peraire}, \binits{J.}}:
\bctitle{Sub-cell shock capturing for discontinuous {Galerkin} methods}.
In: \bbtitle{44th AIAA Aerospace Sciences Meeting and Exhibit},
p. \bfpage{112}
(\byear{2006})
\end{bchapter}
\endbibitem

\bibitem{nemec2008adjoint}
\begin{bchapter}
\bauthor{\bsnm{Nemec}, \binits{M.}},
\bauthor{\bsnm{Aftosmis}, \binits{M.}},
\bauthor{\bsnm{Wintzer}, \binits{M.}}:
\bctitle{Adjoint-based adaptive mesh refinement for complex geometries}.
In: \bbtitle{46th AIAA Aerospace Sciences Meeting and Exhibit},
p. \bfpage{725}
(\byear{2008})
\end{bchapter}
\endbibitem

\bibitem{offermans2023error}
\begin{barticle}
\bauthor{\bsnm{Offermans}, \binits{N.}},
\bauthor{\bsnm{Massaro}, \binits{D.}},
\bauthor{\bsnm{Peplinski}, \binits{A.}},
\bauthor{\bsnm{Schlatter}, \binits{P.}}:
\batitle{Error-driven adaptive mesh refinement for unsteady turbulent flows in
  spectral-element simulations}.
\bjtitle{Computers \& Fluids}
\bvolume{251},
\bfpage{105736}
(\byear{2023})
\end{barticle}
\endbibitem

\bibitem{mavriplis1994adaptive}
\begin{barticle}
\bauthor{\bsnm{Mavriplis}, \binits{C.}}:
\batitle{Adaptive mesh strategies for the spectral element method}.
\bjtitle{Computer methods in applied mechanics and engineering}
\bvolume{116}(\bissue{1-4}),
\bfpage{77}--\blpage{86}
(\byear{1994})
\end{barticle}
\endbibitem

\bibitem{henderson1999adaptive}
\begin{bchapter}
\bauthor{\bsnm{Henderson}, \binits{R.D.}}:
\bctitle{Adaptive spectral element methods for turbulence and transition}.
In: \bbtitle{High-order Methods for Computational Physics},
pp. \bfpage{225}--\blpage{324}.
\bpublisher{Springer}
(\byear{1999})
\end{bchapter}
\endbibitem

\bibitem{kompenhans2016adaptation}
\begin{barticle}
\bauthor{\bsnm{Kompenhans}, \binits{M.}},
\bauthor{\bsnm{Rubio}, \binits{G.}},
\bauthor{\bsnm{Ferrer}, \binits{E.}},
\bauthor{\bsnm{Valero}, \binits{E.}}:
\batitle{Adaptation strategies for high order discontinuous {Galerkin} methods
  based on tau-estimation}.
\bjtitle{Journal of Computational Physics}
\bvolume{306},
\bfpage{216}--\blpage{236}
(\byear{2016})
\end{barticle}
\endbibitem

\bibitem{galley_classical_2013}
\begin{barticle}
\bauthor{\bsnm{Galley}, \binits{C.R.}}:
\batitle{Classical {Mechanics} of {Nonconservative} {Systems}}.
\bjtitle{Physical Review Letters}
\bvolume{110}(\bissue{17}),
\bfpage{174301}
(\byear{2013}).
doi:\doiurl{10.1103/PhysRevLett.110.174301}.
\bcomment{Publisher: American Physical Society}
\end{barticle}
\endbibitem

\bibitem{jost1998calculus}
\begin{bbook}
\bauthor{\bsnm{Jost}, \binits{J.}},
\bauthor{\bsnm{Li-Jost}, \binits{X.}}:
\bbtitle{Calculus of Variations}.
\bsertitle{Cambridge Studies in Advanced Mathematics}.
\bpublisher{Cambridge University Press}
(\byear{1998})
\end{bbook}
\endbibitem

\bibitem{carroll2019spacetime}
\begin{bbook}
\bauthor{\bsnm{Carroll}, \binits{S.M.}}:
\bbtitle{Spacetime and Geometry}.
\bpublisher{Cambridge University Press}
(\byear{2019})
\end{bbook}
\endbibitem

\bibitem{carlip2019general}
\begin{bbook}
\bauthor{\bsnm{Carlip}, \binits{S.}}:
\bbtitle{General Relativity: A Concise Introduction}.
\bpublisher{OUP Oxford}
(\byear{2019})
\end{bbook}
\endbibitem

\bibitem{rizzuti_square_2019}
\begin{botherref}
\oauthor{\bsnm{Rizzuti}, \binits{B.F.}},
\oauthor{\bsnm{Júnior}, \binits{G.F.V.}},
\oauthor{\bsnm{Resende}, \binits{M.A.}}:
To square root the {Lagrangian} or not: an underlying geometrical analysis on
  classical and relativistic mechanical models.
arXiv.
arXiv:1905.01177 [math-ph, physics:physics]
(2019).
doi:\doiurl{10.48550/arXiv.1905.01177}.
\url{http://arxiv.org/abs/1905.01177}
\end{botherref}
\endbibitem

\bibitem{Berges:2007ym}
\begin{barticle}
\bauthor{\bsnm{Berges}, \binits{J.}},
\bauthor{\bsnm{Gasenzer}, \binits{T.}}:
\batitle{{Quantum versus classical statistical dynamics of an ultracold Bose
  gas}}.
\bjtitle{Phys. Rev. A}
\bvolume{76},
\bfpage{033604}
(\byear{2007}).
doi:\doiurl{10.1103/PhysRevA.76.033604}.
\arxivurl{cond-mat/0703163}
\end{barticle}
\endbibitem

\bibitem{carpenter1994time}
\begin{barticle}
\bauthor{\bsnm{Carpenter}, \binits{M.H.}},
\bauthor{\bsnm{Gottlieb}, \binits{D.}},
\bauthor{\bsnm{Abarbanel}, \binits{S.}}:
\batitle{Time-stable boundary conditions for finite-difference schemes solving
  hyperbolic systems: Methodology and application to high-order compact
  schemes}.
\bjtitle{Journal of Computational Physics}
\bvolume{111}(\bissue{2}),
\bfpage{220}--\blpage{236}
(\byear{1994})
\end{barticle}
\endbibitem

\bibitem{zenodoIVP2023}
\begin{botherref}
\oauthor{\bsnm{Rothkopf}, \binits{A.}}:
Mathematica 12 implementation of a symmetry and {Noether} charge preserving
  {IVP} discretization technique
(2023).
doi:\doiurl{10.5281/zenodo.8129657}
\end{botherref}
\endbibitem

\bibitem{svard2019convergence}
\begin{barticle}
\bauthor{\bsnm{Sv{\"a}rd}, \binits{M.}},
\bauthor{\bsnm{Nordstr{\"o}m}, \binits{J.}}:
\batitle{On the convergence rates of energy-stable finite-difference schemes}.
\bjtitle{Journal of Computational Physics}
\bvolume{397},
\bfpage{108819}
(\byear{2019})
\end{barticle}
\endbibitem

\bibitem{linders_properties_2020}
\begin{barticle}
\bauthor{\bsnm{Linders}, \binits{V.}},
\bauthor{\bsnm{Nordstr{\"o}m}, \binits{J.}},
\bauthor{\bsnm{Frankel}, \binits{S.H.}}:
\batitle{Properties of {Runge}-{Kutta}-{Summation}-{By}-{Parts} methods}.
\bjtitle{Journal of Computational Physics}
\bvolume{419},
\bfpage{109684}
(\byear{2020}).
doi:\doiurl{10.1016/j.jcp.2020.109684}
\end{barticle}
\endbibitem

\bibitem{svard_convergence_2021}
\begin{barticle}
\bauthor{\bsnm{Sv{\"a}rd}, \binits{M.}},
\bauthor{\bsnm{Nordstr{\"o}m}, \binits{J.}}:
\batitle{Convergence of energy stable finite-difference schemes with
  interfaces}.
\bjtitle{Journal of Computational Physics}
\bvolume{429},
\bfpage{110020}
(\byear{2021}).
doi:\doiurl{10.1016/j.jcp.2020.110020}
\end{barticle}
\endbibitem

\bibitem{ranocha2021new}
\begin{barticle}
\bauthor{\bsnm{Ranocha}, \binits{H.}},
\bauthor{\bsnm{Nordstr{\"o}m}, \binits{J.}}:
\batitle{A new class of a stable summation by parts time integration schemes
  with strong initial conditions}.
\bjtitle{Journal of Scientific Computing}
\bvolume{87}(\bissue{1}),
\bfpage{1}--\blpage{25}
(\byear{2021})
\end{barticle}
\endbibitem

\bibitem{aalund2016provably}
\begin{botherref}
\oauthor{\bsnm{{\AA}lund}, \binits{O.}},
\oauthor{\bsnm{Nordstr{\"o}m}, \binits{J.}}:
A provably stable, non-iterative domain decomposition technique for the
  advection-diffusion equation.
Link{\"o}ping University Electronic Press
(2016)
\end{botherref}
\endbibitem

\bibitem{ALUND2019209}
\begin{barticle}
\bauthor{\bsnm{Ålund}, \binits{O.}},
\bauthor{\bsnm{Nordström}, \binits{J.}}:
\batitle{Encapsulated high order difference operators on curvilinear
  non-conforming grids}.
\bjtitle{Journal of Computational Physics}
\bvolume{385},
\bfpage{209}--\blpage{224}
(\byear{2019}).
doi:\doiurl{10.1016/j.jcp.2019.02.007}
\end{barticle}
\endbibitem

\bibitem{glaubitz2023summation}
\begin{barticle}
\bauthor{\bsnm{Glaubitz}, \binits{J.}},
\bauthor{\bsnm{Nordstr{\"o}m}, \binits{J.}},
\bauthor{\bsnm{{\"O}ffner}, \binits{P.}}:
\batitle{Summation-by-parts operators for general function spaces}.
\bjtitle{SIAM Journal on Numerical Analysis}
\bvolume{61}(\bissue{2}),
\bfpage{733}--\blpage{754}
(\byear{2023})
\end{barticle}
\endbibitem

\bibitem{ruggiu_eigenvalue_2020}
\begin{barticle}
\bauthor{\bsnm{Ruggiu}, \binits{A.A.}},
\bauthor{\bsnm{Nordstr{\"o}m}, \binits{J.}}:
\batitle{Eigenvalue {Analysis} for {Summation}-by-{Parts} {Finite} {Difference}
  {Time} {Discretizations}}.
\bjtitle{SIAM Journal on Numerical Analysis}
\bvolume{58}(\bissue{2}),
\bfpage{907}--\blpage{928}
(\byear{2020}).
doi:\doiurl{10.1137/19M1256294}.
\bcomment{Publisher: Society for Industrial and Applied Mathematics}
\end{barticle}
\endbibitem

\end{thebibliography}

\newcommand{\BMCxmlcomment}[1]{}

\BMCxmlcomment{

<refgrp>

<bibl id="B1">
  <title><p>Classical Mechanics</p></title>
  <aug>
    <au><snm>Goldstein</snm><fnm>H.</fnm></au>
    <au><snm>Poole</snm><fnm>C.P.</fnm></au>
    <au><snm>Safko</snm><fnm>J.L.</fnm></au>
  </aug>
  <publisher>Addison Wesley</publisher>
  <pubdate>2002</pubdate>
</bibl>

<bibl id="B2">
  <title><p>Mathematical Methods of Classical Mechanics</p></title>
  <aug>
    <au><snm>Arnold</snm><fnm>V.I.</fnm></au>
    <au><snm>Vogtmann</snm><fnm>K.</fnm></au>
    <au><snm>Weinstein</snm><fnm>A.</fnm></au>
  </aug>
  <publisher>Springer New York</publisher>
  <series><title><p>Graduate Texts in Mathematics</p></title></series>
  <pubdate>2013</pubdate>
</bibl>

<bibl id="B3">
  <title><p>{Aspects of Symmetry}: {Selected Erice Lectures}</p></title>
  <aug>
    <au><snm>Coleman</snm><fnm>S</fnm></au>
  </aug>
  <publisher>Cambridge, U.K.: Cambridge University Press</publisher>
  <pubdate>1985</pubdate>
</bibl>

<bibl id="B4">
  <title><p>Invariant variation problems</p></title>
  <aug>
    <au><snm>Noether</snm><fnm>E</fnm></au>
  </aug>
  <source>Transport theory and statistical physics</source>
  <publisher>Taylor \& Francis</publisher>
  <pubdate>1971</pubdate>
  <volume>1</volume>
  <issue>3</issue>
  <fpage>186</fpage>
  <lpage>-207</lpage>
</bibl>

<bibl id="B5">
  <title><p>The Classical Theory of Fields: Volume 2</p></title>
  <aug>
    <au><snm>Landau</snm><fnm>L.D.</fnm></au>
    <au><snm>Lifshitz</snm><fnm>E.M.</fnm></au>
  </aug>
  <publisher>Elsevier Science</publisher>
  <series><title><p>Course of theoretical physics</p></title></series>
  <pubdate>2000</pubdate>
</bibl>

<bibl id="B6">
  <title><p>{Distribution of Stress Tensor around Static Quark--Anti-Quark from
  Yang-Mills Gradient Flow}</p></title>
  <aug>
    <au><snm>Yanagihara</snm><fnm>R</fnm></au>
    <au><snm>Iritani</snm><fnm>T</fnm></au>
    <au><snm>Kitazawa</snm><fnm>M</fnm></au>
    <au><snm>Asakawa</snm><fnm>M</fnm></au>
    <au><snm>Hatsuda</snm><fnm>T</fnm></au>
  </aug>
  <source>Phys. Lett. B</source>
  <pubdate>2019</pubdate>
  <volume>789</volume>
  <fpage>210</fpage>
  <lpage>-214</lpage>
</bibl>

<bibl id="B7">
  <title><p>Discontinuous Galerkin methods: theory, computation and
  applications</p></title>
  <aug>
    <au><snm>Cockburn</snm><fnm>B</fnm></au>
    <au><snm>Karniadakis</snm><fnm>GE</fnm></au>
    <au><snm>Shu</snm><fnm>CW</fnm></au>
  </aug>
  <publisher>Springer Science \& Business Media</publisher>
  <pubdate>2012</pubdate>
  <volume>11</volume>
</bibl>

<bibl id="B8">
  <title><p>Review of summation-by-parts schemes for initial--boundary-value
  problems</p></title>
  <aug>
    <au><snm>Sv{\"a}rd</snm><fnm>M</fnm></au>
    <au><snm>Nordstr{\"o}m</snm><fnm>J</fnm></au>
  </aug>
  <source>Journal of Computational Physics</source>
  <publisher>Elsevier</publisher>
  <pubdate>2014</pubdate>
  <volume>268</volume>
  <fpage>17</fpage>
  <lpage>-38</lpage>
</bibl>

<bibl id="B9">
  <title><p>Review of summation-by-parts operators with simultaneous
  approximation terms for the numerical solution of partial differential
  equations</p></title>
  <aug>
    <au><snm>Fern{\'a}ndez</snm><fnm>DCDR</fnm></au>
    <au><snm>Hicken</snm><fnm>JE</fnm></au>
    <au><snm>Zingg</snm><fnm>DW</fnm></au>
  </aug>
  <source>Computers \& Fluids</source>
  <publisher>Elsevier</publisher>
  <pubdate>2014</pubdate>
  <volume>95</volume>
  <fpage>171</fpage>
  <lpage>-196</lpage>
</bibl>

<bibl id="B10">
  <title><p>The {SBP-SAT} technique for initial value problems</p></title>
  <aug>
    <au><snm>Lundquist</snm><fnm>T</fnm></au>
    <au><snm>Nordstr{\"o}m</snm><fnm>J</fnm></au>
  </aug>
  <source>Journal of Computational Physics</source>
  <publisher>Elsevier</publisher>
  <pubdate>2014</pubdate>
  <volume>270</volume>
  <fpage>86</fpage>
  <lpage>-104</lpage>
</bibl>

<bibl id="B11">
  <title><p>Summation-by-parts in time</p></title>
  <aug>
    <au><snm>Nordstr{\"o}m</snm><fnm>J</fnm></au>
    <au><snm>Lundquist</snm><fnm>T</fnm></au>
  </aug>
  <source>Journal of Computational Physics</source>
  <publisher>Elsevier</publisher>
  <pubdate>2013</pubdate>
  <volume>251</volume>
  <fpage>487</fpage>
  <lpage>-499</lpage>
</bibl>

<bibl id="B12">
  <title><p>Summation-by-parts in time: the second derivative</p></title>
  <aug>
    <au><snm>Nordstr{\"o}m</snm><fnm>J</fnm></au>
    <au><snm>Lundquist</snm><fnm>T</fnm></au>
  </aug>
  <source>SIAM Journal on Scientific Computing</source>
  <publisher>SIAM</publisher>
  <pubdate>2016</pubdate>
  <volume>38</volume>
  <issue>3</issue>
  <fpage>A1561</fpage>
  <lpage>-A1586</lpage>
</bibl>

<bibl id="B13">
  <title><p>Angular momentum on a lattice</p></title>
  <aug>
    <au><snm>Johnson</snm><fnm>R.C.</fnm></au>
  </aug>
  <source>Physics Letters B</source>
  <pubdate>1982</pubdate>
  <volume>114</volume>
  <issue>2</issue>
  <fpage>147</fpage>
  <lpage>151</lpage>
  <url>https://www.sciencedirect.com/science/article/pii/0370269382901344</url>
</bibl>

<bibl id="B14">
  <title><p>VON NEUMANN STABILITY ANALYSIS OF SYMPLECTIC INTEGRATORS APPLIED TO
  HAMILTONIAN PDEs</p></title>
  <aug>
    <au><snm>Regan</snm><fnm>HM</fnm></au>
  </aug>
  <source>Journal of Computational Mathematics</source>
  <publisher>Institute of Computational Mathematics and Scientific/Engineering
  Computing</publisher>
  <pubdate>2002</pubdate>
  <volume>20</volume>
  <issue>6</issue>
  <fpage>611</fpage>
  <lpage>-618</lpage>
  <url>http://www.jstor.org/stable/43693028</url>
</bibl>

<bibl id="B15">
  <title><p>Nonlinear Boundary Conditions for Initial Boundary Value Problems
  with Applications in Computational Fluid Dynamics</p></title>
  <aug>
    <au><snm>Nordström</snm><fnm>J</fnm></au>
  </aug>
  <pubdate>2023</pubdate>
</bibl>

<bibl id="B16">
  <title><p>Computer "{Experiments}" on {Classical} {Fluids}. {I}.
  {Thermodynamical} {Properties} of {Lennard}-{Jones} {Molecules}</p></title>
  <aug>
    <au><snm>Verlet</snm><fnm>L</fnm></au>
  </aug>
  <source>Physical Review</source>
  <pubdate>1967</pubdate>
  <volume>159</volume>
  <issue>1</issue>
  <fpage>98</fpage>
  <lpage>-103</lpage>
  <url>https://link.aps.org/doi/10.1103/PhysRev.159.98</url>
  <note>Publisher: American Physical Society</note>
</bibl>

<bibl id="B17">
  <title><p>Generalized {Hamiltonian} {Dynamics}</p></title>
  <aug>
    <au><snm>Dirac</snm><fnm>PaM</fnm></au>
  </aug>
  <source>Canadian Journal of Mathematics</source>
  <pubdate>1950</pubdate>
  <volume>2</volume>
  <fpage>129</fpage>
  <lpage>-148</lpage>
  <url>https://www.cambridge.org/core/journals/canadian-journal-of-mathematics/article/generalized-hamiltonian-dynamics/F4C30A59B59BEE09E9CB6F07377B8BD3</url>
  <note>Publisher: Cambridge University Press</note>
</bibl>

<bibl id="B18">
  <title><p>Noether’s-type theorems on time scales</p></title>
  <aug>
    <au><snm>Anerot</snm><fnm>B</fnm></au>
    <au><snm>Cresson</snm><fnm>J</fnm></au>
    <au><snm>Hariz Belgacem</snm><fnm>K</fnm></au>
    <au><snm>Pierret</snm><fnm>F</fnm></au>
  </aug>
  <source>Journal of Mathematical Physics</source>
  <pubdate>2020</pubdate>
  <volume>61</volume>
  <issue>11</issue>
  <fpage>113502</fpage>
  <url>http://aip.scitation.org/doi/10.1063/1.5140201</url>
  <note>Number: 11</note>
</bibl>

<bibl id="B19">
  <title><p>Relativity: An Introduction to Special and General
  Relativity</p></title>
  <aug>
    <au><snm>Stephani</snm><fnm>H.</fnm></au>
  </aug>
  <publisher>Cambridge University Press</publisher>
  <pubdate>2004</pubdate>
</bibl>

<bibl id="B20">
  <title><p>{A new variational discretization technique for initial value
  problems bypassing governing equations}</p></title>
  <aug>
    <au><snm>Rothkopf</snm><fnm>A</fnm></au>
    <au><snm>Nordstr\"om</snm><fnm>J</fnm></au>
  </aug>
  <source>J. Comput. Phys.</source>
  <pubdate>2023</pubdate>
  <volume>477</volume>
  <fpage>111942</fpage>
</bibl>

<bibl id="B21">
  <title><p>Adaptive mesh refinement for hyperbolic partial differential
  equations</p></title>
  <aug>
    <au><snm>Berger</snm><fnm>MJ</fnm></au>
    <au><snm>Oliger</snm><fnm>J</fnm></au>
  </aug>
  <source>Journal of computational Physics</source>
  <publisher>Elsevier</publisher>
  <pubdate>1984</pubdate>
  <volume>53</volume>
  <issue>3</issue>
  <fpage>484</fpage>
  <lpage>-512</lpage>
</bibl>

<bibl id="B22">
  <title><p>An adaptive finite element scheme for transient problems in
  {CFD}</p></title>
  <aug>
    <au><snm>L{\"o}hner</snm><fnm>R</fnm></au>
  </aug>
  <source>Computer methods in applied mechanics and engineering</source>
  <publisher>Elsevier</publisher>
  <pubdate>1987</pubdate>
  <volume>61</volume>
  <issue>3</issue>
  <fpage>323</fpage>
  <lpage>-338</lpage>
</bibl>

<bibl id="B23">
  <title><p>Local adaptive mesh refinement for shock hydrodynamics</p></title>
  <aug>
    <au><snm>Berger</snm><fnm>MJ</fnm></au>
    <au><snm>Colella</snm><fnm>P</fnm></au>
  </aug>
  <source>Journal of computational Physics</source>
  <publisher>Elsevier</publisher>
  <pubdate>1989</pubdate>
  <volume>82</volume>
  <issue>1</issue>
  <fpage>64</fpage>
  <lpage>-84</lpage>
</bibl>

<bibl id="B24">
  <title><p>Sub-cell shock capturing for discontinuous {Galerkin}
  methods</p></title>
  <aug>
    <au><snm>Persson</snm><fnm>PO</fnm></au>
    <au><snm>Peraire</snm><fnm>J</fnm></au>
  </aug>
  <source>44th AIAA aerospace sciences meeting and exhibit</source>
  <pubdate>2006</pubdate>
  <fpage>112</fpage>
</bibl>

<bibl id="B25">
  <title><p>Adjoint-based adaptive mesh refinement for complex
  geometries</p></title>
  <aug>
    <au><snm>Nemec</snm><fnm>M</fnm></au>
    <au><snm>Aftosmis</snm><fnm>M</fnm></au>
    <au><snm>Wintzer</snm><fnm>M</fnm></au>
  </aug>
  <source>46th AIAA Aerospace Sciences Meeting and Exhibit</source>
  <pubdate>2008</pubdate>
  <fpage>725</fpage>
</bibl>

<bibl id="B26">
  <title><p>Error-driven adaptive mesh refinement for unsteady turbulent flows
  in spectral-element simulations</p></title>
  <aug>
    <au><snm>Offermans</snm><fnm>N</fnm></au>
    <au><snm>Massaro</snm><fnm>D</fnm></au>
    <au><snm>Peplinski</snm><fnm>A</fnm></au>
    <au><snm>Schlatter</snm><fnm>P</fnm></au>
  </aug>
  <source>Computers \& Fluids</source>
  <publisher>Elsevier</publisher>
  <pubdate>2023</pubdate>
  <volume>251</volume>
  <fpage>105736</fpage>
</bibl>

<bibl id="B27">
  <title><p>Adaptive mesh strategies for the spectral element
  method</p></title>
  <aug>
    <au><snm>Mavriplis</snm><fnm>C</fnm></au>
  </aug>
  <source>Computer methods in applied mechanics and engineering</source>
  <publisher>Elsevier</publisher>
  <pubdate>1994</pubdate>
  <volume>116</volume>
  <issue>1-4</issue>
  <fpage>77</fpage>
  <lpage>-86</lpage>
</bibl>

<bibl id="B28">
  <title><p>Adaptive spectral element methods for turbulence and
  transition</p></title>
  <aug>
    <au><snm>Henderson</snm><fnm>RD</fnm></au>
  </aug>
  <source>High-order methods for computational physics</source>
  <publisher>Springer</publisher>
  <pubdate>1999</pubdate>
  <fpage>225</fpage>
  <lpage>-324</lpage>
</bibl>

<bibl id="B29">
  <title><p>Adaptation strategies for high order discontinuous {Galerkin}
  methods based on Tau-estimation</p></title>
  <aug>
    <au><snm>Kompenhans</snm><fnm>M</fnm></au>
    <au><snm>Rubio</snm><fnm>G</fnm></au>
    <au><snm>Ferrer</snm><fnm>E</fnm></au>
    <au><snm>Valero</snm><fnm>E</fnm></au>
  </aug>
  <source>Journal of Computational Physics</source>
  <publisher>Elsevier</publisher>
  <pubdate>2016</pubdate>
  <volume>306</volume>
  <fpage>216</fpage>
  <lpage>-236</lpage>
</bibl>

<bibl id="B30">
  <title><p>Classical {Mechanics} of {Nonconservative} {Systems}</p></title>
  <aug>
    <au><snm>Galley</snm><fnm>CR</fnm></au>
  </aug>
  <source>Physical Review Letters</source>
  <pubdate>2013</pubdate>
  <volume>110</volume>
  <issue>17</issue>
  <fpage>174301</fpage>
  <url>https://link.aps.org/doi/10.1103/PhysRevLett.110.174301</url>
  <note>Publisher: American Physical Society</note>
</bibl>

<bibl id="B31">
  <title><p>Calculus of Variations</p></title>
  <aug>
    <au><snm>Jost</snm><fnm>J.</fnm></au>
    <au><snm>Li Jost</snm><fnm>X.</fnm></au>
  </aug>
  <publisher>Cambridge University Press</publisher>
  <series><title><p>Cambridge Studies in Advanced
  Mathematics</p></title></series>
  <pubdate>1998</pubdate>
</bibl>

<bibl id="B32">
  <title><p>Spacetime and Geometry</p></title>
  <aug>
    <au><snm>Carroll</snm><fnm>S.M.</fnm></au>
  </aug>
  <publisher>Cambridge University Press</publisher>
  <pubdate>2019</pubdate>
</bibl>

<bibl id="B33">
  <title><p>General Relativity: A Concise Introduction</p></title>
  <aug>
    <au><snm>Carlip</snm><fnm>S.</fnm></au>
  </aug>
  <publisher>OUP Oxford</publisher>
  <pubdate>2019</pubdate>
</bibl>

<bibl id="B34">
  <title><p>To square root the {Lagrangian} or not: an underlying geometrical
  analysis on classical and relativistic mechanical models</p></title>
  <aug>
    <au><snm>Rizzuti</snm><fnm>B. F.</fnm></au>
    <au><snm>Júnior</snm><fnm>GFV</fnm></au>
    <au><snm>Resende</snm><fnm>M. A.</fnm></au>
  </aug>
  <publisher>arXiv</publisher>
  <pubdate>2019</pubdate>
  <url>http://arxiv.org/abs/1905.01177</url>
  <note>arXiv:1905.01177 [math-ph, physics:physics]</note>
</bibl>

<bibl id="B35">
  <title><p>{Quantum versus classical statistical dynamics of an ultracold Bose
  gas}</p></title>
  <aug>
    <au><snm>Berges</snm><fnm>J</fnm></au>
    <au><snm>Gasenzer</snm><fnm>T</fnm></au>
  </aug>
  <source>Phys. Rev. A</source>
  <pubdate>2007</pubdate>
  <volume>76</volume>
  <fpage>033604</fpage>
</bibl>

<bibl id="B36">
  <title><p>Time-stable boundary conditions for finite-difference schemes
  solving hyperbolic systems: Methodology and application to high-order compact
  schemes</p></title>
  <aug>
    <au><snm>Carpenter</snm><fnm>MH</fnm></au>
    <au><snm>Gottlieb</snm><fnm>D</fnm></au>
    <au><snm>Abarbanel</snm><fnm>S</fnm></au>
  </aug>
  <source>Journal of Computational Physics</source>
  <publisher>Elsevier</publisher>
  <pubdate>1994</pubdate>
  <volume>111</volume>
  <issue>2</issue>
  <fpage>220</fpage>
  <lpage>-236</lpage>
</bibl>

<bibl id="B37">
  <title><p>Mathematica 12 implementation of a symmetry and {Noether} charge
  preserving {IVP} discretization technique</p></title>
  <aug>
    <au><snm>Rothkopf</snm><fnm>A</fnm></au>
  </aug>
  <publisher>Zenodo</publisher>
  <pubdate>2023</pubdate>
</bibl>

<bibl id="B38">
  <title><p>On the convergence rates of energy-stable finite-difference
  schemes</p></title>
  <aug>
    <au><snm>Sv{\"a}rd</snm><fnm>M</fnm></au>
    <au><snm>Nordstr{\"o}m</snm><fnm>J</fnm></au>
  </aug>
  <source>Journal of Computational Physics</source>
  <publisher>Elsevier</publisher>
  <pubdate>2019</pubdate>
  <volume>397</volume>
  <fpage>108819</fpage>
</bibl>

<bibl id="B39">
  <title><p>Properties of {Runge}-{Kutta}-{Summation}-{By}-{Parts}
  methods</p></title>
  <aug>
    <au><snm>Linders</snm><fnm>V</fnm></au>
    <au><snm>Nordstr{\"o}m</snm><fnm>J</fnm></au>
    <au><snm>Frankel</snm><fnm>SH</fnm></au>
  </aug>
  <source>Journal of Computational Physics</source>
  <pubdate>2020</pubdate>
  <volume>419</volume>
  <fpage>109684</fpage>
  <url>https://www.sciencedirect.com/science/article/pii/S0021999120304587</url>
</bibl>

<bibl id="B40">
  <title><p>Convergence of energy stable finite-difference schemes with
  interfaces</p></title>
  <aug>
    <au><snm>Sv{\"a}rd</snm><fnm>M</fnm></au>
    <au><snm>Nordstr{\"o}m</snm><fnm>J</fnm></au>
  </aug>
  <source>Journal of Computational Physics</source>
  <pubdate>2021</pubdate>
  <volume>429</volume>
  <fpage>110020</fpage>
  <url>https://www.sciencedirect.com/science/article/pii/S0021999120307944</url>
</bibl>

<bibl id="B41">
  <title><p>A new class of A stable summation by parts time integration schemes
  with strong initial conditions</p></title>
  <aug>
    <au><snm>Ranocha</snm><fnm>H</fnm></au>
    <au><snm>Nordstr{\"o}m</snm><fnm>J</fnm></au>
  </aug>
  <source>Journal of Scientific Computing</source>
  <publisher>Springer</publisher>
  <pubdate>2021</pubdate>
  <volume>87</volume>
  <issue>1</issue>
  <fpage>1</fpage>
  <lpage>-25</lpage>
</bibl>

<bibl id="B42">
  <title><p>A provably stable, non-iterative domain decomposition technique for
  the advection-diffusion equation</p></title>
  <aug>
    <au><snm>{\AA}lund</snm><fnm>O</fnm></au>
    <au><snm>Nordstr{\"o}m</snm><fnm>J</fnm></au>
  </aug>
  <publisher>Link{\"o}ping University Electronic Press</publisher>
  <pubdate>2016</pubdate>
</bibl>

<bibl id="B43">
  <title><p>Encapsulated high order difference operators on curvilinear
  non-conforming grids</p></title>
  <aug>
    <au><snm>Ålund</snm><fnm>O</fnm></au>
    <au><snm>Nordström</snm><fnm>J</fnm></au>
  </aug>
  <source>Journal of Computational Physics</source>
  <pubdate>2019</pubdate>
  <volume>385</volume>
  <fpage>209</fpage>
  <lpage>224</lpage>
  <url>https://www.sciencedirect.com/science/article/pii/S0021999119301184</url>
</bibl>

<bibl id="B44">
  <title><p>Summation-by-parts operators for general function
  spaces</p></title>
  <aug>
    <au><snm>Glaubitz</snm><fnm>J</fnm></au>
    <au><snm>Nordstr{\"o}m</snm><fnm>J</fnm></au>
    <au><snm>{\"O}ffner</snm><fnm>P</fnm></au>
  </aug>
  <source>SIAM Journal on Numerical Analysis</source>
  <publisher>SIAM</publisher>
  <pubdate>2023</pubdate>
  <volume>61</volume>
  <issue>2</issue>
  <fpage>733</fpage>
  <lpage>-754</lpage>
</bibl>

<bibl id="B45">
  <title><p>Eigenvalue {Analysis} for {Summation}-by-{Parts} {Finite}
  {Difference} {Time} {Discretizations}</p></title>
  <aug>
    <au><snm>Ruggiu</snm><fnm>AA</fnm></au>
    <au><snm>Nordstr{\"o}m</snm><fnm>J</fnm></au>
  </aug>
  <source>SIAM Journal on Numerical Analysis</source>
  <pubdate>2020</pubdate>
  <volume>58</volume>
  <issue>2</issue>
  <fpage>907</fpage>
  <lpage>-928</lpage>
  <url>https://epubs.siam.org/doi/abs/10.1137/19M1256294</url>
  <note>Publisher: Society for Industrial and Applied Mathematics</note>
</bibl>

</refgrp>
} 


\end{backmatter}


\end{document}